\documentclass{amsart}
\usepackage{amscd}
\usepackage{amsmath}
\usepackage{latexsym}
\usepackage{graphicx}
\usepackage{amsfonts}
\usepackage{amssymb}%

\newcounter{fig}

\theoremstyle{plain}
\newtheorem{Theorem}{Theorem}[section]

\newtheorem{theorem}[Theorem]{Theorem}
\newtheorem{conjecture}[Theorem]{Conjecture}
\newtheorem{corollary}[Theorem]{Corollary}
\newtheorem{proposition}[Theorem]{Proposition}
\newtheorem{lemma}[Theorem]{Lemma}
\newtheorem{Alternative Version}{Alternative Version}
\theoremstyle{definition}

\newtheorem{definition}[Theorem]{Definition}
\newtheorem{Convention}[Theorem]{Convention}

\newtheorem{example}[Theorem]{Example}

\newtheorem{remark}[Theorem]{Remark}

\newtheorem{Remarks}[Theorem]{Remarks}
\newtheorem{Examples}[Theorem]{Examples}

\newtheorem{ex}[Theorem]{Example}

\begin{document}
\title{Periodicity of hermitian $K$-groups}
\author[BKO \qquad\today]{A.\thinspace J.~Berrick, M.~Karoubi and P.\thinspace A.~{\O }stv{\ae }r}
\thanks{First and second authors partially supported by the National University of
Singapore R-146-000-097-112. Third author partially supported by RCN 185335/V30.}
\date{December 10, 2010}
\maketitle

\pagestyle{myheadings} \markright{{Hermitian periodicity}  \quad\today}
\setcounter{section}{-1}%

\section{Introduction and statements of main results}

By the fundamental work of Bott \cite{Bott} it is known that the homotopy
groups of classical Lie groups are periodic, of period $2$ or $8$. For
instance, the general linear and symplectic groups satisfy the isomorphisms:%
\[
\pi_{n}(\mathrm{GL}(\mathbb{R)})\cong\pi_{n+8}(\mathrm{GL}(\mathbb{R)})
\]%
\[
\pi_{n}(\mathrm{Sp}(\mathbb{C)})\cong\pi_{n+8}(\mathrm{Sp}(\mathbb{C)})
\]%
\[
\pi_{n}(\mathrm{GL}(\mathbb{C)})\cong\pi_{n+2}(\mathrm{GL}(\mathbb{C)})
\]
These periodicity statements were interpreted by Atiyah, Hirzebruch and others
in the framework of topological $K$-theory of a Banach algebra $A$: recall
that there are isomorphisms%
\[
K_{n}^{\mathrm{top}}(A)\cong K_{n+p}^{\mathrm{top}}(A)\text{,}%
\]
where $K_{n}^{\mathrm{top}}(A)=\pi_{n-1}(\mathrm{GL}(A))$ if $n>0$ and
$K_{0}^{top}(A)=K(A)$ is the usual Grothendieck group. Here $p$ is the period
which is $2$ or $8$ according as $A$ is complex or real. We refer to
\cite{KaroubiHandbook} and \cite{Rosenberg} for an overview of the subject,
both algebraically and topologically.

A few years later, after higher algebraic $K$-theory was introduced by
Quillen, an analogous periodicity statement was sought, of the form%
\[
K_{n}(A)\cong K_{n+p}(A)\text{,}%
\]
where $A$ is now a discrete ring. The first computations showed that a
periodicity isomorphism of this form is far from true in basic examples.
However, if we consider $K$-theory with finite coefficients, and $n$ is at
least a certain bound $d$, then some periodicity conjectures appeared
feasible, at least for certain rings of a geometric nature. These conjectures
were formulated for different prime power coefficient groups, and are
essentially of the following type ($n\geq d$)
\[
K_{n}(A;\,\mathbb{Z}/m)\cong K_{n+p}(A;\,\mathbb{Z}/m)\text{.}%
\]

The relationship between the prime power $m$ and the associated smallest
period $p$ is given by the following convention, which we maintain throughout
the paper.

\begin{Convention}
\label{convention: (p,p')}For $\mathbb{Z}/m$ coefficients, where $m=\ell^{\nu
}$ with $\ell$ prime, the smallest period $p$ is given by%
\[
p=\left\{
\begin{array}
[c]{lll}%
\mathrm{sup}\left(  8,\,\ell^{\nu-1}\right)  & \quad & \text{if }%
\ell=2\text{,}\\
2(\ell-1)\ell^{\nu-1} &  & \text{otherwise.}%
\end{array}
\right.
\]

\end{Convention}

Using techniques of algebraic geometry and a comparison
theorem with \'{e}tale $K$-theory, numerous examples listed below showed that
these conjectures hold. In the case of a $2$-power, the first three are
particular cases of Theorem 2 in \cite{Ostvar-Rosenschon}, based on the
fundamental work of Voevodsky \cite{Voevodskymod2}. In the case of an odd
prime power, the first four examples are consequences of the Bloch-Kato
conjecture.

Before giving these examples, we define the mod\ $2$ virtual \'{e}tale
cohomological dimension \textrm{vcd}$_{2}(A)$ of a commutative ring $A$ as the
mod\ $2$ \'{e}tale cohomological dimension of $A\otimes_{\mathbb{Z}}%
\mathbb{Z}[\mu_{4}]$ obtained by adjoining a primitive fourth root of unity to
$A$. For convenience, if $\ell$ is odd, then \textrm{vcd}$_{\ell}(A)$ denotes
the mod\ $\ell$ \'{e}tale cohomological dimension \textrm{cd}$_{\ell}(A)$ of
$A$. For $\ell$ fixed, here are the examples we consider.

\begin{enumerate}
\item Any field $k$ of characteristic $\operatorname{char}(k)\neq\ell$ for
which $\mathrm{vcd}_{\ell}(k)<\infty$. In this case, $d=\mathrm{vcd}_{\ell
}(k)-1$ if $\mathrm{vcd}_{\ell}(k)\neq0$ and $d=0$ otherwise.

\item The ring $\mathcal{O}_{F}\left[  1/\ell\right]  $ of $\ell$-integers in
any number field $F$. In this case, $d=\mathrm{vcd}_{\ell}(F)-1=1$
(\textsl{cf.}~\cite{Ostvar} when $\ell=2$).

\item Any finitely generated and regular $\mathbb{Z}\left[  1/\ell\right]
$-algebra $A$ with finite mod\ $\ell$ virtual \'{e}tale cohomological
dimension. In this case $d=\sup\{\mathrm{vcd}_{\ell}(k(s))-1,0\}$, where
$k(s)$ is the residue field at any point $s\in\text{Spec}(A)$. The same
statement holds when replacing $\mathbb{Z}\left[  1/\ell\right]  $ by
$\mathbb{Q}$ or by any other field $k$ of characteristic $\neq\ell$.
\newline\noindent The regularity assumption on $A$ can be dispensed with when
working with negative $K$-theory \cite{Bass}, \cite{Karoubi7}, \cite{Karoubi9}%
, \cite{TT}. As shown in \cite[Theorem 4.5]{Ostvar-Rosenschon2}, this does not
change the bound $d$.

\item Group rings $R\left[  G\right]  $, where $G$ is finite and $R$ is a ring
of $\ell$-integers in a number field, as shown in \cite{WeibelGrouprings}.
Here $d=1$. For some explicit computations see \cite{Ostvargrouprings}.

\item The ring $C(X)$ of real or complex continuous functions on a compact
space $X$, as shown in \cite{Fisher}, \cite{Prasolov}. In this case $d=1$.
\end{enumerate}

In these examples, the periodicity isomorphism between the groups
$K_{n}(A;\,\mathbb{Z}/m)$ and $K_{n+p}(A;\,\mathbb{Z}/m)$ is defined by taking
cup-product with a \textquotedblleft Bott element\textquotedblright\ $b_{K}$.
For $p=2^{\nu-1}$ with $\nu\geq4$, one can construct this element in the group
$K_{p}(\mathbb{Z};\,\mathbb{Z}/2p)$, such that its image in the topological
$K$-group\footnote{We shall write $K_{n}$ instead of $K_{n}^{\mathrm{top}}$
when dealing with the field $\mathbb{R}$ of real numbers or the field
$\mathbb{C}$ of complex numbers with their usual topology, and likewise for
spectra.} $K_{p}(\mathbb{R};\,\mathbb{Z}/2p)\cong\mathbb{Z}/2p$ is the class
$\operatorname{mod}2p$ of a generator in $K_{p}(\mathbb{R})\cong\mathbb{Z}$.
The cup-product alluded to above is a pairing%
\[
_{\cup}:K_{n}(A;\,\mathbb{Z}/m)\times K_{p}(\mathbb{Z};\,\mathbb{Z}%
/2p)\longrightarrow K_{n+p}(A;\,\mathbb{Z}/m)\text{.}%
\]
We refer to Section \ref{section: Bott elements} for precise definitions and
the extension to odd prime powers.

As we can see in these examples, a key role is played by the infinite general
linear group $\mathrm{GL}(A)$. However, it was already shown in the works of
Bott and Borel \cite{Borel}, and also in topological applications, that other
infinite series associated to classical Lie groups may be considered as well.
More precisely, if we consider a ring with involution $A$ and a sign of
symmetry $\varepsilon=\pm1$ generalizing the orthogonal ($\varepsilon=1$) or
symplectic ($\varepsilon=-1$) case, one defines higher hermitian $K$-groups,
denoted in this paper ${}_{\varepsilon}KQ_{n}(A)$, in a parallel way to
algebraic $K$-groups $K_{n}(A)$. These groups are associated to the infinite
$\varepsilon$-orthogonal group {}$_{\varepsilon}O(A)$. A typical example is
when $A$ is commutative and $\varepsilon=-1$, in which case one recovers the
infinite symplectic group on $A$. We refer to the survey paper
\cite{KaroubiHandbook} already mentioned above for precise definitions.

The main purpose of this paper is to show that a periodicity statement in
algebraic $K$-theory implies a similar one in $KQ$-theory, when $1/2\in A$.
Since $KQ$-theory with coefficients $\mathbb{Z}/m$, with $m=2^{\nu}$, is the
most important and difficult case, we state the main theorems in this context,
leaving the case of odd prime power coefficients to the end of this
Introduction and to Section \ref{section: odd prime} of the main body of the paper.

For the first step in the argument, we introduce a parameter $q$ that is
essentially $p$, apart from a slight modification in the case $m=16$.
Specifically, we make the following convention.

\begin{Convention}
\label{stronger convention}%
\[
q=\left\{
\begin{array}
[c]{lll}%
8 & \quad & \text{if }m\leq8\text{,}\\
16 &  & \text{if }m=16\text{,}\\
m/2 &  & \text{otherwise.}%
\end{array}
\right.
\]

\end{Convention}

In other words, $q=p$ except when $m=16$, in which case $q=2p$. It is
meaningful to speak of periodicity maps raising dimension by $q$, since $q$ is
a multiple of $p$.

As a convenient notation, we write $\overline{KQ}$
(resp. $\overline{K}$) for the $KQ$-theory (resp. $K$-theory) with
coefficients in $\mathbb{Z}/m$, the relationship between $m$ and the period
$p$ being as in Convention \ref{convention: (p,p')}. One of our main theorems
is the following.

\begin{theorem}
\label{abstract KQ-periodicity} With the above definitions, assume that there
exists an integer $d$ such that the cup-product map
\[
_{\cup}b_{K}:\overline{K}_{n}(A)\longrightarrow\overline{K}_{n+p}(A)
\]
with the Bott element in {}$K_{p}(\mathbb{Z};\,\mathbb{Z}/2p)$ is an
isomorphism whenever $n\geq d$. Then, for $n\geq d+q-1$, there is also an
isomorphism
\[
{}_{\varepsilon}\overline{KQ}_{n}(A)\cong{}_{\varepsilon}\overline{KQ}%
_{n+p}(A)\text{.}%
\]

\end{theorem}

Surprisingly, the isomorphism between the $KQ$-groups is in general not given
by cup-product with a Bott element (see Remark \ref{remark: theta not iso} in
Section 4). This relates to the fact that hermitian $K$-theory possesses more
than one Bott element, as we now describe. Whereas in algebraic $K$-theory
universal Bott elements are to be found in the $K$-groups of the integers
$\mathbb{Z}$, here, because we are working with rings containing $1/2$, our
Bott elements are to be found in the hermitian $K$-groups of the ring of
$2$-integers $\mathbb{Z}^{\prime}=\mathbb{Z}\left[  1/2\right]  $.

As in algebraic $K$-theory, using the methods of \cite{BK}, in this paper we
prove the existence of a ``positive\ Bott element''\ $b^{+}$\ in {}$_{1}%
KQ_{p}(\mathbb{Z}^{\prime};\,\mathbb{Z}/2p)$ whose image in {}$K_{p}%
(\mathbb{Z}^{\prime};\,\mathbb{Z}/2p)\cong K_{p}(\mathbb{Z};\,\mathbb{Z}/2p)$
is the Bott element in $K$-theory alluded to above.

On the other hand, one of the main differences between algebraic and hermitian
$K$-theory in our context is the existence of another element\footnote{We
recall that the negative $K$-groups of a regular noetherian ring (for instance
$\mathbb{Z}$ or $\mathbb{Z}^{\prime}$) are trivial.} $u$\ in $_{-1}%
KQ_{-2}(\mathbb{Z}^{\prime})$, which plays an important role in the
fundamental theorem in hermitian $K$-theory \cite{KaroubiAnnals2}. We now
define the \emph{negative\ Bott element}\ $b^{-}$ in
hermitian $K$-theory to be the image of the element $u^{p/2}$ in the group
${}_{1}KQ_{-p}(\mathbb{Z}^{\prime};\,\mathbb{Z}/2p)$.

To make the statement of Theorem \ref{abstract KQ-periodicity} more precise,
we note that the cup-product with the positive Bott element in {}$_{1}%
KQ_{p}(\mathbb{Z}^{\prime};\,\mathbb{Z}/2p)$ determines a direct system of
abelian groups
\[
{}_{\varepsilon}\overline{KQ}_{n}(A)\longrightarrow{}_{\varepsilon}%
\overline{KQ}_{n+p}(A)\longrightarrow{}_{\varepsilon}\overline{KQ}%
_{n+2p}(A)\longrightarrow\cdots\text{ .}%
\]
Symmetrically, cup-product with the negative Bott element in {}$_{1}%
KQ_{-p}(\mathbb{Z}^{\prime};\,\mathbb{Z}/2p)$ determines an inverse system of
abelian groups
\[
\cdots\longrightarrow{}_{\varepsilon}\overline{KQ}_{n+2p}(A)\longrightarrow
{}_{\varepsilon}\overline{KQ}_{n+p}(A)\longrightarrow{}_{\varepsilon}%
\overline{KQ}_{n}(A)\text{.}%
\]

The theorem above can now be restated in a more precise form. (Recall that the
overbar denotes $\mathbb{Z}/m$ coefficients.)

\begin{theorem}
\label{Thm: invlimKQdirlim}Let $A$ be any ring (with $1/2\in A$), $m$, $p$ and
$q$ be $2$-powers as in Conventions \ref{convention: (p,p')} and
\ref{stronger convention}, and let $d\in\mathbb{Z}$, such that the cup-product
with the Bott element $b_{K}$ in $K_{p}(\mathbb{Z};\,\mathbb{Z}/2p)$ induces
an isomorphism%
\[
\overline{K}_{n}(A)\overset{\cong}{\longrightarrow}\overline{K}_{n+p}(A)
\]
whenever $n\geq d$. Then, for $n\geq d$, there is an exact
sequence%
\[
\cdots\overset{\theta^{-}}{\longrightarrow}{}_{\varepsilon}\overline{KQ}%
_{n+1}(A)\overset{\theta^{+}}{\longrightarrow}{}\underrightarrow{\lim}%
{}_{\varepsilon}\overline{KQ}_{n+1+ps}(A)
\]%
\[
\rightarrow\underleftarrow{\lim}{}_{\varepsilon}\overline{KQ}_{n+ps}%
(A)\overset{\theta^{-}}{\longrightarrow}{}_{\varepsilon}\overline{KQ}%
_{n}(A)\overset{\theta^{+}}{\longrightarrow}{}\underrightarrow{\lim}%
{}_{\varepsilon}\overline{KQ}_{n+ps}(A)
\]
where $\theta^{+}$ (respectively $\theta^{-}$) is induced from the cup-product
with the positive Bott element $b^{+}$ (resp. the negative Bott element
$b^{-}$). Moreover, for $n\geq d+q-1$, there is a short split exact sequence%
\[
0\rightarrow\underleftarrow{\lim}{}_{\varepsilon}\overline{KQ}_{n+ps}%
(A)\overset{\theta^{-}}{\longrightarrow}{}_{\varepsilon}\overline{KQ}%
_{n}(A)\overset{\theta^{+}}{\longrightarrow}{}\underrightarrow{\lim}%
{}_{\varepsilon}\overline{KQ}_{n+ps}(A)\rightarrow0.
\]

\end{theorem}

It turns out that the inverse limit is not always trivial. This point is
discussed in Section
\ref{Proof of the periodicity theorem for totally real 2-regular number fields}
(where the inverse limit vanishes) and Section \ref{section: proof of main}
(where it does not).

However, for rings of geometric nature and of finite mod\ $2$ virtual
\'{e}tale cohomological dimension, we conjecture that the inverse limit is trivial.

\begin{definition}
\label{defn: hermitian regular}We say that a ring $A$
is \emph{hermitian regular}\ if $\underleftarrow{\lim}\,_{\varepsilon
}\overline{KQ}_{n+ps}(A)$ and $\underleftarrow{\lim}^{1}\,_{\varepsilon
}\overline{KQ}_{n+ps}(A)$ are trivial\footnote{As a matter of fact, with our
hypothesis about periodicity of the $K$-groups, we always have $\lim^{1}=0,$
since the inverse system satisfies the Mittag-Leffler condition as we shall
see in Sections \ref{Higher KSC-theories} and \ref{section: proof of main}.}.
\end{definition}

\begin{remark} It should be noted that subsequent to the original 
submission of this paper at the
beginning of February 2010, as a consequence of a more recent theorem of Hu,
Kriz and Ormsby \cite{HuKrizOrmsby} in characteristic $0$, the authors and
M.\thinspace Schlichting proved independently that a field of characteristic
$0$ that is of finite mod\ $2$ virtual \'{e}tale cohomological definition is
hermitian regular. Furthermore, Schlichting extended this theorem for fields
of characteristic $p>0$ with the same cohomological properties. This affirms
Conjecture \ref{conjecture:comparisonforfields} of the present paper, which
implies in turn our Conjecture \ref{conjectrure} and therefore considerably
extends the number of examples of commutative rings (and schemes) that are
hermitian regular. The details of the proofs will appear in a forthcoming
joint paper of the authors and Schlichting \cite{Schlichting2}. A particular
example quoted below is given by suitable rings of integers in a number field.
In Theorem \ref{totally real theorem}, we state the periodicity theorem in
this case with an independent proof which will be given in Section
\ref{Proof of the periodicity theorem for totally real 2-regular number fields}. 
A more general theorem is as follows.
\end{remark}

\begin{theorem}
\label{W-regularity}Let $A$ be a ring which is hermitian regular and satisfies
the hypothesis of the previous theorem for its $K$-groups. Then for $n\geq d$,
the cup-product with the positive Bott element induces an isomorphism
\[
{}_{\varepsilon}\overline{KQ}_{n}(A)\overset{\cong}{\longrightarrow}{}%
{}_{\varepsilon}\overline{KQ}_{n+p}(A).
\]

\end{theorem}

More generally, in order to fully exploit the spectrum approach and to improve
the previous theorems, we may consider a pointed CW-complex $X$ and define the
group $K_{X}(A)$ as the group of homotopy classes of pointed maps from $X$ to
$\mathcal{K}(A)$, where $\mathcal{K}$ denotes the $K$-theory spectrum. If $X$
is a pointed sphere $S^{n}$, we recover Quillen's $K$-group $K_{n}(A)$. For
brevity, we shall also write $K_{X+t}(A)$ instead of $K_{X\wedge S^{t}}(A)$,
and $K_{X-t}(A)$ instead of $K_{X}(S^{t}A)$, where $S^{t}A$ denotes the
$t$-iterated suspension of $A$ (see for instance \cite{KaroubiHandbook} for
the definition of the suspension and the basic definitions of various
$K$-theories). We adopt the same conventions for hermitian $K$-theory, and
also for algebraic or hermitian $K$-theory with coefficients, and finally for spectra.

The previous theorem can now be generalized as follows.

\begin{theorem}
\label{periodicity X}Let $A$ be any ring (with $1/2\in A$), $m$, $p$ and $q$
be $2$-powers as in Conventions \ref{convention: (p,p')} and
\ref{stronger convention}, and let $d\in\mathbb{Z}$, such that the cup-product
with the Bott element $b_{K}$ in $K_{p}(\mathbb{Z};\,\mathbb{Z}/2p)$ induces
an isomorphism%
\[
\overline{K}_{n}(A)\overset{\cong}{\longrightarrow}\overline{K}_{n+p}(A)
\]
whenever $n\geq d$. Then, if $X$ is a $(d-1)$-connected space, there is an
exact sequence%
\[
{}_{\varepsilon}\overline{KQ}_{X+1}(A)\overset{\theta^{+}}{\longrightarrow}%
{}\underrightarrow{\lim}{}_{\varepsilon}\overline{KQ}_{X+1+ps}(A)
\]%
\[
\longrightarrow\underleftarrow{\lim}{}{}_{\varepsilon}\overline{KQ}%
_{X+ps}(A)\overset{\theta^{-}}{\longrightarrow}{}_{\varepsilon}\overline
{KQ}_{X}(A)\overset{\theta^{+}}{\longrightarrow}{}\underrightarrow{\lim}%
{}_{\varepsilon}\overline{KQ}_{X+ps}(A)\rightarrow\cdots.
\]
If $X$ is $(d+q-2)$-connected, there is a split short exact sequence%
\[
0\rightarrow\underleftarrow{\lim}{}{}_{\varepsilon}\overline{KQ}%
_{X+ps}(A)\overset{\theta^{-}}{\longrightarrow}{}_{\varepsilon}\overline
{KQ}_{X}(A)\overset{\theta^{+}}{\longrightarrow}{}\underrightarrow{\lim}%
{}_{\varepsilon}\overline{KQ}_{X+ps}(A)\rightarrow0\text{.}%
\]
Finally, if $A$ is hermitian regular and if $X$ is $(d-1)$-connected, the
cup-product with the positive Bott element induces an isomorphism%
\[
{}_{\varepsilon}\overline{KQ}_{X}(A)\cong{}_{\varepsilon}\overline{KQ}%
_{X+p}(A)\text{.}%
\]

\end{theorem}

\begin{corollary}
For any $(d+q-2)$-connected space $X$ and $A$ as above (not necessarily
hermitian regular), there is a periodicity isomorphism%
\[
_{\varepsilon}\overline{KQ}_{X}(A)\cong{}_{\varepsilon}\overline{KQ}%
_{X+p}(A)\text{.}%
\]

\end{corollary}

For suitable subrings $A_{S}$ in a number field $F$, the previous results may
be stated more precisely, by using the methods of \cite{BKO}. The rings
$A_{S}$, defined in Section
\ref{Proof of the periodicity theorem for totally real 2-regular number fields}
below, generalize both the ring of $S$-integers (when $S$ is finite) and the
number field $F$ itself. (More general examples are considered in Section
\ref{Generalization to schemes} and in \cite{Schlichting2}.)

\begin{theorem}
\label{totally real theorem}Let $F$ be a totally real $2$-regular number field
as considered in \cite{BKO}; also, let $m$ and $p$ be $2$-powers as in
Convention \ref{convention: (p,p')}. Then, for all integers $n>0$, the inverse
limit $\underleftarrow{\lim}{}_{\varepsilon}\overline{KQ}_{n+ps}(A_{S})$ is
trivial (\textsl{i.e.} $A_{S}$ is hermitian regular) and the \textquotedblleft
positive\textquotedblright\ Bott map%
\[
\beta_{n}=\,_{\cup}b^{+}:{}_{\varepsilon}\overline{KQ}_{n}(A_{S}%
)\longrightarrow{}_{\varepsilon}\overline{KQ}_{n+p}(A_{S})
\]
is an isomorphism. More generally, if $X$ is any connected CW-complex, the
Bott map
\[
\beta_{X}:{}_{\varepsilon}\overline{KQ}_{X}(A_{S})\longrightarrow
{}_{\varepsilon}\overline{KQ}_{X+p}(A_{S})
\]
is an isomorphism.
\end{theorem}

For completeness we mention the odd-primary analog of Theorem
\ref{Thm: invlimKQdirlim}, which is proved in Section \ref{section: odd prime}%
. Its applications are related to the Bloch-Kato conjecture as we mentioned at
the beginning. We note that the hypothesis $1/2\in A$ may be dropped in this case.

\begin{theorem}
Let $p$ and $m$ be odd prime powers as in Convention \ref{convention: (p,p')}.
Let $b_{K}$ be the associated Bott element\ in $K_{p}(\mathbb{Z}%
;\,\mathbb{Z}/m)$ (see Section \ref{section: Bott elements} for details). Now
let $A$ be any ring and assume that, whenever $n\geq d$, cup-product with
$b_{K}$ induces an isomorphism
\[
\overline{K}_{n}(A)\cong\overline{K}_{n+p}(A)\text{.}%
\]
Then there exists a \textquotedblleft mixed Bott element\textquotedblright%
\ $b$ in $_{1}\overline{KQ}_{p}(\mathbb{Z}^{\prime})$ such that for $n\geq d$,
the cup-product with $b$ induces an isomorphism between the related
$KQ$-groups
\[
\beta_{n}:{}_{\varepsilon}\overline{KQ}_{n}(A)\overset{\cong}{\longrightarrow
}{}_{\varepsilon}\overline{KQ}_{n+p}(A)\text{.}%
\]
More generally, if $X$ is a ($d-1$)-connected CW complex, then the cup-product
map with $b$ induces an isomorphism%
\[
\beta_{X}:{}_{\varepsilon}\overline{KQ}_{X}(A)\overset{\cong}{\longrightarrow
}{}_{\varepsilon}\overline{KQ}_{X+p}(A)\text{.}%
\]

\end{theorem}

In Section \ref{Generalization to schemes} we note that work in progress by
Schlichting \cite{Schlichting} allows us to extend our results from
commutative rings to schemes $S$ that are separated, noetherian and of finite
Krull dimension. More precisely, following Jardine's method for algebraic
$K$-theory \cite{Jardinebook} we define an \textquotedblleft%
\'{e}tale\textquotedblright\ $\overline{KQ}$-theory, denoted by ${}%
_{\varepsilon}\overline{KQ}_{n}^{\mathrm{\acute{e}t}}(S)$, where the
coefficient groups are prime powers. The \'{e}tale $\overline{KQ}$-theory
shares many properties with the \'{e}tale $\overline{K}$-theory introduced by
Dwyer and Friedlander \cite{DF}. For example, there exists a comparison map
\[
\sigma:{}_{\varepsilon}\overline{KQ}_{n}(S)\longrightarrow{}_{\varepsilon
}\overline{KQ}_{n}^{\mathrm{\acute{e}t}}(S)\text{.}%
\]

For odd prime powers, there is an involution on the odd torsion group
${}_{\varepsilon}\overline{KQ}_{n}(S)$. Let
${}_{\varepsilon}\overline{KQ}_{n}^{\mathrm{\acute{e}t}}(S)_{+}$ and
${}_{\varepsilon}\overline{KQ}_{n}^{\mathrm{\acute{e}t}}(S)_{-}$ denote the
corresponding eigenspaces. On the other hand, for any prime power (odd or
even), the cup-product map with the Bott element $b^{+}$ defined above induces
a direct system of groups, whose colimit we shall denote by ${}_{\varepsilon
}\overline{KQ}_{n}(S)\left[  \beta^{-1}\right]  $, in the notation of
\cite{Thomason}. Next, we state two theorems and a conjecture in this context.

\begin{theorem}
With the coefficient group $\mathbb{Z}/\ell^{\nu}$, where $\ell$ is an odd
prime, there is an isomorphism
\[
{}_{\varepsilon}\overline{KQ}_{n}(S)\left[  \beta^{-1}\right]  \cong%
{}_{\varepsilon}\overline{KQ}_{n}^{\mathrm{\acute{e}t}}(S)
\]
for all $n$ if $\mathrm{cd}_{\ell}(S)<\infty$. Moreover, the comparison map
$\sigma$ induces an isomorphism
\[
{}_{\varepsilon}\overline{KQ}_{n}(S)_{+}\cong{}_{\varepsilon}\overline{KQ}%
_{n}^{\mathrm{\acute{e}t}}(S)
\]
for $n\geq\sup\{\mathrm{cd}_{\ell}(k(s))-1\}_{s\in S}$.
\end{theorem}

Recall that $S$ is uniformly $\ell$-bounded with bound $d$ if for all residue
fields $k(s)$ we have $\mathrm{cd}_{\ell}(k(s))\leq d$. In the event that $S$
is uniformly $\ell$-bounded with bound $d$, then $\mathrm{cd}_{\ell}(S)\leq
n+d$ where $n$ denotes the Krull dimension of $S$; an elegant proof for this
inequality is given in \cite[Theorem 2.8]{Mitchellhypercohomology}.

At the prime $2$ we prove the following theorem, reminiscent of the main
results in \cite{DFST} and in \cite{Thomason}.

\begin{theorem}
\label{theorem:splotsurjectivity} With the coefficient group $\mathbb{Z}%
/2^{\nu}$, there is an isomorphism
\[
_{\varepsilon}\overline{KQ}_{n}(S)[\beta^{-1}]\cong{}_{\varepsilon}%
\overline{KQ}_{n}^{\mathrm{\acute{e}t}}(S)
\]
for all $n$ if $\mathrm{vcd}_{2}(S)<\infty$. Moreover, the comparison map
\[
\sigma:{}_{\varepsilon}\overline{KQ}_{n}(S)\longrightarrow{}_{\varepsilon
}\overline{KQ}_{n}^{\mathrm{\acute{e}t}}(S)
\]
is a split surjection for $n\geq\sup\{\mathrm{vcd}_{2}(k(s))-1\}_{s\in S}+q-1$.
\end{theorem}

More generally, we make the following conjecture.

\begin{conjecture}
\label{conjectrure} With the coefficient group $\mathbb{Z}/2^{\nu}$ the map
$\sigma$ is bijective whenever $n\geq\sup\{\mathrm{vcd}_{2}(k(s))-1\}_{s\in
S}$.
\end{conjecture}

Using algebro-geometric methods, in Theorem
\ref{comparisontheorem} below we show how to reduce this conjecture to the
case of fields.
As mentioned above, the characteristic $0$ case was solved independently by
the authors and M.\thinspace Schlichting, while the positive characteristic
case was solved by Schlichting. A proof of this conjecture in general will
appear in a joint paper with Schlichting \cite{Schlichting2}.

\medskip

Let us now briefly discuss the contents of the paper.

In Section \ref{section: Bott elements}, for $2$-power coefficients we
carefully construct the Bott elements that play an important role in this
work, as referred to above.

Section
\ref{Proof of the periodicity theorem for totally real 2-regular number fields}
is somewhat independent of the other sections. In
particular, we prove a refined version of our theorems in the case $A$ is the
ring of integers in a totally real $2$-regular number field. 
(This version is a particular case of the considerations in 
Section \ref{Generalization to schemes} for schemes. 
Assuming Conjecture \ref{conjectrure},  
which will be proven in \cite{Schlichting2},
Theorem \ref{2-regular main} may be given an independent proof in a much more general framework.)

In Section \ref{Higher KSC-theories}, we introduce what we call
\textquotedblleft higher $\overline{KSC}$-theories\textquotedblright. These
theories in some sense measure the deviation of \textquotedblleft
negative\textquotedblright\ periodicity of the $\overline{KQ}$-groups. On the
other hand, they are built by successive extensions of the $\overline{K}%
$-groups. Therefore, they are periodic if the $\overline{K}$-groups are periodic.

Section \ref{section: proof of main} is devoted to the proof of our main
Theorems \ref{Thm: invlimKQdirlim} and \ref{Theorem of W-regularity} (for
arbitrary rings with $2$ invertible and mod\ $2^{\nu}$ coefficients). The
proof is roughly divided into two steps as follows. In the first one, we prove
a cruder periodicity statement for $n\geq d+q-1$. In the second, we use the
$\overline{KQ}$-spectrum and an argument about cohomology theories to prove
the periodicity theorems in full generality. We conclude this section with an
upper bound of the $\overline{KQ}$-groups in terms of the $\overline{K}$-groups.

In Section \ref{section: odd prime}, we study the case of odd prime powers,
which is paradoxically simpler in our framework. The main observation is that
the $\overline{KQ}$-ring spectrum splits naturally as the product of two ring
spectra, the first one being the \textquotedblleft symmetric\textquotedblright%
\ part of the $\overline{K}$-theory spectrum.

Section \ref{Generalization to schemes} is more geometric in nature and
generalizes the previous considerations (when $A$ is commutative) to
noetherian separated schemes of finite Krull dimension. Here we rely heavily
on the fundamental theorem in hermitian $K$-theory proved in the scheme
framework by Schlichting \cite{Schlichting}.

Finally, Sections \ref{section:applications} and \ref{section:applications:2}
are devoted to selected applications: rings of integers in number fields,
smooth complex algebraic varieties, and rings of continuous functions on
compact spaces. Another application, to hermitian $\overline{KQ}$-theory of
group rings, is a consequence of an appendix to this paper by C.\thinspace
Weibel \cite{WeibelGrouprings}.\smallskip

\textbf{Acknowledgements}. We warmly thank the referee for very relevant
comments on a previous version of this paper. We extend our thanks to Marco
Schlichting for  discussions resulting in our joint
work \cite{Schlichting2}.

\section{Bott elements in $K$- and $KQ$-theories}

\label{section: Bott elements}

Let $\ell$ be a prime number and $\mathcal{S}^{0}/\ell^{\nu}$ the
mod\ $\ell^{\nu}$ Moore spectrum. In \cite[\S 12]{Adams}, Adams constructed
$\mathcal{KO}_{\ast}$-equivalences
\[
A_{\ell^{\nu}}:\Sigma^{p}\mathcal{S}^{0}/\ell^{\nu}\longrightarrow
\mathcal{S}^{0}/\ell^{\nu}\text{.}%
\]
The dimension shift $p$ is $\sup\{8,2^{\nu-1}\}$ if $\ell=2$ and
$2(\ell-1)\ell^{\nu-1}$ if $\ell$ is odd.
As shown by Bousfield in \cite[\S 4]{Bousfield}, work of Mahowald and Miller
implies that a spectrum $\mathcal{E}$ is $\mathcal{KO}$-local if and only if
its mod\ $\ell$ homotopy groups are periodic via $A_{\ell}$ for every prime
$\ell$. We shall refer to the periodicity manifested in $\mathcal{KO}$-local
spectra as Bott periodicity. Note that $\mathcal{KO}$-localizations are the
same as $\mathcal{KU}$-localizations \cite[\S 4]{Bousfield}.


In general there are several choices of an element $A_{\ell^{\nu}}$ as above
if the only criterion is that it induces a $\mathcal{KO}$-isomorphism. We are
interested in particular choices of elements pertaining to classical Bott
periodicity. Let $u$ denote a generator of the infinite cyclic group $\pi
_{2}(BU)$. Then for $r\geq1$ the Bott element $u^{2r}$ in $\pi_{4r}(BU)$ is
independent of the choice of $u$. We denote by $v$ the element of $\pi
_{8r}(BO)$ mapping to the Bott element in $\pi_{8r}(BU)$ under the map induced
by complexification $c:BO\rightarrow BU$.

The mod\ $2^{\nu}$ Bott element in degree $8r>0$ is the generator
\[
\overline{v}=\operatorname{id}_{\mathcal{S}^{0}/2^{\nu}}\wedge v\in
KO_{8r}(\mathcal{S}^{0}/2^{\nu};\,\mathbb{Z}/2^{\nu})=[\mathcal{S}^{0}/2^{\nu
},\,\mathcal{KO}\wedge\mathcal{S}^{0}/2^{\nu}]_{8r}\text{.}%
\]
The element $A_{2^{\nu}}$ is called an Adams periodicity operator if it maps
to the mod\ $2^{\nu}$ Bott element in degree $p$ under the naturally induced
$\mathcal{KO}$-Hurewicz map
\[
\pi_{\ast}(\mathcal{S}^{0}/2^{\nu};\,\mathbb{Z}/2^{\nu})\rightarrow
\mathcal{KO}_{\ast}(\mathcal{S}^{0}/2^{\nu};\,\mathbb{Z}/2^{\nu})
\]
for $\mathcal{S}^{0}/2^{\nu}$. When $\ell\neq2$, the definition of a
mod\ $\ell^{\nu}$ Bott element is the same as above, except that
$\mathcal{KO}$ is replaced by $\mathcal{KU}$. Crabb and Knapp
\cite{CrabbKnapp} have shown that there exist Adams periodicity operators for
all $\ell$ and $\nu\geq1$.

By smashing the unit map $\mathcal{S}^{0}\rightarrow\mathcal{E}$ of a ring
spectrum $\mathcal{E}$ with $\mathcal{S}^{0}/\ell^{\nu}$ and pushing forward
the class in $\pi_{p}(\mathcal{S}^{0}/\ell^{\nu};\,\mathbb{Z}/\ell^{\nu})$
represented by the map $A_{\ell^{\nu}}$, one obtains a class in the group
$\pi_{p}(\mathcal{E};\,\mathbb{Z}/\ell^{\nu})$ that we call a Bott element.

Next, for $m=2p,$ where $p=2^{\upsilon-1}$ is a $2$-power $\geq8$, we study
mod\ $m$ Bott elements in more detail for $K$- and $KQ$-theory in the example
of $\mathbb{Z}^{\prime}$. The case of an odd prime is dealt with in Section
\ref{section: odd prime}.

To begin, we shall consider \textquotedblleft Bott elements\textquotedblright%
\ in $K_{p}(\mathbb{Z}^{\prime};\,\mathbb{Z}/m)$ and ${}_{1}KQ_{p}%
(\mathbb{Z}^{\prime};\,\mathbb{Z}/m)$, whose images in $K_{p}(\mathbb{R}%
;\,\mathbb{Z}/m)$ and ${}_{1}KQ_{p}(\mathbb{R};\,\mathbb{Z}/m)$, respectively,
are generators deduced from classical Bott periodicity for the real numbers
(as $KQ$-modules). This is well-known for the algebraic
$K$-groups; it is included here for the sake of completeness.

B{\"{o}}kstedt's square of algebraic $K$-theory spectra introduced in
\cite{Bokstedt}
\[%
\begin{array}
[c]{ccc}%
\mathcal{K(}\mathbb{Z}^{\prime})_{\#} & \longrightarrow & \mathcal{K(}%
\mathbb{R})_{\#}^{c}\\
\downarrow &  & \downarrow\\
\mathcal{K(}\mathbb{F}_{3})_{\#} & \longrightarrow & \mathcal{K(}%
\mathbb{C})_{\#}^{c}%
\end{array}
\]
was verified to be homotopy cartesian by Rognes-Weibel in \cite{RW},
\cite{WeibelCRAS}, as a consequence of Voevodsky's proof of the Milnor
conjecture. Here $_{\#}$ means $2$-adic completions and $^{c}$ means
connective cover. Smashing with $\mathcal{S}^{0}/2^{\nu}$ yields a homotopy
cartesian square (an overbar indicates reduction mod\ $m$):
\[%
\begin{array}
[c]{ccc}%
\mathcal{\overline{K}(}\mathbb{Z}^{\prime}) & \longrightarrow &
\mathcal{\overline{K}(}\mathbb{R})^{c}\\
\downarrow &  & \downarrow\\
\mathcal{\overline{K}(}\mathbb{F}_{3}) & \longrightarrow & \mathcal{\overline
{K}(}\mathbb{C})^{c}%
\end{array}
\]
Denote by $\overline{K}$ the corresponding mod\ $m$ homotopy groups. By Bott
periodicity and the isomorphism $\overline{K}_{p-1}(\mathbb{Z}^{\prime
})\rightarrow\overline{K}_{p-1}(\mathbb{F}_{3})$, there is a split short exact
sequence
\[
0\longrightarrow\overline{K}_{p}(\mathbb{Z}^{\prime})\longrightarrow
\overline{K}_{p}(\mathbb{R})\oplus\overline{K}_{p}(\mathbb{F}_{3}%
)\longrightarrow\overline{K}_{p}(\mathbb{C})\longrightarrow0\text{.}%
\]
On the other hand, Quillen's homotopy fibration%
\[
\Omega\mathcal{\overline{K}}\mathbb{(C)}\overset{\Psi^{3}-1}{\longrightarrow
}\Omega\mathcal{\overline{K}}\mathbb{(C)}\longrightarrow\mathcal{\overline
{K}(}\mathbb{F}_{3})\longrightarrow\mathcal{\overline{K}}\mathbb{(C)}%
\overset{\Psi^{3}-1}{\longrightarrow}\mathcal{\overline{K}}\mathbb{(C)}%
\]
yields an exact sequence%
\[
\overline{K}_{p+1}(\mathbb{C})\overset{\cdot m}{\longrightarrow}\overline
{K}_{p+1}(\mathbb{C})\longrightarrow\overline{K}_{p}(\mathbb{F}_{3}%
)\longrightarrow\overline{K}_{p}(\mathbb{C})\overset{\cdot m}{\longrightarrow
}\overline{K}_{p}(\mathbb{C})\text{,}%
\]
and hence the isomorphisms
\[
\overline{K}_{p}(\mathbb{F}_{3})\cong\,_{m}\overline{K}_{p}(\mathbb{C}%
)\cong\mathbb{Z}/m\text{.}%
\]
Here $_{n}A$ denotes the kernel of the multiplication by $n$ map on an abelian
group $A$. Hence, diagram chasing shows there are isomorphisms
\[
\overline{K}_{p}(\mathbb{Z}^{\prime})\cong\overline{K}_{p}(\mathbb{R}%
)\text{\quad and\quad}\overline{K}_{p}(\mathbb{Z}^{\prime})\cong\overline
{K}_{p}(\mathbb{F}_{3})\text{.}%
\]
More precisely, there exists a Bott element $b_{K}$ in $\overline{K}%
_{p}(\mathbb{Z}^{\prime})$ mapping at the same time to a generator of
$\overline{K}_{p}(\mathbb{R})$ and to a generator of $\overline{K}%
_{p}(\mathbb{F}_{3})$.

We proceed in the same manner in order to explicate Bott elements in hermitian
$K$-theory, having almost the exact same properties as their namesakes in
algebraic $K$-theory. More precisely, we shall prove the following theorem:

\begin{theorem}
\label{positive Bott element in KQ}Let $p\geq8$ a $2$-power and $m=2p$. Then
the group ${}_{1}KQ_{p}(\mathbb{Z}^{\prime};\,\mathbb{Z}/m)$ is isomorphic to
$\mathbb{Z}/m\oplus\mathbb{Z}/m\oplus\mathbb{Z}/2$. There is a Bott element
$b^{+}$ in ${}_{1}KQ_{p}(\mathbb{Z}^{\prime};\,\mathbb{Z}/m)$ that maps at the
same time to a generator of {}$\mathbb{Z}/m$ in ${}_{1}KQ_{p}(\mathbb{F}%
_{3};\,\mathbb{Z}/m)\cong\mathbb{Z}/m\oplus\mathbb{Z}/2$ and to a generator of
${}_{1}KQ_{p}(\mathbb{R};\,\mathbb{Z}/m)$, viewed as a module\footnote{The
ring structure for $KQ$-theory with mod$\ 2^{\nu}$ coefficients is
well-defined if $\nu\geq4$.\label{footnote: KQ ring structure}} over
$_{1}KQ_{0}(\mathbb{R};\,\mathbb{Z}/m)$.
\end{theorem}

\noindent\textbf{Proof.} In the following proof, we are going to use the
results of \cite[Theorem 6.1]{BK} and \cite[Theorems 1.2, 1.5]{BKO}. In
\cite{BK}, it is shown that the square of hermitian $K$-theory completed
connective spectra
\[%
\begin{array}
[c]{ccc}%
{}_{1}\mathcal{KQ(}\mathbb{Z}^{\prime})_{\#}^{c} & \longrightarrow & {}%
_{1}\mathcal{KQ(}\mathbb{R})_{\#}^{c}\\
\downarrow &  & \downarrow\\
{}_{1}\mathcal{KQ(}\mathbb{F}_{3})_{\#}^{c} & \longrightarrow & {}%
_{1}\mathcal{KQ(}\mathbb{C})_{\#}^{c}%
\end{array}
\]
is homotopy cartesian (Recall that ${}_{1}\mathcal{KQ}(\mathbb{C})$ is just
$\mathcal{KO}$.) Reducing mod\ $m$ yields another homotopy cartesian square:
\[%
\begin{array}
[c]{ccc}%
_{1}\overline{\mathcal{KQ}}\mathcal{(}\mathbb{Z}^{\prime})^{c} &
\longrightarrow & _{1}\overline{\mathcal{KQ}}\mathcal{(}\mathbb{R})^{c}\\
\downarrow &  & \downarrow\\
_{1}\overline{\mathcal{KQ}}\mathcal{(}\mathbb{F}_{3})^{c} & \longrightarrow &
_{1}\overline{\mathcal{KQ}}\mathcal{(}\mathbb{C})^{c}%
\end{array}
\]
This in turn gives rise to a short exact sequence
\begin{equation}
0\longrightarrow{}_{1}\overline{KQ}_{p}(\mathbb{Z}^{\prime})\longrightarrow
{}_{1}\overline{KQ}_{p}(\mathbb{R})\oplus{}_{1}\overline{KQ}_{p}%
(\mathbb{F}_{3})\longrightarrow{}_{1}\overline{KQ}_{p}(\mathbb{C}%
)\longrightarrow0\text{,} \label{split ses for KQ(Z')}%
\end{equation}
which splits since ${}_{1}\overline{KQ}_{p}(\mathbb{R})$ is a direct sum of
two copies of ${}_{1}\overline{KQ}_{p}(\mathbb{C})$, say $G\oplus G$. The
first copy of $G$, say $G_{1}$, is generated by the image of $1$ under the
Bott isomorphism ${}_{1}\overline{KQ}_{0}(\mathbb{R})\cong{}_{1}\overline
{KQ}_{p}(\mathbb{R})$ (see Appendix $B$ in \cite{BKO}). The splitting is given
by the isomorphism between ${}_{1}\overline{KQ}_{p}(\mathbb{C})$ and the
second copy of $G$, say $G_{2}$. Therefore, we get an isomorphism
\[
{}_{1}\overline{KQ}_{p}(\mathbb{Z}^{\prime})\cong G_{1}\oplus{}_{1}%
\overline{KQ}_{p}(\mathbb{F}_{3})
\]

In order to finish the proof of the theorem, we need to compute ${}%
_{1}\overline{KQ}_{p}(\mathbb{F}_{3})$. By \cite{Friedlander}, there is a
Bockstein exact sequence
\begin{equation}
0\longrightarrow\mathbb{Z}/2\longrightarrow{}_{1}\overline{KQ}_{p}%
(\mathbb{F}_{3})\longrightarrow\mathbb{Z}/m\longrightarrow0\text{.}
\label{KQ-_8(F_3)}%
\end{equation}
In order to resolve this extension problem, consider the map
\[
{}_{1}KQ_{0}(\mathbb{F}_{3})/m=\mathbb{Z}/m\oplus\mathbb{Z}/2\longrightarrow
{}_{1}\overline{KQ}_{p}(\mathbb{F}_{3})
\]
given by cup-product with any element that maps to the generator of ${}%
_{1}\overline{KQ}_{p}(\mathbb{C})\cong\mathbb{Z}/m$ under the Brauer lift
${}_{1}\overline{KQ}_{p}(\mathbb{F}_{3})\rightarrow{}_{1}\overline{KQ}%
_{p}(\mathbb{C})$. This gives a splitting of the exact sequence
(\ref{KQ-_8(F_3)}), and therefore ${}_{1}\overline{KQ}_{p}(\mathbb{F}%
_{3})\cong\mathbb{Z}/m\oplus\mathbb{Z}/2$. $\hfill\Box\smallskip$

\begin{Remarks}
By considering the forgetful map from the hermitian $K$ sequence
(\ref{split ses for KQ(Z')}) to its algebraic $K$ counterpart, one sees that
the Bott element of ${}_{1}KQ_{p}(\mathbb{Z}^{\prime};\,\mathbb{Z}/m)$ maps to
the Bott element in the corresponding algebraic $K$-theory group under the map
induced by the forgetful functor. Moreover, all the results for $\mathbb{F}%
_{3}$ in the above also hold for any finite field $\mathbb{F}_{t}$ with $t$
elements, provided $t\equiv\pm3\;(\mathrm{mod}$\ $8)$.
\end{Remarks}

\section{Proof of the periodicity theorem for totally real $2$-regular number
fields\label{Proof of the periodicity theorem for totally real 2-regular number fields}%
}

Let $A$ be the ring of $2$-integers in a totally real $2$-regular number field
$F$ with $r$ real embeddings. In \cite{BKO}, we proved that the square of
hermitian $K$-theory $2$-completed connective spectra
\[%
\begin{array}
[c]{ccc}%
_{\varepsilon}\mathcal{KQ(}A)_{\#}^{c} & \longrightarrow & \bigvee^{r}%
{}_{\varepsilon}\mathcal{KQ(}\mathbb{R})_{\#}^{c}\\
\downarrow &  & \downarrow\\
_{\varepsilon}\mathcal{KQ(}\mathbb{F}_{t})_{\#}^{c} & \longrightarrow &
\bigvee^{r}{}_{\varepsilon}\mathcal{KQ(}\mathbb{C})_{\#}^{c}%
\end{array}
\]
is homotopy cartesian (with $t$ a carefully chosen odd prime and where $_{\#}$
denotes $2$-adic completion). Therefore, the mod\ $2^{\nu}$ reduction of this
square, namely
\[%
\begin{array}
[c]{ccc}%
_{\varepsilon}\overline{\mathcal{KQ}}\mathcal{(}A)^{c} & \longrightarrow &
\bigvee^{r}{}_{\varepsilon}\overline{\mathcal{KQ}}\mathcal{(}\mathbb{R})^{c}\\
\downarrow &  & \downarrow\\
_{\varepsilon}\overline{\mathcal{KQ}}\mathcal{(}\mathbb{F}_{t})^{c} &
\longrightarrow & \bigvee^{r}{}_{\varepsilon}\overline{\mathcal{KQ}%
}\mathcal{(}\mathbb{C})^{c}%
\end{array}
\]
is also homotopy cartesian, since ${}_{\varepsilon}KQ_{-1}(A)=0$ by Lemmas
3.11 and 3.12 in \cite{BKO}. Using this square, we deduce an enhanced version
of our periodicity theorem.

\begin{theorem}
\label{2-regular main}For $n\geq0$ and $p=\sup\{8,2^{\nu-1}\}$ for $\nu\geq1$,
taking cup-product with the positive Bott element in $_{1}KQ_{p}%
(\mathbb{Z}^{\prime};\,\mathbb{Z}/2^{\nu})$ induces an isomorphism
\[
{}_{\varepsilon}KQ_{n}(A;\,\mathbb{Z}/2^{\nu})\cong{}_{\varepsilon}%
KQ_{n+p}(A;\,\mathbb{Z}/2^{\nu}).
\]

\end{theorem}

\noindent\textbf{Proof.} Cup-product with the Bott element in ${}_{1}%
KQ_{p}(\mathbb{Z}^{\prime};\,\mathbb{Z}/2^{\nu})$ induces an isomorphism of
${}_{\varepsilon}\overline{KQ}$-groups for the rings $\mathbb{F}_{t}$,
$\mathbb{R}$ and $\mathbb{C}$, where $\mathbb{F}_{t}$ is the finite field with
$t$ elements. This is due to $\mathcal{KO}$-localness of the corresponding
hermitian $K$-theory spectra, and an induction on the order of the coefficient
group based on the five lemma applied to the Bockstein exact sequence.
Therefore, the result follows from the five lemma together with the homotopy
cartesian square above. $\hfill\Box\smallskip$

\begin{remark}
The isomorphism for $n=0$ reflects the fact that $2$-regularity implies that
there is no nontrivial $2$-torsion in the Picard group of $A$.
\end{remark}

We note that the number $\nu$ was related to the choice of $t\equiv
\pm3\;(\mathrm{mod}$\ $8)$ in Theorem \ref{positive Bott element in KQ} (for
$A=\mathbb{Z}^{\prime}$). However, the number $t$ that makes the diagrams
above homotopy cartesian (for $A$ totally real $2$-regular) is different in general.

Therefore, we can improve the previous result by replacing $m=2^{\nu}$ by $M$,
which is the number $m$ multiplied by the $2$-primary factor $m^{\prime
}=\left(  \frac{t^{2}-1}{8}\right)  _{2}$ of $(t^{2}-1)/8$ (compare with Lemma
2.9 in \cite{BKO}). More precisely, we have the following proposition.

\begin{proposition}
\label{2-regular Bott element}Let $p=\sup\{8,2^{\nu-1}\}$ and consider the
canonically induced map
\[
{}_{1}KQ_{p}(\mathbb{Z}^{\prime};\,\mathbb{Z}/2^{\nu})\longrightarrow{}%
_{1}KQ_{p}(A;\,\mathbb{Z}/2^{\nu})\text{.}%
\]
The image of the Bott element $b$ in $_{1}KQ_{p}(A;\,\mathbb{Z}/2^{\nu})$ is
the reduction $\operatorname{mod}\ 2^{\nu}$ of a class $\overline{b}$
$\operatorname{mod}\ M$, with $M=m\cdot\left(  \frac{t^{2}-1}{8}\right)  _{2}$
as defined above.
\end{proposition}

\noindent\textbf{Proof.} For brevity, we use the above notation, whereby
$m=2^{\nu}$ and $M=mm^{\prime}$. As in the case of $\mathbb{Z}^{\prime}$
considered in Section \ref{section: Bott elements}, we can write the following
diagram of exact sequences (where $KQ={}_{1}KQ$):%
\[%
\begin{array}
[c]{ccc}
&  & 0\\
&  & \downarrow\\
& KQ_{p}(\mathbb{R}^{r};\,\mathbb{Z}/m^{\prime})\oplus KQ_{p}(\mathbb{F}%
_{t};\,\mathbb{Z}/m^{\prime}) & \rightarrow KQ_{p}(\mathbb{C}^{r}%
;\,\mathbb{Z}/m^{\prime})\rightarrow0\\
& \downarrow & \downarrow\\
KQ_{p}(A;\,\mathbb{Z}/M)\rightarrow & KQ_{p}(\mathbb{R}^{r};\,\mathbb{Z}%
/M)\oplus KQ_{p}(\mathbb{F}_{t};\,\mathbb{Z}/M) & \rightarrow KQ_{p}%
(\mathbb{C}^{r};\,\mathbb{Z}/M)\rightarrow0\\
\downarrow & \downarrow & \downarrow\\
KQ_{p}(A;\,\mathbb{Z}/m)\rightarrow & KQ_{p}(\mathbb{R}^{r};\,\mathbb{Z}%
/m)\oplus KQ_{p}(\mathbb{F}_{t};\,\mathbb{Z}/m) & \rightarrow KQ_{p}%
(\mathbb{C}^{r};\,\mathbb{Z}/m)\rightarrow0\\
& \downarrow & \downarrow\\
& 0 & 0
\end{array}
\]
Chasing in this diagram shows the reduction map ${}KQ_{p}(A;\,\mathbb{Z}%
/M)\rightarrow{}KQ_{p}(A;\,\mathbb{Z}/m)$ is surjective. Therefore, the Bott
element $b$ in $KQ_{p}(A;\,\mathbb{Z}/m)$ is the reduction $\operatorname{mod}%
\ m$ of a class $\overline{b}$ $\operatorname{mod}\ M$ which we shall call an
\emph{exotic Bott element}\ (we do not claim, however, that $\overline{b}$ is
unique).\hfill$\Box$\smallskip

\begin{theorem}
Let $\overline{b}$ be an exotic\ Bott element in the group $_{1}%
KQ_{p}(A;\,\mathbb{Z}/M)$ defined above. Then cup-product with $\overline{b}$
induces an isomorphism%
\[
\overline{\beta}:{}_{\varepsilon}KQ_{n}(A;\,\mathbb{Z}/M^{\prime}%
)\overset{\cong}{\longrightarrow}{}_{\varepsilon}KQ_{n+p}(A;\,\mathbb{Z}%
/M^{\prime})
\]
for every $n\geq0$ and divisor $M^{\prime}$ of $M$.
\end{theorem}

\noindent\textbf{Proof.} We just copy the proof of Theorem
\ref{2-regular main}, using the five lemma, since this periodicity statement
holds for the rings $\mathbb{R},\mathbb{C}$ and $\mathbb{F}_{t}$ (see the
independent lemma below for the field $\mathbb{F}_{t}$).\hfill$\Box$\smallskip

\begin{lemma}
Let $\mathbb{F}_{t}$ be a finite field with $t$ elements and let
$p=\sup\{8,2^{\nu-1}\}$ and $m=2^{\nu}$. Then the image of the Bott element by
the canonical map%
\[
{}_{1}KQ_{p}(\mathbb{Z}^{\prime};\,\mathbb{Z}/m)\longrightarrow{}_{1}%
KQ_{p}(\mathbb{F}_{t};\,\mathbb{Z}/m)
\]
is the reduction $\operatorname{mod}\ m$ of a class $\operatorname{mod}%
\ M^{\prime}$ (with $m\,|\,M^{\prime}$) if and only if
\[
(M^{\prime})_{2}\leq m\cdot((t^{2}-1)/8)_{2}\text{,}%
\]
where $(i)_{2}$ is the $2$-primary part of $i$.
\end{lemma}

\noindent\textbf{Proof.} We look at the following commutative diagram, with
exact rows:%
\[%
\begin{array}
[c]{cccccc}%
\longrightarrow & {}_{1}KQ_{p}(\mathbb{F}_{t};\,\mathbb{Z}/M^{\prime}) &
\longrightarrow & _{1}KQ_{p-1}(\mathbb{F}_{t}) & \overset{\cdot M^{\prime}%
}{\longrightarrow} & _{1}KQ_{p-1}(\mathbb{F}_{t})\\
& \downarrow\alpha_{M^{\prime}} &  & \downarrow\cdot M^{\prime}2^{-\nu} &  &
\downarrow\mathrm{id}\\
\longrightarrow & {}_{1}KQ_{p}(\mathbb{F}_{t};\,\mathbb{Z}/2^{\nu}) &
\longrightarrow & _{1}KQ_{p-1}(\mathbb{F}_{t}) & \overset{\cdot2^{\nu}%
}{\longrightarrow} & _{1}KQ_{p-1}(\mathbb{F}_{t})
\end{array}
\]
The Bott element in the group $_{1}KQ_{p}(\mathbb{F}_{t};\,\mathbb{Z}/2^{\nu
})$ maps nontrivially into the group $_{1}KQ_{p-1}(\mathbb{F}_{t})$ and its
image is divisible by $M^{\prime}2^{-\nu}$. On the other hand, we know by
\cite{Friedlander} that $_{1}KQ_{p-1}(\mathbb{F}_{t})$ is cyclic of order $M$,
where $M$ is the $2$-primary part of $(t^{p/2}-1)$, which is also the
$2$-primary part of $(t^{2}-1)\cdot p/4$ by Lemma 2.7 in \cite{BKO}. This
number is also the $2$-primary part of $2^{\nu}\cdot(t^{2}-1)/8$. Therefore, a
simple diagram chase shows that $\alpha_{M^{\prime}}$ is surjective if and
only if $M^{\prime}\,|\,M$. \hfill$\Box$\smallskip

Now let us consider a nonzero prime ideal $\mathfrak{p}$ in $A$, and the
quotient field $A/\mathfrak{p}$. There is a commutative diagram%
\[%
\begin{array}
[c]{ccc}%
\mathbb{Z}^{\prime} & \longrightarrow & A/\mathfrak{p}\\
\downarrow &  & \downarrow\\
A & \longrightarrow & A/\mathfrak{p}%
\end{array}
\]
where the right vertical arrow is the identity map. Since the Bott element in
the $KQ$-group ${}_{1}KQ_{p}(A/\mathfrak{p};\,\mathbb{Z}/2^{^{\nu}})$ is the
reduction $\operatorname{mod}2^{\nu}$ of a class $\operatorname{mod}\ M$,
where $M$ is a power of $2$, we have an isomorphism
\[
_{1}KQ_{p}(A/\mathfrak{p;}\,\mathbb{Z}/M)\cong\mathbb{Z}/M\oplus\mathbb{Z}/2
\]
according to the computations of the $KQ$-theory of finite fields in
\cite{Friedlander} and Section \ref{section: Bott elements}. It follows that
there is a periodicity isomorphism%
\[
_{\varepsilon}KQ_{n}(A/\mathfrak{p};\,\mathbb{Z}/M^{\prime})\cong%
{}_{\varepsilon}KQ_{n+p}(A/\mathfrak{p};\,\mathbb{Z}/M^{\prime})
\]
for $n\geq0$ and any $M^{\prime}\,|\,M$, given by the cup-product with an
exotic Bott element. For the next two results, we recall from
\cite[Proposition 2.1]{BKO} that $F$ contains a unique dyadic prime (that is,
prime ideal lying over the rational prime $(2)$). For any set $S$ of
valuations in $F$ including the dyadic valuation and the infinite ones, we
define $A_{S}$ to consist of the elements in $F$ whose valuations not in $S$
are non-negative. Thus, when $S$ is finite, $A_{S}$ is just the ring of
$S$-integers. When $S$ comprises only the dyadic valuation and the infinite
ones, $A_{S}=A$; while, when $S$ comprises all valuations, $A_{S}=F$.

\begin{theorem}
\label{KQAS}Let $p=\sup\{8,2^{\nu-1}\}$ for $\nu\geq1$. Then, for $n>0$,
cup-product with an exotic Bott element in $_{1}KQ_{p}(A;\,\mathbb{Z}/M)$
induces an isomorphism
\[
\overline{\beta}:{}_{\varepsilon}KQ_{n}(A_{S};\,\mathbb{Z}/M^{\prime}%
)\overset{\cong}{\longrightarrow}{}_{\varepsilon}KQ_{n+p}(A_{S};\,\mathbb{Z}%
/M^{\prime})
\]
for any $M^{\prime}$ such that $2\,|\,M^{\prime}\,|\,M$.
\end{theorem}

\noindent\textbf{Proof.} We use the homotopy fibration
\[
\bigvee{}_{\varepsilon}\mathcal{U}(A/\mathfrak{p})\longrightarrow
{}_{\varepsilon}\mathcal{KQ}(A)\longrightarrow{}_{\varepsilon}\mathcal{KQ}%
(A_{S})
\]
noted in \cite{Hornbostel}, where $\mathfrak{p}$ runs through all nonzero
prime ideals in $S$. For the corresponding mod\ $M^{\prime}$ reductions
(indicated as usual by an overbar) where $M^{\prime}\,|\,M$, there is a
homotopy fibration
\[
\bigvee{}_{\varepsilon}\overline{\mathcal{U}}(A/\mathfrak{p})\longrightarrow
{}_{\varepsilon}\overline{KQ}(A)\longrightarrow{}_{\varepsilon}\overline
{KQ}(A_{S}).
\]
The maps in this fibration are compatible with cup-products with elements of
$\overline{KQ}_{\ast}(A)$. The $\overline{U}$-theory spectra of finite fields
are $\mathcal{KO}$-local as we showed more precisely above (this is a
consequence of the same property for the $\overline{K}$ and $\overline{KQ}%
$-theories). From these facts, the five lemma implies the Bott periodicity
isomorphism
\[
{}_{\varepsilon}KQ_{n}(A_{S};\,\mathbb{Z}/M^{\prime})\cong{}_{\varepsilon
}KQ_{n+p}(A_{S};\,\mathbb{Z}/M^{\prime})
\]
for $n>0$, given by the cup-product with an exotic Bott element. \hfill$\Box
$\smallskip

\begin{theorem}
\label{KQXAS}Let $A_{S}$ be as before, and let $\overline{b}$ be an
exotic\ Bott element in the group $_{1}KQ_{p}(A;\,\mathbb{Z}/M)$. Then, for
any connected CW-complex $X$, cup-product with $\overline{b}$ induces an
isomorphism%
\[
\overline{\beta}:{}_{\varepsilon}KQ_{X}(A_{S};\,\mathbb{Z}/M^{\prime})\cong%
{}_{\varepsilon}KQ_{X+p}(A_{S};\,\mathbb{Z}/M^{\prime})
\]
for any $M^{\prime}$ such that $2\,|\,M^{\prime}\,|\,M$. Moreover, when
$A_{S}=A$, the previous isomorphism holds for any CW-complex, not necessarily
connected\textsf{.}
\end{theorem}

\noindent\textbf{Proof.} By Theorem \ref{KQAS}, the Bott map $\overline{\beta
}$ is an isomorphism when $X$ is a sphere $S^{n}$ for $n\geq1$. According to
general facts about representable cohomology theories \cite{Boardman}, it
follows that $\overline{\beta}$ is also an isomorphism if $X$ is a connected
CW-complex, (finite or infinite, thanks to Milnor's $\lim^{1}$ exact
sequence). If $A_{S}=A$, then the Bott map $\overline{\beta}$ is also an
isomorphism when $X=$ $S^{0}$. Therefore, the previous isomorphism holds also
for not necessarily connected CW-complexes.\hfill$\Box$\smallskip

\section{Higher KSC-theories\label{Higher KSC-theories}}

The useful concept of topological $K$-theory based upon self conjugate vector
bundles $KSC$ was introduced by Anderson \cite{Anderson} and Green
\cite{Green}. In \cite[p. 281]{KaroubiAnnals2}, for a ring $A$ with
involution, the spectrum $\mathcal{KSC}(A)$ was defined as the homotopy fiber
of $1-\tau$, where $\tau$ is the duality functor in algebraic $K$-theory
\[
\tau:\mathcal{K}(A)\longrightarrow\mathcal{K}(A)\text{.}%
\]
The importance of $KSC$-theory becomes evident from the homotopy fibration
\cite[p.~282]{KaroubiAnnals2}
\[
\mathcal{KSC}(A)\longrightarrow\Omega{}_{\varepsilon}\mathcal{KQ}%
(A)\overset{\sigma^{(2)}}{\longrightarrow}\Omega^{-1}{}_{-\varepsilon
}\mathcal{KQ}(A)\text{,}%
\]
which implies a long exact sequence (for legibility we omit the ring $A$ in
the notation)%
\[
\cdots\rightarrow{}_{\varepsilon}KQ_{n+2}\overset{s^{(2)}}{\rightarrow}%
{}_{-\varepsilon}KQ_{n}\rightarrow KSC_{n}\rightarrow{}_{\varepsilon}%
KQ_{n+1}{}\overset{s^{(2)}}{\rightarrow}{}_{-\varepsilon}KQ_{n-1}%
\rightarrow\cdots.
\]
The morphism $s^{(2)}$ between the $KQ$-groups is the periodicity map made
explicit in \cite{KaroubiAnnals2}. It is defined by taking cup-product with a
generator of the free part of the group
\[
_{-1}KQ_{-2}(\mathbb{Z}^{\prime})\cong{}_{1}W_{0}(\mathbb{Z}^{\prime}%
)\cong\mathbb{Z}\oplus\mathbb{Z}/2\text{.}%
\]
(Recall our assumption that $1/2\in A$.) We should note that this cup-product
induces a morphism between cohomology theories, and thence the associated
$KQ$-spectra and $\overline{KQ}$-spectra, according to Brown's
representability theorem.

It turns out that the $KSC$-groups measure the failure of negative Bott
periodicity for the $KQ$-groups. To keep track of the degree shift we let
$KSC^{(2)}$ (resp.~$\mathcal{KSC}^{(2)}$) denote the $KSC$-groups (resp.~$KSC$-spectrum).

There exist higher analogs of this spectrum corresponding to degree shifts by
$4$, $8$ and higher $2$-powers.

The next version, denoted\footnote{\textsl{A priori},
this theory depends on $\varepsilon$. A proof of this statement may be found
in Lemma \ref{epsilon dependance} below. For $KSC$ -theory with coefficients,
we can also argue by contradiction as in Lemma \ref{decomposition}.} by
${}_{\varepsilon}\mathcal{KSC}^{(4)}\mathcal{(}A\mathcal{)}$, is the homotopy
fiber of the composite map%
\[
\sigma^{(4)}:\Omega{}_{\varepsilon}\mathcal{KQ}(A)\overset{\sigma^{(2)}%
}{\longrightarrow}\Omega^{-1}{}_{-\varepsilon}\mathcal{KQ}(A)\overset
{\Omega^{(-2)}\sigma^{(2)}}{\longrightarrow}\Omega^{-3}{}_{\varepsilon
}\mathcal{KQ}(A)\text{.}%
\]

\begin{proposition}
\label{KSCiteration}There exists a homotopy fibration of spectra
\[
{}_{\varepsilon}\mathcal{KSC}^{(2)}\mathcal{(}A\mathcal{)}\longrightarrow
{}_{\varepsilon}\mathcal{KSC}^{(4)}\mathcal{(}A\mathcal{)}\longrightarrow
\Omega^{-2}({}_{\varepsilon}\mathcal{KSC}^{(2)}\mathcal{(}A\mathcal{))}%
\]
and a long exact sequence (we again omit the ring $A$ for convenience)
\[
\cdots\rightarrow{}_{\varepsilon}KQ_{n+2}\overset{s^{(4)}}{\rightarrow}%
{}_{\varepsilon}KQ_{n-2}\rightarrow{}_{\varepsilon}KSC_{n}^{(4)}\rightarrow
{}_{\varepsilon}KQ_{n+1}\overset{s^{(4)}}{\rightarrow}{}_{\varepsilon}%
KQ_{n-3}\rightarrow\cdots\text{.}%
\]

\end{proposition}

\noindent\textbf{Proof.} This is just the observation that for two composable
maps $u$ and $v$, there is a homotopy fibration $\mathcal{F}(u)\longrightarrow
\mathcal{F}(v\circ u)\longrightarrow\mathcal{F}(v)$, where $\mathcal{F}(f)$
denotes the homotopy fiber of some map $f$. $\hfill\Box\smallskip$

Iterating, for $r>4$ a $2$-power, we proceed
similarly and define ${}_{\varepsilon}\mathcal{KSC}^{(r)}\mathcal{(}%
A\mathcal{)}$ as the homotopy fiber of the map
\[
\sigma^{(r)}:\Omega{}_{\varepsilon}\mathcal{KQ}(A)\ \longrightarrow
\Omega^{-r+1}{}_{\varepsilon}\mathcal{KQ}(A)
\]
where $\sigma^{(r)}=\Omega^{-r/2}\sigma^{(r/2)}\circ\sigma^{(r/2)}$. In the
other direction, if we allow the convention ${}_{\varepsilon}\mathcal{KSC}%
^{(1)}(A)={\Omega K}(A)$, then from the original definition
of $\mathcal{KSC}$ above the following also holds for $r=2$.

\begin{proposition}
\label{KSC/KSC}For a $2$-power $r\geq2$, there is a homotopy fibration of
spectra
\[
{}_{\varepsilon}\mathcal{KSC}^{(r/2)}\mathcal{(}A\mathcal{)}\longrightarrow
{}_{\varepsilon}\mathcal{KSC}^{(r)}\mathcal{(}A\mathcal{)}\longrightarrow
\Omega^{-r/2}{}_{\varepsilon}\mathcal{KSC}^{(r/2)}(A)
\]
and an associated long exact sequence
\[
\cdots\rightarrow{}_{\varepsilon}KQ_{n+2}\overset{s^{(r)}}{\rightarrow}%
{}_{\varepsilon^{\prime}}KQ_{n+2-r}\rightarrow{}_{\varepsilon}KSC_{n}%
^{(r)}\rightarrow{}_{\varepsilon}KQ_{n+1}\overset{s^{(r)}}{\rightarrow}%
{}_{\varepsilon^{\prime}}KQ_{n+1-r}\rightarrow\cdots
\]
where $\varepsilon^{\prime}=-\varepsilon$ if $r=2$ and $\varepsilon^{\prime
}=\varepsilon$ if $r>2$.$\hfill\Box$
\end{proposition}

Finally, we show that the higher $KSC$-theories depends on the sign of
symmetry $\varepsilon$.

\begin{lemma}
\label{epsilon dependance}Let $F$ be a finite field of characteristic $\neq2$.
Then the group ${}_{1}KSC_{1}^{(4)}(F)$ is isomorphic to $\mathbb{Z}/2$, while
${}_{-1}KSC_{1}^{(4)}(F)=0$.
\end{lemma}

\noindent\textbf{Proof.\ }Let us drop the field $F$ for notational
convenience. Then the group ${}_{-1}KSC_{1}^{(4)}$ fits into the exact
sequence%
\[
{}_{-1}KQ_{3}\longrightarrow{}_{-1}KQ_{-1}\longrightarrow{}_{-1}KSC_{1}%
^{(4)}\longrightarrow{}_{-1}KQ_{2}\longrightarrow{}_{-1}KQ_{-2}\text{.}%
\]
We have ${}_{-1}KQ_{-1}=0$ by the same argument used in the proof of Lemma
3.11 in \cite{BKO}, where we should replace $R_{F}$ by $F$. We also have
${}_{-1}KQ_{2}=0$ by a result of Friedlander \cite{Friedlander}. Therefore,
${}_{-1}KSC_{1}^{(4)}(F)=0$.

On the other hand, the group ${}_{1}KSC_{1}^{(4)}$ fits into the exact
sequence%
\[
{}_{1}KQ_{3}\longrightarrow{}_{1}KQ_{-1}\longrightarrow{}_{1}KSC_{1}%
^{(4)}\longrightarrow{}_{1}KQ_{2}\overset{\alpha}{\longrightarrow}{}%
_{1}KQ_{-2}\text{.}%
\]
For the same reason as above, we have ${}_{1}KQ_{-1}=0$. We also have ${}%
_{1}KQ_{2}=\mathbb{Z}/2$ by an analogous result of Friedlander
\cite{Friedlander}. The periodicity map $\alpha$ can be factored though the
group ${}_{-1}KQ_{0}$ which is isomorphic to $\mathbb{Z}$. Therefore,
\textrm{Ker}($\alpha)=\mathbb{Z}/2$ and the group ${}_{1}KSC^{(4)}$ is
isomorphic to $\mathbb{Z}/2$.$\hfill\Box\smallskip$

We can mimic the previous definitions by taking spectra or groups
\textrm{mod\ }$m$, where $m$ is related to $p$ according to our convention
\ref{convention: (p,p')}. In that case, we shall write $\overline
{\mathcal{KQ}}$ instead of $\mathcal{KQ}$, $\overline{KSC}$ instead of $KSC$,
\textsl{etc.}

\begin{proposition}
\label{KSCperiodicity}Let $d$ be the number defined in the Introduction
(\textsl{i.e.}~the starting point of periodicity for the $\overline{K}%
$-groups). Then for any $2$-power $r\geq2$, the
positive Bott map%
\[
\sigma:{}_{\varepsilon}\overline{KSC}_{n}^{(r)}(A)\rightarrow{}_{\varepsilon
}\overline{KSC}_{n+p}^{(r)}(A)
\]
is an isomorphism if $n\geq d+r-2$.
\end{proposition}

\noindent\textbf{Proof.} We argue by iteration on the
$2$-power $r$, using the following diagram of exact sequences from
(\ref{KSC/KSC}):%
\begin{equation}%
\begin{array}
[c]{ccccc}%
\overline{KSC}_{n+1-r/2}^{(r/2)}{\tiny \rightarrow} & \overline{KSC}%
_{n}^{(r/2)}{\tiny \rightarrow} & \overline{KSC}_{n}^{(r)}{\tiny \rightarrow}
& \overline{KSC}_{n-r/2}^{(r/2)}{\tiny \rightarrow} & \overline{KSC}%
_{n-1}^{(r/2)}\\
{\tiny \downarrow} & {\tiny \downarrow} & {\tiny \downarrow} &
{\tiny \downarrow} & {\tiny \downarrow}\\
\overline{KSC}_{n++p+1-r/2}^{(r/2)}{\tiny \rightarrow} & \overline{KSC}%
_{n+p}^{(r/2)}{\tiny \rightarrow} & \overline{KSC}_{n+p}^{(r)}%
{\tiny \rightarrow} & \overline{KSC}_{n+p-r/2}^{(r/2)}{\tiny \rightarrow} &
\overline{KSC}_{n+p-1}^{(r/2)}%
\end{array}
\label{KSC/KSC exact seq}%
\end{equation}
Commutativity of this diagram follows from the fact that the vertical maps are
induced by cup-product with the positive Bott element in $KQ$-theory as
constructed in Section \ref{section: Bott elements}, and that all maps are
$KQ$-module maps. By induction, we know that the vertical maps, with the
possible exception of the middle one, are isomorphisms if $n-r/2\geq d+r/2-2$,
that is $n\geq d+r-2$. We conclude thanks to the five lemma. \hfill$\Box
$\smallskip

\begin{remark}
A variant of this proposition is to consider a parameter space $X$ instead of
a sphere $S^{n}$. More precisely, by the method of proof of Theorem
\ref{KQXAS}, the Bott map%
\[
\sigma:{}_{\varepsilon}\overline{KSC}_{X}^{(r)}(A)\rightarrow{}_{\varepsilon
}\overline{KSC}_{X+p}^{(r)}(A)
\]
is an isomorphism if the space $X$ is $(d+r-3)$-connected.
\end{remark}

Next we consider the failure of positive $p$-periodicity in hermitian
$\overline{K}$-theory. This is encoded in the homotopy fiber of the
periodicity map given by the cup-product with the positive Bott element
\[
_{\varepsilon}\overline{\mathcal{KQ}}(A)\longrightarrow\Omega^{p}%
{}_{\varepsilon}\overline{\mathcal{KQ}}(A)\text{,}%
\]
which we shall denote by $\mathcal{P{}}_{\varepsilon}\overline{\mathcal{KQ}%
}(A)$. In the same way, we denote by $P\overline{\mathcal{K}}(A)$ the homotopy
fiber of the periodicity map in $K$-theory
\[
\overline{\mathcal{K}}(A)\longrightarrow\Omega^{p}{}\overline{\mathcal{K}%
}(A)\text{.}%
\]
According to our general assumptions, the homotopy groups $P\overline{K_{n}%
}(A)$ of $\mathcal{P}\overline{\mathcal{K}}(A)$ vanish if $n\geq d$. On the
other hand, we can introduce, cf.~\cite{KaroubiAnnals2}, the homotopy fibers
$\mathcal{P}_{\varepsilon}\overline{\mathcal{U}}(A)$ and $\mathcal{P}%
_{-\varepsilon}\overline{\mathcal{V}}(A)$ of the hyperbolic and forgetful maps
$\mathcal{P}\overline{\mathcal{K}}(A)\rightarrow\mathcal{P}_{\varepsilon
}\overline{\mathcal{KQ}}(A)$ and $\mathcal{P}_{-\varepsilon}\overline
{\mathcal{KQ}}(A)\rightarrow\mathcal{P}\overline{\mathcal{K}}(A)$, respectively.

\begin{proposition}
There is a homotopy equivalence%
\[
\mathcal{P{}}_{-\varepsilon}\overline{\mathcal{V}}(A)\simeq\Omega
(\mathcal{P{}}_{\varepsilon}\overline{\mathcal{U}}(A))\text{.}%
\]
Moreover, the composition%
\[
\Omega^{2}(\mathcal{P}_{\varepsilon}\overline{\mathcal{KQ}}(A))\longrightarrow
\Omega(\mathcal{P{}}_{\varepsilon}\overline{\mathcal{U}}(A))\simeq
\mathcal{P{}}_{-\varepsilon}\overline{\mathcal{V}}(A)\longrightarrow
\mathcal{P}_{-\varepsilon}\overline{\mathcal{KQ}}(A)
\]
is induced by cup-product with the negative Bott element in the group ${}%
_{-1}KQ_{-2}(\mathbb{Z}^{\prime})$.
\end{proposition}

\noindent\textbf{Proof.} The fundamental theorem of hermitian $K$-theory
\cite{KaroubiAnnals2} exhibits an explicit homotopy equivalence (given both
ways) between the spectra $\mathcal{{}}$$_{-\varepsilon}\mathcal{V}(A)$ and
$\Omega\mathcal{{}}_{\varepsilon}\mathcal{U}(A)$.

Moreover, with the notation $\mathcal{KQ}(\mathbb{Z}^{\prime})={}%
_{1}\mathcal{KQ}(\mathbb{Z}^{\prime})\vee{}_{-1}\mathcal{KQ(}\mathbb{Z}%
^{\prime})$, these maps are $\mathcal{KQ}(\mathbb{Z}^{\prime})$-module maps.
Therefore, the reduction \textrm{mod\ }$m$ of these spectra is also a homotopy
equivalence. Since the composition
\[
\Omega^{2}(\mathcal{P}_{\varepsilon}\overline{\mathcal{KQ}}(A))\longrightarrow
\Omega(\mathcal{P{}}_{\varepsilon}\overline{\mathcal{U}}(A))\simeq
\mathcal{P{}}_{-\varepsilon}\overline{\mathcal{V}}(A)\longrightarrow
\mathcal{P}_{-\varepsilon}\overline{\mathcal{KQ}}(A)
\]
is a $\mathcal{KQ}(\mathbb{Z}^{\prime})$-module map, it is defined by the
cup-product with the negative Bott element, as proved in \cite{KaroubiAnnals2}%
. $\hfill\Box\smallskip$

\begin{lemma}
\label{negative perioidicity}If $P\overline{K}_{n}(A)=0$ for $n\geq d$, then
the negative Bott map
\[
P{}_{\varepsilon}\overline{KQ}_{n+2}(A)\rightarrow P\mathcal{{}}%
_{-\varepsilon}\overline{KQ}_{n}(A)
\]
is an isomorphism for $n\geq d$ and is a monomorphism for $n=d-1$.
\end{lemma}

\noindent\textbf{Proof.\ }This follows from the diagram ($A$ omitted) with
exact rows%
\[%
\begin{array}
[c]{ccccccc}%
P\overline{K}_{n+2}\longrightarrow & P{}_{\varepsilon}\overline{KQ}_{n+2} &
\longrightarrow & P\mathcal{{}}_{\varepsilon}\bar{U}_{n+1} & \longrightarrow &
P\overline{K}_{n+1} & \\
&  &  & \downarrow\cong &  &  & \\
& P\overline{K}_{n+1} & \longrightarrow & P\mathcal{{}}_{-\varepsilon}\bar
{V}_{n} & \longrightarrow & P\mathcal{{}}_{-\varepsilon}\overline{KQ}_{n} &
\longrightarrow P\overline{K}_{n}%
\end{array}
\]
the vertical isomorphism being a consequence of the fundamental theorem in
hermitian $K$-theory.$\hfill\Box\smallskip$

Iteration of this Bott map induces a further isomorphism
\[
P{}_{\varepsilon}\overline{KQ}_{n+4}(A)\overset{\cong}{\longrightarrow}%
P{}_{\varepsilon}\overline{KQ}_{n}(A)\text{.}%
\]
The classical induction method \cite[(3.5)]{BK}, adapted to this case, enables
us to prove the following theorem. We recall that the overbar over the $KQ$
indicates reduction $\operatorname{mod}m$, where $m$ was defined in the
Introduction. Strictly speaking, one has to take $m\geq16$, in the theorem, so
that $\overline{KQ}_{\ast}(\mathbb{Z}^{\prime})$ is an associative ring ---
see Footnote \ref{footnote: KQ ring structure}. However, if $m<16$, we can
consider all these groups as modules over $KQ_{\ast}(\mathbb{Z}^{\prime
};\,\mathbb{Z}/16)$ and the Bott map still makes sense.

\begin{theorem}
\label{KQperioidcity}Assume that $P\overline{K}_{n}(A)=0$ for $n\geq d$.
Assume moreover that $P{}_{\varepsilon}\overline{KQ}_{n}(A)=0$ for
$\varepsilon=\pm1$ and for $n=d$ and $d+1$. Then the cup-product with the Bott
element in $_{1}\overline{KQ}_{p}(\mathbb{Z}^{\prime})$ induces a morphism
\[
\beta_{n}:{}_{\varepsilon}\overline{KQ}_{n}(A)\overset{\cong}{\longrightarrow
}{}_{\varepsilon}\overline{KQ}_{n+p}(A)
\]
which is an isomorphism for $n\geq d+1$ and a monomorphism for $n=d$.
\end{theorem}

\textbf{Proof.} The $2$-periodicity of the $P\overline{KQ}$-groups shown above
(with a change of symmetry) implies that $P{}_{\varepsilon}\overline{KQ}%
_{n}(A)=0$ for $n\geq d$ and $\varepsilon=\pm1$. From the exact sequence
\[
P{}_{\varepsilon}\overline{KQ}_{n}(A)\longrightarrow{}_{\varepsilon}%
\overline{KQ}_{n}(A)\overset{\beta_{n}}{\longrightarrow}{}_{\varepsilon
}\overline{KQ}_{n+p}(A)\longrightarrow P{}_{\varepsilon}\overline{KQ}%
_{n-1}(A),
\]
we deduce the required isomorphism (starting from $n=d+1$) and a monomorphism
for $n=d$.$\hfill\Box\smallskip$

Unfortunately, this strategy is not efficient to establish Bott periodicity
because the starting point of the induction is not always valid (see the end
of Section \ref{section: proof of main} for counterexamples and Section
\ref{Proof of the periodicity theorem for totally real 2-regular number fields}
for examples). Therefore, we are going to take another approach towards Bott
periodicity. From now on, we often assume implicitly that the $\overline{K}%
$-groups are periodic starting in degree $d$. More precisely, the Bott map%
\[
\overline{K}_{n}(A)\longrightarrow\overline{K}_{n+p}(A)
\]
is an isomorphism for $n\geq d$, which implies that $P\overline{K}_{n}(A)=0$
in the same range. Our aim is to prove a similar periodicity assertion for
$KQ$-theory, as we announced in the Introduction.

Let us investigate in detail the composition of the two \textquotedblleft
opposite\textquotedblright\ periodicity maps
\[
{}_{\varepsilon}\overline{\mathcal{KQ}}(A)\overset{u}{\longrightarrow}%
\Omega^{q}{}_{\varepsilon}\overline{\mathcal{KQ}}(A)\overset{v}%
{\longrightarrow}{}_{\varepsilon}\overline{\mathcal{KQ}}(A)\text{.}%
\]

\begin{proposition}
\label{negative cup positive}Let $m$ and $q$ be $2$-powers as in Convention
\ref{stronger convention}. Then the cup-product between the images of\textsf{
}the negative and positive Bott elements in%
\[
{}_{1}KQ_{\pm q}(\mathbb{Z}^{\prime};\,\mathbb{Z}/m)
\]
is reduced to $0$. Therefore, the compositions $v\circ u$ and $u\circ v$ are nullhomotopic.
\end{proposition}

\noindent\textbf{Proof.} Since, by \textsl{e.g. }\cite[pp. 797, 799]{BK},
${}_{1}KQ_{0}(\mathbb{Z}^{\prime};\,\mathbb{Z}/m)$ embeds in ${}_{1}%
KQ_{0}(\mathbb{R};\,\mathbb{Z}/m)\oplus{}_{1}KQ_{0}(\mathbb{F}_{3}%
;\,\mathbb{Z}/m)$, we consider separately the projections of the composites to
each summand. In each case, the key point is that the negative Bott element is
a power of an element in degree $-2$.\smallskip

\noindent\textbf{Projection to }${}_{1}KQ_{0}(\mathbb{R};\,\mathbb{Z}/m)$.\ We
compute the cup-product of the two Bott elements, using the first step in
\cite[Lemma 1.1]{Karoubifiltration}\footnote{Erratum: in the statement of
Lemma 1.1 of \cite{Karoubifiltration}, one should replace $8y$ by $16y$
because in the proof the inclusion {}$_{-1}W_{-6}^{\mathrm{top}}%
\hookrightarrow${}$_{1}W_{-8}^{\prime\mathrm{top}}$ is strict.}, that is, the
twelve-term exact sequences of \cite{KaroubiAnnals2} for both $\mathbb{Z}%
^{\prime}$ and the topological ring $\mathbb{R}$. As in \cite[Lemma
1.1]{Karoubifiltration}, they show that the map%
\[
\mathbb{Z}\oplus\mathbb{Z}/2\cong{}_{1}W_{-4}(\mathbb{Z}^{\prime
})\longrightarrow{}_{1}W_{-4}(\mathbb{R})\cong\mathbb{Z}%
\]
is given by $(w,\alpha)\mapsto2z$ where $z$ generates ${}_{1}W_{-4}%
(\mathbb{R})$. We now use some standard facts:

\begin{enumerate}
\item[(i)] there is a multiplicative isomorphism between $K_{n}(\mathbb{R})$
and ${}_{1}W_{n}(\mathbb{R})$ for all $n\in\mathbb{Z}$ \cite[Th\'{e}or\`{e}me
2.3]{K:AnnM112hgo}, and when $n$ is a multiple of $8$, each group is
$\mathbb{Z}$, generator $y_{n}$ say;

\item[(ii)] the cup-square $z^{2}$ of the generator $z$ of $K_{-4}%
(\mathbb{R})\cong{}_{1}W_{-4}(\mathbb{R})$ is $4$ times a generator $y_{-8}$
of $K_{-8}(\mathbb{R})\cong\mathbb{Z}$;

\item[(iii)] the cup-square of any generator of the free part of ${}_{1}%
W_{-4}(\mathbb{Z}^{\prime})$ projects to a generator of the free part of
${}_{1}W_{-8}(\mathbb{Z}^{\prime})\cong\mathbb{Z}\oplus\mathbb{Z}/2$.
\end{enumerate}

Let us write $q=2^{i+3}$ where $i\geq0$. From these facts, under the map%
\[
\mathbb{Z}\oplus\mathbb{Z}/2\cong{}_{1}W_{-q}(\mathbb{Z}^{\prime
})\longrightarrow{}_{1}W_{-q}(\mathbb{R})\cong\mathbb{Z}\text{,}%
\]
$(w,\alpha)^{2^{i+1}}$ is sent to $(2z)^{2^{i+1}}=2^{2^{i+1}}\cdot(4)^{2^{i}%
}y_{-q}=2^{q/2}y_{-q}$. Now consider the commuting diagram
\[%
\begin{array}
[c]{ccccc}%
{}_{1}W_{-q}(\mathbb{Z}^{\prime})\cong\mathbb{Z}\oplus\mathbb{Z}/2 &
\overset{\cong}{\longleftarrow} & {}_{1}KQ_{-q}(\mathbb{Z}^{\prime}%
)\cong\mathbb{Z}\oplus\mathbb{Z}/2 & \longrightarrow & K_{-q}(\mathbb{Z}%
^{\prime})=0\\
\downarrow &  & \downarrow &  & \downarrow\\
{}_{1}W_{-q}(\mathbb{R})\cong\mathbb{Z} & \overset{}{\longleftarrow} & {}%
_{1}KQ_{-q}(\mathbb{R})\cong\mathbb{Z}\oplus\mathbb{Z} & \longrightarrow &
K_{-q}(\mathbb{R})\cong\mathbb{Z}%
\end{array}
\]
Since the two lower horizontal maps correspond to the cokernel of the
hyperbolic map (on the left) and to the signature map (on the right), they
send a pair $(u,v)\in\mathbb{Z}\oplus\mathbb{Z}$ to $(u-v)y_{-q}\in{}%
_{1}W_{-q}(\mathbb{R})$ and $(u+v)y_{-q}\in K_{-q}(\mathbb{R})$ respectively.
Now, by (iii) above, the element $1\in\mathbb{Z}\subset{}_{1}KQ_{-q}%
(\mathbb{Z}^{\prime})$ may be taken as (up to sign) the projection of
$(w,\alpha)^{2^{i+1}}$, and therefore maps both on the left to $\pm
2^{q/2}y_{-q}\in{}_{1}W_{-q}(\mathbb{R})$ and on the right to $0\in
K_{-q}(\mathbb{R})$. Hence, the image $(u,v)\in{}_{1}KQ_{-q}(\mathbb{R})$ of
$1$ must have the form $\pm(2^{q/2-1},\,-2^{q/2-1})$.

Therefore, since always $m\leq2^{q/2-1}$ (the reason for the change from $p$
to $q$), if we now take $KQ$-theory with coefficients in $\mathbb{Z}/m$, then
the cup-product of the two Bott elements may be written
\[
\sigma=(0,0,\eta)\in{}_{1}KQ_{0}(\mathbb{Z}^{\prime};\,\mathbb{Z}%
/m)\cong\mathbb{Z}/m\oplus\mathbb{Z}/m\oplus\mathbb{Z}/2\text{,}%
\]
as calculated from the Bockstein exact sequence%
\[
{}_{1}KQ_{0}(\mathbb{Z}^{\prime})\longrightarrow{}_{1}KQ_{0}(\mathbb{Z}%
^{\prime})\longrightarrow{}_{1}KQ_{0}(\mathbb{Z}^{\prime};\,\mathbb{Z}%
/m)\longrightarrow{}_{1}KQ_{-1}(\mathbb{Z}^{\prime})=0
\]
and Lemma 3.11 of \cite{BKO}. Its image in ${}_{1}KQ_{0}(\mathbb{R}%
;\,\mathbb{Z}/m)\cong\mathbb{Z}/m\oplus\mathbb{Z}/m$ is thus $\pm
(2^{q/2-1},\,-2^{q/2-1})=(0,0)$.\smallskip

\noindent\textbf{Projection to }${}_{1}KQ_{0}(\mathbb{F}_{3};\,\mathbb{Z}%
/m)$.\ To compute the cup-product $\gamma$ of the images of the two Bott
elements in ${}_{1}KQ_{\pm q}(\mathbb{F}_{3};\,\mathbb{Z}/m)$, we exploit the
definition of the negative Bott element as the iterated power of an element in
${}_{-1}KQ_{-2}(\mathbb{Z}^{\prime})$. Therefore, $\gamma$ is the image of the
positive Bott element in ${}_{1}KQ_{q}(\mathbb{F}_{3};\,\mathbb{Z}/m)$ under
the following composition:%
\begin{align*}
&  {}_{1}KQ_{q}(\mathbb{F}_{3};\,\mathbb{Z}/m)\rightarrow{}_{-1}%
KQ_{q-2}(\mathbb{F}_{3};\,\mathbb{Z}/m)\rightarrow\cdots\rightarrow{}%
_{1}KQ_{4}(\mathbb{F}_{3};\,\mathbb{Z}/m)\rightarrow\\
\rightarrow &  {}_{-1}KQ_{2}(\mathbb{F}_{3};\,\mathbb{Z}/m)\rightarrow{}%
_{1}KQ_{0}(\mathbb{F}_{3};\,\mathbb{Z}/m)
\end{align*}
According to Friedlander \cite{Friedlander}, we have ${}_{-1}KQ_{2}%
(\mathbb{F}_{3})={}_{-1}KQ_{1}(\mathbb{F}_{3})=0$. Therefore, from another
Bockstein exact sequence, ${}_{-1}KQ_{2}(\mathbb{F}_{3};\,\mathbb{Z}/m)=0$,
and hence $\gamma=0$.

Finally, for the last part of the proposition, we use well-known facts in
cohomology theories \textrm{mod\ }$2^{k}$ \cite[I.~p.~75]{Araki-Toda} to prove
that the composite maps $v\circ u$ and $u\circ v$ are nullhomotopic. More
specifically, the multiplication by $2^{s}$ on cohomology theories
$\operatorname{mod}2^{k}$ is null-homotopic if $s\geq k$ and $s\geq2$%
.\hfill$\Box$\smallskip

For the next step, we need the following well-known Lemma (\textsl{cf.~}%
\cite[I.~p.~75]{Araki-Toda} again) which is a consequence of the splitting of
the multiplication by $m^{\prime}$ on the spectrum $S^{0}/m$, where $m$ and
$m^{\prime}$ are $2$-powers defined below.

\begin{lemma}
\label{cohomology splitting}Let $h^{\ast}$ be a cohomology theory represented
by a spectrum $\mathcal{S}$ and $m$ be a $2$-power. Let $h^{\ast}(-;$
$\mathbb{Z}/m)$ be the associated cohomology theory represented by the
spectrum $\mathcal{S}/m=\mathcal{S}\wedge S^{0}/m$. Finally, let
$\mathcal{T}_{m^{\prime}}$ be the homotopy fiber of the map%
\[
\mathcal{S}/m\longrightarrow\mathcal{S}/m
\]
defined by the multiplication by a $2$-power $m^{\prime}$, where $m^{\prime
}\geq\sup\{4,m\}$. Then we have a canonical splitting%
\[
\mathcal{T}_{m^{\prime}}\sim\mathcal{S}/m\times\Omega(\mathcal{S}/m)\text{.}%
\]

\end{lemma}

Let us denote by $\mathcal{F}$ the generic homotopy fiber of the maps
described before the lemma. There is a homotopy fibration
\[
\mathcal{F}(u)\longrightarrow\mathcal{F}(v\circ u)\longrightarrow
\mathcal{F}(v)\text{.}%
\]
According to the above considerations, $\mathcal{F}(u)$ is the spectrum
$\mathcal{P}_{\varepsilon}\overline{\mathcal{KQ}}(A)$, while $\mathcal{F}(v)$
is the spectrum $\Omega^{q-1}{}_{\varepsilon}\mathcal{KSC}^{(q)}\mathcal{(}%
A)$. On the other hand, as a consequence of Proposition
\ref{negative cup positive} and the previous lemma applied to the spectrum of
hermitian $K$-theory, $\mathcal{F}(v\circ u)$ may be canonically identified
with the product of spectra ${}_{\varepsilon}\overline{\mathcal{KQ}}%
(A)\times\Omega{}_{\varepsilon}\overline{\mathcal{KQ}}(A)$. Therefore, by
taking homotopy groups of the previous fibration, we get the exact sequence
\[
{}\overline{KSC}_{n+q}^{(q)}\rightarrow P\overline{KQ}_{n}\rightarrow
{}\overline{KQ}_{n}\oplus{}\overline{KQ}_{n+1}\rightarrow{}\overline
{KSC}_{n+q-1}^{(q)}\rightarrow P\overline{KQ}_{n-1}\,\text{.}%
\]
As a piece of convenient notation, set%
\[
\overline{KQ}_{n,n+1}=\overline{KQ}_{n}\oplus\overline{KQ}_{n+1}\,\text{.}%
\]

More generally, we shall also use the notation $\overline{KQ}_{X,X+1}$ for the
direct sum $\overline{KQ}_{X}\oplus\overline{KQ}_{X+1}$.

\begin{proposition}
\label{PKQ/KQ/KSC}We have the following two diagrams of exact sequences where
the vertical maps are induced respectively by the cup-product with the
negative or positive Bott element:%
\[%
\begin{tabular}
[c]{ccccc}%
$P\overline{KQ}_{n}$ & $\rightarrow$ & $\overline{KQ}_{n,n+1}$ & $\rightarrow$
& ${}\overline{KSC}_{n+q-1}^{(q)}$\\
$\uparrow P\beta_{n}^{\prime}$ &  & $\uparrow$ &  & $\uparrow\sigma
_{n}^{\prime}$\\
$P\overline{KQ}_{n+q}$ & $\rightarrow$ & $\overline{KQ}_{n+q,n+q+1}$ &
$\rightarrow$ & $\overline{KSC}_{n+2q-1}^{(q)}$%
\end{tabular}
\ \ \ \ \ \ \ \ \ \ \ \ \
\]
and%
\[%
\begin{tabular}
[c]{ccccc}%
$P{}\overline{KQ}_{n}$ & $\rightarrow$ & $\overline{KQ}_{n,n+1}$ &
$\rightarrow$ & $\overline{KSC}_{n+q-1}^{(q)}$\\
$\downarrow P\beta_{n}$ &  & $\downarrow$ &  & $\downarrow\sigma_{n}$\\
$P{}\overline{KQ}_{n+q}$ & $\rightarrow$ & $\overline{KQ}_{n+q,n+q+1}$ &
$\rightarrow$ & $\overline{KSC}_{n+2q-1}^{(q)}$%
\end{tabular}
\ \ \ \ \ \ \ \ \ \ \ \ \
\]
In these diagrams, $P\beta_{n}^{\prime}$ is an isomorphism if $n\geq d$ and a
monomorphism if $n=d-1$. We also have $P\beta_{n}=0$ if $n\geq d-1$. Finally,
$\sigma_{n}$ (resp.~$\sigma_{n}^{\prime}$) is an isomorphism (resp.~the zero
map) if $n\geq d-1$.
\end{proposition}

\noindent\textbf{Proof.} The second diagram is included in a bigger one with
horizontal exact sequences of $\overline{KQ}(\mathbb{Z}^{\prime})$-modules,
where we recall that $\overline{KQ}(\mathbb{Z}^{\prime})={}_{1}\overline
{KQ}(\mathbb{Z}^{\prime})\oplus{}_{-1}\overline{KQ}(\mathbb{Z}^{\prime})$:%
\[%
\begin{array}
[c]{lllll}%
\overline{KSC}_{n+q}^{(q)} & {\small \rightarrow}P\overline{KQ}_{n} &
{\small \rightarrow}\overline{KQ}_{n,n+1} & {\small \rightarrow}\overline
{KSC}_{n+q-1}^{(q)} & {\small \rightarrow}P\overline{KQ}_{n-1}\\
\multicolumn{1}{c}{{\small \downarrow}} &
\multicolumn{1}{c}{{\small \downarrow P}\beta_{n}} &
\multicolumn{1}{c}{{\small \downarrow}} &
\multicolumn{1}{c}{{\small \downarrow}} &
\multicolumn{1}{c}{{\small \downarrow}}\\
\overline{KSC}_{n+2q}^{(q)} & {\small \rightarrow}P\overline{KQ}_{n+q} &
{\small \rightarrow}\overline{KQ}_{n+q,n+q+1} & {\small \rightarrow}%
\overline{KSC}_{n+2q-1}^{(q)} & {\small \rightarrow}P\overline{KQ}_{n+q-1}%
\end{array}
\]
We claim that the second vertical map $P\beta_{n}$ is reduced to $0$ if $n\geq
d-1$. In order to see this, we consider the reverse map
\[
P\beta_{n}^{\prime}:P\overline{KQ}_{n+q}\longrightarrow P\overline{KQ}_{n}%
\]
given by the cup-product with the negative Bott element. This is a
monomorphism for $n\geq d-1$ according to Lemma \ref{negative perioidicity}.
Since the cup-product between the positive and the negative Bott elements is
equal to $0$ according to Proposition \ref{negative cup positive}, we have
necessarily $P\beta_{n}^{\prime}=0$. On the other hand, as we have seen in
Proposition \ref{KSCperiodicity}, the positive Bott map $\sigma_{n}$ is an
isomorphism%
\[
\overline{KSC}_{n+q-1}^{(q)}\cong\overline{KSC}_{n+2q-1}^{(q)}%
\]
for $n+q-1\geq d+q-2$, \textsl{i.e.} for $n\geq d-1$.

The reverse map $\sigma_{n}^{\prime}$ is reduced to $0$ since its composite
with $\sigma_{n}$ is trivial (note that all maps in these diagrams are
$KQ(\mathbb{Z}^{\prime})$-module maps).. \hfill$\Box$\smallskip

\begin{proposition}
\label{splitPKQ/KSC} If $n\geq d+q-1$, we have a split short exact sequence%
\[
0\rightarrow P{}_{\varepsilon}\overline{KQ}_{n}(A)\rightarrow{\small {}%
}_{\varepsilon}\overline{KQ}_{n,n+1}{}(A)\rightarrow{}_{\varepsilon}%
\overline{KSC}_{n+q-1}^{(q)}(A)\rightarrow0\text{.}%
\]

\end{proposition}

\noindent\textbf{Proof.} Let us consider a bigger diagram, where we now choose
$n\geq d+q-1$:%
\[%
\begin{tabular}
[c]{ccccccc}%
${}\overline{KSC}_{n}^{(q)}$ & $\rightarrow$ & $P{}\overline{KQ}_{n-q}$ &
$\rightarrow$ & ${\small {}}\overline{KQ}_{n-q,n-q+1}$ & $\rightarrow$ &
${}\overline{KSC}_{n-1}^{(q)}$\\
$\sigma_{n-q+1}\downarrow\cong$ &  & $\downarrow0$ &  & $\downarrow
\alpha_{n-q}$ & $\swarrow\gamma$ & $\sigma_{n-q}\downarrow\cong$\\
${}\overline{KSC}_{n+q}^{(q)}$ & $\rightarrow$ & $P{}\overline{KQ}_{n}$ &
$\rightarrow$ & ${\small {}}\overline{KQ}_{n,n+1}$ & $\rightarrow$ &
${}\overline{KSC}_{n+q-1}^{(q)}$\\
$\sigma_{n+1}\downarrow\cong$ &  & $\downarrow0$ &  & $\downarrow\alpha_{n}$ &
& $\sigma_{n}\downarrow\cong$\\
${}\overline{KSC}_{n+2q}^{(q)}$ & $\rightarrow$ & $P{}\overline{KQ}_{n+q}$ &
$\rightarrow$ & ${}\overline{KQ}_{n+q,n+q+1}$ & $\rightarrow$ & ${}%
\overline{KSC}_{n+2q-1}^{(q)}$%
\end{tabular}
\ \ \ \ \ \ \ \ \ \ \ \ \ \ \ \ \ \ \ \ \ \ \ \ \label{PKSC}%
\]
We would like to insert a map $\gamma:\overline{KSC}_{n-1}^{(q)}%
\longrightarrow$ $\overline{KQ}_{n,n+1}$ that renders this diagram
commutative. For this, we consider the other composition%
\[
{}\Omega^{q}{}_{\varepsilon}\overline{\mathcal{KQ}}(A)\overset{v}%
{\longrightarrow}{}_{\varepsilon}\overline{\mathcal{KQ}}(A)\overset
{u}{\longrightarrow}\Omega^{q}{}_{\varepsilon}\overline{\mathcal{KQ}%
}(A)\text{.}%
\]
We get a map from $\Omega^{-q}(\mathcal{F}(v))$ to $\Omega^{-q}(\mathcal{F}%
(u\circ v))$ that induces the required map $\gamma$ since $u\circ v$ is
nullhomotopic. The commutativity of the above diagram with $\gamma$ inserted
is a consequence of the homotopy commutative square%
\[%
\begin{tabular}
[c]{ccc}%
$\Omega^{q}{}_{\varepsilon}\overline{\mathcal{KQ}}(A)$ & $\overset
{v}{\longrightarrow}$ & $_{\varepsilon}\overline{\mathcal{KQ}}(A)$\\
$\downarrow\Omega^{q}(u)$ &  & $\downarrow u$\\
$\Omega^{2q}{}_{\varepsilon}\overline{\mathcal{KQ}}(A)$ & $\overset{\Omega
^{q}(v)}{\longrightarrow}$ & $\Omega^{q}{}_{\varepsilon}\overline
{\mathcal{KQ}}(A)$%
\end{tabular}
\ \ \ \ \ \
\]
where the vertical (resp.~horizontal) maps are defined by the cup-product with
the positive (resp.~negative) Bott element. Therefore, we get the short split
exact sequence%
\[
0\longrightarrow P{}\overline{KQ}_{n}\longrightarrow\overline{KQ}%
_{n,n+1}\longrightarrow\overline{KSC}_{n+q-1}^{(q)}\longrightarrow0\text{,}%
\]
which is an abbreviated formulation of our Proposition. \hfill$\Box$\smallskip

\begin{remark}
\label{decomposition}The group ${}_{\varepsilon}\overline{KQ}_{n,n+1}$
therefore decomposes in two distinct ways as the direct sum of two groups,
that is%
\[
{}_{\varepsilon}\overline{KQ}_{n,n+1}={}_{\varepsilon}\overline{KQ}_{n}%
\oplus{}_{\varepsilon}\overline{KQ}_{n+1}%
\]
and%
\[
{}_{\varepsilon}\overline{KQ}_{n,n+1}\cong P{}_{\varepsilon}\overline{KQ}%
_{n}\oplus{}_{\varepsilon}\overline{KSC}_{n+q-1}^{(q)}\text{.}%
\]
In the appendix to \cite{BKO} and also in the theory of stabilized Witt groups
\cite{Karoubistab.Witt}, we give many examples of rings $A$ such that
$K_{n}(A)=0$ for all $n\in\mathbb{Z}$ and therefore ${}_{\varepsilon}%
\overline{KSC}_{n+q-1}^{(q)}(A)=0$. This implies that our two direct sum
decompositions are not the same in general, since we may choose $A$ such that
the two groups ${}_{\varepsilon}\overline{KQ}_{n}(A)$ and ${}_{\varepsilon
}\overline{KQ}_{n+1}(A)$ are not $0$. Moreover, this remark may be used to
show that the higher ${}_{\varepsilon}\overline{KSC}$-theory depends \textsl{a
priori} on the sign of symmetry $\varepsilon$. An example of this fact is
$A=\mathbb{Z}^{\prime}$, where we know that $P{}_{\varepsilon}\overline
{KQ}_{n}(\mathbb{Z}^{\prime})=0$. On the other hand, from the table of the
$KQ$-groups of $\mathbb{Z}^{\prime}$ \cite{BK}, it is easy to see that ${}%
_{1}\overline{KQ}_{n,n+1}(\mathbb{Z}^{\prime})\neq{}_{-1}\overline{KQ}%
_{n,n+1}(\mathbb{Z}^{\prime})$ in general.
\end{remark}

For a better understanding of the periodicity statements we shall prove in
full generality in Section \ref{section: proof of main}, let us consider the
case where the $2$-primary abelian groups ${}_{\varepsilon}\overline{KQ}%
_{n}(A)$ are finite. The category of finite $2$-primary abelian groups is of
course well understood: its Grothendieck group (with respect to direct sums)
is freely generated by the groups $\mathbb{Z}/2^{k}$. On the other hand, it
follows from Proposition \ref{splitPKQ/KSC} that the groups ${}_{\varepsilon
}\overline{KQ}_{n,n+1}(A)={}_{\varepsilon}\overline{KQ}_{n}(A)\oplus
{}_{\varepsilon}\overline{KQ}_{n+1}(A)$ are periodic of period $q$ with
respect to $n$ for $n\geq d+q-1$. More precisely, $P{}_{\varepsilon}%
\overline{KQ}_{n}$ is periodic for $n\geq d$ according to Lemma
\ref{negative perioidicity} and ${}_{\varepsilon}\overline{KSC}_{n+q-1}^{(q)}$
is periodic for $n+q-1\geq d+q-2$, \textsl{i.e.} for $n\geq d-1$, according to
Proposition \ref{KSCperiodicity}.

Let us write $\alpha_{r}$ for the class of the group ${}_{\varepsilon
}\overline{KQ}_{r+d+q-1}$ in the Grothendieck group and put $\tau=\alpha
_{q}-\alpha_{0}$. We have the identities $\alpha_{q}=\alpha_{0}+\tau$,
$\alpha_{q+1}=\alpha_{1}-\tau$, $\alpha_{q+2}=\alpha_{2}+\tau$, \textsl{etc.}
In general, we may prove by induction on $s$ the formula%
\[
\alpha_{r+qs}=\alpha_{r}+(-1)^{r}s\tau
\]
when $r\geq0$.

\begin{proposition}
\label{finite KQ}Let us assume that the cup-product with the Bott element
induces an isomorphism%
\[
K_{n}(A;\mathbb{Z}/m)\rightarrow K_{n+p}(A;\mathbb{Z}/m)
\]
for $n\geq d$, and that the hermitian $K$-groups
${}_{\varepsilon}KQ_{n}(A;\,\mathbb{Z}/m)$ are finite for $n\geq d+q-1$. Then
these groups are periodic with respect to $n$ when $n\geq d+q-1$, of period
$q$. In particular, we have $\tau=0$ in the formulas above.
\end{proposition}

\noindent\textbf{Proof.} In the previous computation, let us write%
\[
\tau=%
{\textstyle\sum\limits_{i=1}^{k}}
u_{i}-%
{\textstyle\sum\limits_{j=1}^{k^{\prime}}}
v_{j}%
\]
where $u_{i}$ and $v_{i}$ are classes of nonzero irreducible modules and where
$k$ and $k^{\prime}$ are chosen minimal. We have the
two identities, valid for all $s\geq0$ :
\[
\alpha_{qs}=\alpha_{0}+s\tau
\]%
\[
\alpha_{qs+1}=\alpha_{1}-s\tau
\]
From the former, for all such $s$, the module $s%
{\textstyle\sum\limits_{j=1}^{k^{\prime}}}
v_{j}$ is always a summand of ${}_{\varepsilon}KQ_{d+q-1}(A;\,\mathbb{Z}/m)$.
Thus, $k^{\prime}=0$. Likewise, from the latter identity, $k=0$. Hence,
$\tau=0$.${}$ \hfill$\Box$\smallskip

Another example of a periodicity statement in hermitian $K$-theory is to
consider the case where $A$ is the ring with involution $B\times
B^{\mathrm{op}}$. Here $B^{\mathrm{op}}$ is the opposite algebra of $B$, the
involution on $A$ being defined by $(b,b^{\prime})\mapsto(b^{\prime},b)$. It
is easy to see that $_{\varepsilon}KQ_{n}(A)\cong K_{n}(B)$ and that the
negative periodicity map $_{\varepsilon}KQ_{n}%
(A)\rightarrow{}_{-\varepsilon}KQ_{n-2}(A)$ is reduced to $0$. Therefore, we
have the isomorphisms%
\[
_{\varepsilon}KSC_{n}(A)={}_{\varepsilon}KSC_{n}^{(2)}(A)\cong K_{n+1}%
(B)\oplus K_{n}(B)\text{,}%
\]
and, more generally%
\[
{}_{\varepsilon}KSC_{n}^{(p)}(A)\cong K_{n+1}(B)\oplus K_{n-p+2}(B)
\]
for $p$ a $2$-power. If we now reduce these theories \textrm{mod\ }$m$ and
assume the positive $p$-periodicity of the associated $\overline{K}$-groups
(for $n\geq d$), then the last identity can also be written as
\begin{align*}
{}_{\varepsilon}\overline{KSC}_{n+p-1}^{(p)}(A)  &  \cong\overline{K}%
_{n+p}(B)\oplus\overline{K}_{n+1}(B)\\
&  \cong\overline{K}_{n}(B)\oplus\overline{K}_{n+1}(B)\cong{}_{\varepsilon
}\overline{KQ}_{n}(A)\oplus{}_{\varepsilon}\overline{KQ}_{n+1}(A)\text{,}%
\end{align*}
which is a particular case of (\ref{splitPKQ/KSC}), since $P_{\varepsilon
}\overline{KQ}_{n}(A)=0$. Note that the positive Bott map%
\[
{}_{\varepsilon}KQ_{n}(A)\rightarrow{}_{\varepsilon}KQ_{n+p}(A)
\]
is here an isomorphism for $n\geq d$.

In the spirit of Thomason, one may consider the periodized $KQ$-theory, which
we shall denote by ${}_{\varepsilon}\overline{KQ}_{n}(A)\left[  \beta
^{-1}\right]  =\underrightarrow{\lim}{}_{\varepsilon}\overline{KQ}_{n+ps}(A)$.
We have the following theorem, quite similar to Connes' exact sequence
relating cyclic and Hochschild homologies \cite{Connes}\textsf{.}

\begin{proposition}
\label{KSC/splitting/KQ}Let us take $n\geq d+r-2$. Then (with the ring $A$
omitted from notation) we have the exact sequence%
\[%
\begin{array}
[c]{c}%
\cdots\rightarrow{}_{\varepsilon}\overline{KQ}_{n+2}\left[  \beta^{-1}\right]
\rightarrow{}_{\varepsilon^{\prime}}\overline{KQ}_{n+2-r}\left[  \beta
^{-1}\right]  \rightarrow{}_{\varepsilon}\overline{KSC}_{n}^{(r)}\rightarrow
{}_{\varepsilon}\overline{KQ}_{n+1}\left[  \beta^{-1}\right] \\
\rightarrow{}_{\varepsilon^{\prime}}\overline{KQ}_{n+1-r}\left[  \beta
^{-1}\right]  \rightarrow\cdots\rightarrow{}_{\varepsilon}\overline
{KQ}_{d-1+r}\left[  \beta^{-1}\right]  \rightarrow{}_{\varepsilon^{\prime}%
}\overline{KQ}_{d-1}\left[  \beta^{-1}\right]  ,
\end{array}
\]
where $\varepsilon^{\prime}=-\varepsilon$ if $r=2$ and $\varepsilon^{\prime
}=\varepsilon$ if $r>2$.\textsf{ }Moreover, if also $r\geq q$, we have a
splitting%
\[
{}_{\varepsilon}\overline{KSC}_{n}^{(r)}(A)\cong{}_{\varepsilon}\overline
{KQ}_{n+1}(A)\left[  \beta^{-1}\right]  \oplus{}_{\varepsilon}\overline
{KQ}_{n+2}(A)\left[  \beta^{-1}\right]  \text{.}%
\]

\end{proposition}

Note that for $r=q$, this splitting is also proved in
the beginning of the proof of Theorem\textsf{ }\ref{inv.dir.limit.theorem}%
\textsf{.}\smallskip

\noindent\textbf{Proof.\ }The first part of the proposition is a direct
consequence of Proposition \ref{KSCperiodicity} showing that the
$\overline{KSC}^{(r)}$-groups are periodic for $n\geq d+r-2$.

For the second part, we notice that the map between the $\overline{KQ}$-groups
is $0$, since the cup-product between the positive and negative Bott elements
is $0$, according to Proposition \ref{negative cup positive}. Thus, the
sequence decomposes into short exact sequences%
\[
0\rightarrow{}_{\varepsilon^{\prime}}\overline{KQ}_{n+2-r}(A)\left[
\beta^{-1}\right]  \longrightarrow{}_{\varepsilon}\overline{KSC}_{n}%
^{(r)}(A)\longrightarrow{}_{\varepsilon}\overline{KQ}_{n+1}(A)\left[
\beta^{-1}\right]  \rightarrow0\text{.}%
\]
Now, the inversion of $\beta$ yields an isomorphism
\[
{}_{\varepsilon^{\prime}}\overline{KQ}_{n+2-r}(A)\left[  \beta^{-1}\right]
\cong{}_{\varepsilon}\overline{KQ}_{n+2}(A)\left[  \beta^{-1}\right]
\]
since $q\,|\,r$. Finally, the splitting of these
sequences is as before a consequence of a general statement on cohomology
theories. One has to replace $n$ by a parameter space $X$, as we shall also do
in the next section.$\hfill\Box\smallskip$

\begin{Remarks}
The most interesting cases of this proposition are when one has $r=2$ or
$r=q$. The proposition also shows that the groups ${}_{\varepsilon}%
\overline{KSC}_{n}^{(r)}(A)$ are isomorphic for $r\geq q$.
\end{Remarks}

\section{Proof of the periodicity theorems\label{section: proof of main}}

Our aim in this Section is essentially to prove Theorems
\ref{Thm: invlimKQdirlim} and \ref{W-regularity}. We begin with a lemma
showing that it suffices to consider only the parameter $q$ of Convention
\ref{stronger convention}, rather than the desired period $p$ of Convention
\ref{convention: (p,p')} when dealing with direct or inverse limits.

\begin{lemma}
Let $m$ be a $2$-power and $p,q$ be as in Conventions \ref{convention: (p,p')}
and \ref{stronger convention}; that is,
\[%
\begin{array}
[c]{ccc}%
\underline{m} & \underline{p} & \underline{q}\\
\leq8 & 8 & 8\\
16 & 8 & 16\\
\geq32 & m/2 & m/2
\end{array}
\]
Then
\[
{}\underrightarrow{\lim}{}_{\varepsilon}KQ_{n+ps}(A;\,\mathbb{Z}%
/m)=\underrightarrow{\lim}{}_{\varepsilon}KQ_{n+qs}(A;\,\mathbb{Z}/m)
\]
and
\[
\underleftarrow{\lim}{}_{\varepsilon}KQ_{n+ps}(A;\,\mathbb{Z}%
/m)=\underleftarrow{\lim}{}_{\varepsilon}KQ_{n+qs}(A;\,\mathbb{Z}/m)\text{.}%
\]

\end{lemma}

\noindent\textbf{Proof.} We may focus on the exceptional case where $m=q=16$
and $p=8$. Here, by Theorem \ref{positive Bott element in KQ} there is a
positive Bott element $b^{+}\in{}_{1}\overline{KQ}_{8}(\mathbb{Z}^{\prime})$,
multiplication by which gives rise to the direct system of abelian groups%
\[
{}_{\varepsilon}\overline{KQ}_{n}(A)\longrightarrow{}_{\varepsilon}%
\overline{KQ}_{n+8}(A)\longrightarrow{}_{\varepsilon}\overline{KQ}%
_{n+16}(A)\longrightarrow{}_{\varepsilon}\overline{KQ}_{n+24}%
(A)\longrightarrow\cdots\,\text{.}%
\]
The direct limit of its subsystem%
\[
{}_{\varepsilon}\overline{KQ}_{n}(A)\longrightarrow{}_{\varepsilon}%
\overline{KQ}_{n+16}(A)\longrightarrow{}_{\varepsilon}\overline{KQ}%
_{n+32}(A)\longrightarrow\cdots
\]
appears in the exact sequence of Theorem \ref{inv.dir.limit.theorem} below.
However, since this subsystem is cofinal, its direct limit is precisely that
of the original system. In other words, we may replace the term
$\underrightarrow{\lim}{}_{\varepsilon}\overline{KQ}_{n+16s}(A)$ by
$\underrightarrow{\lim}{}_{\varepsilon}\overline{KQ}_{n+8s}(A)$.

Since the negative Bott element originates in ${}_{-1}KQ_{-2}(\mathbb{Z}%
^{\prime})$, a similar argument shows that the term $\underleftarrow{\lim}%
{}_{\varepsilon}\overline{KQ}_{n+16s}(A)$ may be replaced by $\underleftarrow
{\lim}{}_{\varepsilon}\overline{KQ}_{n+8s}(A)$. \hfill$\Box$\smallskip

We now start the proof of the periodicity theorems which will be a consequence
of our considerations in Section \ref{Higher KSC-theories}.

\begin{theorem}
\label{inv.dir.limit.theorem}\label{Periodicity large range}Let $A$ be a ring
with involution such that $1/2\in A$ and let $m$ and $p$ be $2$-powers
according to Convention \ref{convention: (p,p')}. We assume the existence of
an integer $d$, such that the cup-product with the Bott element in
$K_{p}(\mathbb{Z};\,\mathbb{Z}/m)$ induces an isomorphism%
\[
K_{n}(A;\,\mathbb{Z}/m)\overset{\cong}{\longrightarrow}K_{n+p}(A;\,\mathbb{Z}%
/m)\text{.}%
\]
for $n\geq d$. For such $n$, there is an exact sequence%
\[
\cdots\overset{\theta_{n}^{+}}{\longrightarrow}{}\underrightarrow{\lim}%
{}_{\varepsilon}KQ_{n+1+ps}(A;\,\mathbb{Z}/m)\rightarrow\underleftarrow{\lim
}{}_{\varepsilon}KQ_{n+ps}(A;\,\mathbb{Z}/m)
\]%
\[
\overset{\theta_{n}^{-}}{\longrightarrow}{}_{\varepsilon}KQ_{n}(A;\,\mathbb{Z}%
/m)\overset{\theta_{n}^{+}}{\longrightarrow}{}\underrightarrow{\lim}%
{}_{\varepsilon}KQ_{n+ps}(A;\,\mathbb{Z}/m)\text{,}%
\]
which for $n\geq d+q-1$ gives a split short exact sequence%
\[
0\rightarrow\underleftarrow{\lim}{}_{\varepsilon}KQ_{n+ps}(A;\,\mathbb{Z}%
/m)\overset{\theta_{n}^{-}}{\longrightarrow}{}_{\varepsilon}KQ_{n}%
(A;\,\mathbb{Z}/m)\overset{\theta_{n}^{+}}{\longrightarrow}{}\underrightarrow
{\lim}{}_{\varepsilon}KQ_{n+ps}(A;\,\mathbb{Z}/m)\rightarrow0\text{.}%
\]

\end{theorem}

\noindent\textbf{Proof.} For $n\geq d$, we consider the diagram of exact
sequences of Propositions \ref{PKQ/KQ/KSC} and \ref{splitPKQ/KSC} (for
convenience, we again drop the ring $A$ and the sign of symmetry $\varepsilon$
in the notation):%
\[%
\begin{array}
[c]{cclllcc}%
\overline{KSC}_{n+q}^{(q)} & \rightarrow & P\overline{KQ}_{n} &
{\small \rightarrow}\overline{KQ}_{n,n+1} & {\small \rightarrow}\overline
{KSC}_{n+q-1}^{(q)} & \rightarrow & P\overline{KQ}_{n-1}\\
&  & \multicolumn{1}{c}{P\beta_{n}{\small \downarrow0}} &
\multicolumn{1}{c}{{\small \downarrow\alpha}_{n}} & \multicolumn{1}{c}{\sigma
_{n}{\small \downarrow\cong}} &  & \\
0 & \rightarrow & P\overline{KQ}_{n+q} & {\small \rightarrow}\overline
{KQ}_{n+q,n+q+1} & {\small \rightarrow}\overline{KSC}_{n+2q-1}^{(q)} &
\rightarrow & 0
\end{array}
\]
Since the direct limit of the $P\overline{KQ}_{n+qs}$ is equal to $0$, we see
that%
\[
\underrightarrow{\lim}{}{}\overline{KQ}_{n+qs,n+1+qs}\cong{}\overline
{KSC}_{n+q-1}^{(q)}\cong\operatorname{Im}(\alpha_{n+q})
\]
which has already been proven in Proposition\textsf{ }\ref{KSC/splitting/KQ}%
\textsf{.}

We also have a reverse diagram of exact sequences%
\begin{equation}%
\begin{array}
[c]{ccccc}%
P{}\overline{KQ}_{n} & \rightarrow & {}\overline{KQ}_{n,n+1} & \rightarrow &
{}\overline{KSC}_{n+q-1}^{(q)}\rightarrow\\
P\beta_{n}^{\prime}\uparrow\cong &  & \uparrow\alpha_{n+q}^{\prime} &  &
\sigma_{n}^{\prime}\uparrow0\\
0\rightarrow P{}\overline{KQ}_{n+q} & \rightarrow & {}\overline{KQ}%
_{n+q,n+q+1} & \rightarrow & {}\overline{KSC}_{n+2q-1}^{(q)}\rightarrow0
\end{array}
\end{equation}
where the vertical maps are now induced by the cup-product with the negative
Bott element. The first vertical map is an isomorphism, while the last one is
reduced to $0$, by Proposition \ref{PKQ/KQ/KSC}.\textsf{ }

From the splitting of exact sequences afforded by Proposition
\ref{splitPKQ/KSC}, we have%
\[
\underleftarrow{\lim}{}\overline{KQ}_{n+qs,n+1+qs}\cong\underleftarrow{\lim}%
{}P{}\overline{KQ}_{n+qs}\cong\operatorname{Im}(\alpha_{n+q}^{\prime})\cong
P{}\overline{KQ}_{n}%
\]
by Lemma \ref{negative perioidicity}.

The first exact sequence in the first diagram above implies the exactness of
the middle row of the diagram%
\begin{equation}%
\begin{tabular}
[c]{ccccc}%
$\underleftarrow{\lim}{}\overline{KQ}_{n+1+qs}$ & $\rightarrow$ & ${}%
\overline{KQ}_{n+1}$ & $\overset{\theta_{n+1}^{+}}{\rightarrow}$ &
$\underrightarrow{\lim}{}\overline{KQ}_{n+1+qs}$\\
$\downarrow$ &  & $\downarrow$ &  & $\downarrow$\\
$\underleftarrow{\lim}{}\overline{KQ}_{n+qs,n+1+qs}$ & $\overset{\chi^{-}%
}{\rightarrow}$ & $\overline{KQ}_{n,n+1}$ & $\overset{\chi^{+}}{\rightarrow}$
& $\underrightarrow{\lim}{}\overline{KQ}_{n+qs,n+1+qs}$\\
$\downarrow$ &  & $\downarrow$ &  & $\downarrow$\\
$\underleftarrow{\lim}{}\overline{KQ}_{n+qs}$ & $\overset{\theta_{n}^{-}%
}{\rightarrow}$ & $\overline{KQ}_{n}$ & $\overset{\theta_{n}^{+}}{\rightarrow
}$ & $\underrightarrow{\lim}{}\overline{KQ}_{n+qs}$%
\end{tabular}
\ \ \ \ \ \ \ \
\end{equation}
Since the middle row is the direct sum of first and third rows, the exactness
of the middle row implies the exactness of the third row. We apply the same
argument for the left part of the required exact sequence.

Now suppose that $n\geq d+q-1$. Then Proposition \ref{splitPKQ/KSC} implies
that in the diagram above $\chi^{+}$ is a split epimorphism. The result
follows.\hfill$\Box$\smallskip

Recall from Definition \ref{defn: hermitian regular} in
the Introduction that a ring $A$ is \emph{hermitian regular} if the inverse
limits%
\[
\underset{\longleftarrow}{\lim}\overline{KQ}_{n+ps}(A)\quad\text{and\quad
}\underset{\longleftarrow}{\lim}^{1}\overline{KQ}_{n+ps}(A)
\]
are reduced to $0$ for all $n$.

\smallskip

\begin{Examples}
We remark that the second condition (with $\lim^{1}$) is always fulfilled if
$A$ has a periodic $\overline{K}$-group after a certain range, since in
Section \ref{Higher KSC-theories} we have shown that the inverse system
$\left\{  \overline{KQ}_{n+ps}(A)\right\}  $ satisfies the Mittag-Leffler
property. We have also seen in Section
\ref{Proof of the periodicity theorem for totally real 2-regular number fields}
that suitable rings of integers in a number field are hermitian regular. On
the other hand, according to Hu, Kriz and Ormsby \cite{HuKrizOrmsby}
(resp. Schlichting), if $k$ is a field of finite
mod\ $2$ virtual {\'{e}}tale cohomological dimension and of characteristic $0$
(resp. $p$), we have a homotopy equivalence%
\[
{}_{1}\overline{\mathcal{KQ}}(k)\simeq\overline{\mathcal{K}}(k)^{h\mathbb{Z}%
/2}\text{.}%
\]
Since $\overline{K}_{n}(k)\cong\overline{K}_{n+p}(k)$ for $n$ large enough by
\cite{Ostvar-Rosenschon} and \cite{Voevodskymod2}, the positive Bott map%
\[
{}_{1}KQ_{n}(k)\longrightarrow{}_{1}KQ_{n+p}(k)
\]
is also an isomorphism for $n$ large enough. The same statement is true for
the groups ${}{}_{-1}KQ_{n}$, as we see by relating the groups ${}_{1}%
KQ,{}_{-1}KQ$ and $KSC$ (see Section %
\ref{Higher KSC-theories}). As a conclusion, the negative Bott map%
\[
{}_{\varepsilon}\overline{KQ}_{n+p}(k)\longrightarrow{}_{\varepsilon}%
\overline{KQ}_{n}(k)
\]
is trivial for $n$ large enough, which implies that $k$ is hermitian regular
by Theorem \ref{inv.dir.limit.theorem}. More generally, as a special case of
the results in Section \ref{Generalization to schemes} and a
planned joint paper with Schlichting
\cite{Schlichting2}, any commutative algebra $A$ whose residual fields are all
of finite mod\ $2$ virtual {\'{e}}tale cohomological dimension is hermitian
regular. In \cite{Schlichting2}, these results will be
generalized in the scheme framework.
\end{Examples}

\begin{remark}
It is easy to see that the inverse systems of hermitian $K$-groups $\left\{
{}_{\varepsilon}\overline{KQ}_{n+ps}(A)\right\}  $ and Witt groups $\left\{
{}_{\varepsilon}\overline{W}_{n+ps}(A)\right\}  $ are equivalent since we have
the following factorization of the negative Bott map:%
\[
{}_{\varepsilon}\overline{KQ}_{n+8}(A)\rightarrow{}_{\varepsilon}\overline
{W}_{n+4}(A)\rightarrow{}_{\varepsilon}\overline{KQ}_{n}(A)\rightarrow
{}_{\varepsilon}\overline{W}_{n-4}(A)\text{.}%
\]
Therefore, the $\lim$ and $\lim^{1}$ groups may as well be computed with
higher Witt groups.
\end{remark}

\begin{theorem}
\label{Theorem of W-regularity}With the same hypotheses
as in the previous theorem, let us assume moreover that the ring $A$ is
hermitian regular. Then, for $n\geq d$, the positive Bott map%
\[
{}_{\varepsilon}\overline{KQ}_{n}(A){}\longrightarrow{}_{\varepsilon}%
\overline{KQ}_{n+p}(A)
\]
is an isomorphism.
\end{theorem}

\noindent\textbf{Proof.} According to Theorem \ref{inv.dir.limit.theorem}, it
is enough to show that $\theta_{d}^{+}$ is surjective. From the long exact
sequence of Proposition \ref{KSC/KSC}, we obtain the map of exact sequences
\begin{equation}%
\begin{tabular}
[c]{ccccccc}%
$\overline{KQ}_{d+q}$ & $\overset{s_{d+q}}{\rightarrow}$ & $\overline{KQ}_{d}$
& $\rightarrow$ & $\overline{KSC}_{d+q-2}^{(q)}$ & $\rightarrow$ &
$\overline{KQ}_{d+q-1}$\\
$\downarrow$ &  & $\downarrow\theta_{d}^{+}$ &  & $\downarrow\cong$ &  &
$\downarrow\gamma$\\
$\underrightarrow{\lim}{}\overline{KQ}_{d+qs}$ & $\overset{u}{\rightarrow}$ &
$\underrightarrow{\lim}{}\overline{KQ}_{d+qs}$ & $\overset{}{\rightarrow}$ &
$\overline{KSC}_{d+2q-2}^{(q)}$ & $\rightarrow$ & $\underrightarrow{\lim}%
{}\overline{KQ}_{d-1+qs}$%
\end{tabular}
\ \ \ \ \ \ \ \ \ \ \
\end{equation}
Observe from Theorem \ref{inv.dir.limit.theorem} that $\gamma$ is injective
because the inverse limits are reduced to $0$. Since $u$ is defined by the
cup-product with the negative Bott element, it is reduced to $0$. An
elementary diagram chase now shows that $\theta_{d}^{+}$ is surjective.\hfill
$\Box$\smallskip

Although for simplicity we have presented the arguments only in the case where
$X$ is a sphere, we observe that the groups obtained in the exact sequences of
the previous theorems can be considered as cohomology theories with respect to
pointed spaces $X$, if we replace the various spectra involved by their
sufficiently connected associated spectra (so that the low-dimensional
cohomology groups are trivial).

In order to pass from exactness of the sequence to split exactness as we did
at the end of Section \ref{Higher KSC-theories}, we may appeal to a general
fact about cohomology theories, that any surjective morphism like%
\[
{}_{\varepsilon}\theta_{X}^{+}:{}_{\varepsilon}\overline{KQ}_{X}%
(A)\longrightarrow{}\underrightarrow{\lim}{}_{\varepsilon}\overline{KQ}%
_{X+qs}(A)
\]
always admits a section.

The following theorems are the analogs of the previous ones with a parameter
space $X.$

\begin{theorem}
\label{Split exactness}Let $A$ be a ring with involution such that $1/2\in A$,
$m$, $p$ and $q$ be $2$-powers according to Convention
\ref{convention: (p,p')} and \ref{stronger convention}. We assume the
existence of an integer $d$ such that for $n\geq d$ the cup-product with the
Bott element in $K_{p}(\mathbb{Z};\,\mathbb{Z}/m)$ induces an isomorphism%
\[
K_{n}(A;\,\mathbb{Z}/m)\overset{\cong}{\longrightarrow}K_{n+p}(A;\,\mathbb{Z}%
/m)\text{.}%
\]
\emph{(a)} If $X$ is $(d+q-2)$-connected, we have a split short exact
sequence%
\[
0\rightarrow\underleftarrow{\lim}{}_{\varepsilon}KQ_{X+ps}(A;\,\mathbb{Z}%
/m)\overset{\theta_{n}^{-}}{\longrightarrow}{}_{\varepsilon}KQ_{X}%
(A;\,\mathbb{Z}/m)\overset{\theta_{n}^{+}}{\longrightarrow}{}\underrightarrow
{\lim}{}_{\varepsilon}KQ_{X+ps}(A;\,\mathbb{Z}/m)\rightarrow0\text{.}%
\]
As a consequence, the groups ${}_{\varepsilon}KQ_{X}(A;\,\mathbb{Z}/m)$ are
\textquotedblleft periodic\textquotedblright\ with respect to $X$ of period
$p$, more precisely%
\[
_{\varepsilon}KQ_{X}(A;\,\mathbb{Z}/m)\cong{}_{\varepsilon}KQ_{X+p}%
(A;\,\mathbb{Z}/m)\text{.}%
\]
In particular, if $n\geq d+q-1$, there is an isomorphism
\[
_{\varepsilon}KQ_{n}(A;\,\mathbb{Z}/m)\cong{}_{\varepsilon}KQ_{n+p}%
(A;\,\mathbb{Z}/m)\text{.}%
\]
\emph{(b) }Moreover, if $A$ is hermitian regular, the previous statements are
still true if we replace the number $q$ by $1$.
\end{theorem}

\begin{corollary}
For $A,X$ as in Theorem \ref{Split exactness}(a), and
$n\geq d+q-1$, the positive Bott map%
\[
\beta_{n}:{}_{\varepsilon}KQ_{n}(A;\,\mathbb{Z}/m)\longrightarrow
{}_{\varepsilon}KQ_{n+p}(A;\,\mathbb{Z}/m)
\]
has

\begin{enumerate}
\item[(i)] its image naturally isomorphic to the periodized $KQ$-theory, that
is%
\[
\operatorname{Im}(\beta_{n})\cong\underrightarrow{\lim}{}_{\varepsilon
}KQ_{n+ps}(A;\,\mathbb{Z}/m)\text{,}%
\]
and

\item[(ii)] its kernel and cokernel naturally isomorphic to $\underleftarrow
{\lim}{}_{\varepsilon}KQ_{n+ps}(A;\,\mathbb{Z}/m)$. Consequently, if
$\beta_{n}$ is either injective or surjective then it is an isomorphism.
\end{enumerate}

Moreover, if $A$ is hermitian regular, the same statements remain true on
replacing $q$ by $1$.
\end{corollary}

\noindent\textbf{Proof.} Here, we chase the following commutative diagram,
where the vertical maps are given by the cup-product with the positive Bott
element:%
\[%
\begin{tabular}
[c]{ccccccc}%
$0\rightarrow$ & $\underleftarrow{\lim}{}_{\varepsilon}\overline{KQ}%
_{n+ps}(A)$ & $\overset{\theta^{-}}{\longrightarrow}$ & ${}_{\varepsilon
}\overline{KQ}_{n}(A)$ & $\overset{\theta^{+}}{\longrightarrow}$ &
$\underrightarrow{\lim}{}_{\varepsilon}\overline{KQ}_{n+ps}(A)$ &
$\rightarrow0$\\
& $\downarrow0$ &  & $\downarrow\beta_{n}$ &  & $\downarrow\cong$ & \\
$0\rightarrow$ & $\underleftarrow{\lim}{}_{\varepsilon}\overline{KQ}%
_{n+p+ps}(A)$ & $\overset{\theta^{-}}{\longrightarrow}$ & ${}_{\varepsilon
}\overline{KQ}_{n+p}(A)$ & $\overset{\theta^{+}}{\longrightarrow}$ &
$\underrightarrow{\lim}{}_{\varepsilon}\overline{KQ}_{n+p+ps}(A)$ &
$\rightarrow0$%
\end{tabular}
\vspace{-10pt}%
\]
$\hfill\Box\medskip$

\begin{remark}
\label{remark: theta not iso}One may ask if the map $\theta_{n}^{+}$ in
Theorem \ref{Periodicity large range} is an isomorphism in general for $n$
sufficiently large. This is not the case however, as is shown in the Appendix
C to \cite{BKO} and in \cite{Karoubistab.Witt}, where we give many examples of
rings with trivial $K$-theory and nontrivial $KQ$-theory. It follows from the
$12$-term exact sequence of \cite[p. 278]{KaroubiAnnals2}, that for such rings
$\theta_{n}^{-}$ is an isomorphism. From the short exact sequence of Theorem
\ref{Periodicity large range}, $\theta_{n}^{+}$ must therefore vanish,
although the $KQ$-theory is nontrivial. Thus, $\theta_{n}^{+}$ fails to be an
isomorphism. Other examples may be found in the paper of Hu, Kriz and Ormsby
\cite{HuKrizOrmsby} for commutative rings and schemes.

However, one may hope it is so for the examples of commutative rings $A$
considered in the Introduction, which are of \textquotedblleft geometric
nature\textquotedblright. See also Section \ref{Generalization to schemes} for
an analogous conjecture in the category of schemes.
\end{remark}

We finish this section with an application to the computation of the
$\overline{KQ}$-groups in terms of the $\overline{K}$-groups when these groups
are finite. (This result formally dates from the
December 2010 resubmission of the paper.) For reading convenience, we again
suppress  the index $\varepsilon\in\{\pm1\}$.

\begin{theorem}
\label{kq k inequality}Let us assume the hypotheses of 
Theorem \ref{Split exactness} and that for $n\geq d$ the $\overline{K}_{n}%
$-groups are finite, of order $k_{n}\,$. Then for $n\geq d$ the $\overline
{KQ}_{n}$-groups are finite, of order $\mathrm{kq}_{n}$ subject to the
inequality%
\[
\mathrm{kq}_{n}+\mathrm{kq}_{n+1}\leq k_{d}+k_{d+1}+\cdots+k_{d+q-1}\,
\]
and equality $\mathrm{kq}_{n}=\mathrm{kq}_{n+q}\,$.
\end{theorem}

\noindent\textbf{Proof. }According to (\ref{KSC/KSC exact seq}), there is an
exact sequence
\[
\cdots\longrightarrow\overline{KSC}_{n}^{(r/2)}\longrightarrow\overline
{KSC}_{n}^{(r)}\longrightarrow\overline{KSC}_{n-r/2}^{(r/2)}\longrightarrow
\cdots\,\text{.}%
\]
In the case $r=2$, the outer two groups are respectively $\overline{K}_{n+1}$
and $\overline{K}_{n}$, and so assumed finite when $n\geq d$. In general, if
we denote by $s_{n}^{r}$ the order of the group $\overline{KSC}_{n}^{(r)}$
when it is finite, we therefore have the inequality%
\[
s_{n}^{r}\leq s_{n-r/2}^{r/2}+s_{n}^{r/2}\text{.}%
\]
For instance, with $n\geq d$,%
\[%
\begin{array}
[c]{ll}%
s_{n+1}^{2}\leq & k_{n+1}+k_{n+2}\\
s_{n+3}^{4}\leq & s_{n+1}^{2}+s_{n+3}^{2}\leq k_{n+1}+k_{n+2}+k_{n+3}%
+k_{n+4}\\
\quad\vdots & \\
s_{n+q-1}^{q}\leq & k_{n+1}+k_{n+2}+\cdots+k_{n+q}=k_{d}+\cdots+k_{d+q-1}%
\,\text{,}%
\end{array}
\]
where the equality follows from the assumption of $\overline{K}$-periodicity.
On the other hand, since the ring $A$ is hermitian regular, according to
Theorem \ref{inv.dir.limit.theorem} and its proof we have%
\[
s_{n+q-1}^{q}=\mathrm{kq}_{n}+\mathrm{kq}_{n+1}%
\]
for $n\geq d$. This gives the required inequality. Finally, the equality comes
from Theorem \ref{Theorem of W-regularity}.\hfill$\Box$\smallskip

\begin{example}
If $m=8$, we have $q=8$. Therefore, for $n\geq d$, we have the inequality%
\[
\mathrm{kq}_{n}+\mathrm{kq}_{n+1}\leq k_{d}+k_{d+1}+\cdots+k_{d+7}\text{.}%
\]
In particular, if we assume that the groups
${}_{\varepsilon}KQ_{n}(A)$ and${}_{\varepsilon}KQ_{n+1}(A)$ are finitely
generated, and that for some $r$ and $\varepsilon$ the Witt group
${}_{\varepsilon}W_{r}(A)$ has an infinite free summand, then by periodicity
of Witt groups after tensoring with $\mathbb{Q}$, one of the four groups
${}_{\varepsilon}\overline{KQ}_{n},$ ${}_{\varepsilon}\overline{KQ}_{n+1}$
($\varepsilon\in\{\pm1\}$) must be nonzero. It now follows that at least one
of the groups $K_{d}(A),\ldots,K_{d+7}(A)$ must also be nonzero.
\end{example}

\section{The case of odd prime power coefficients\label{section: odd prime}}

For the sake of completeness, we should also study $KQ$-theory with odd prime
power coefficients, which is much easier to handle, starting from known
results in $K$-theory. As is well known, if $\ell$ is an odd prime, there is a
remarkable Bott element $b_{K}$ in the group $K_{2(\ell-1)\ell^{\nu-1}%
}(\mathbb{Z};\,\mathbb{Z}/\ell^{\nu})$ (see Section
\ref{section: Bott elements} of this paper). In particular, its image in the
topological $K$-group%
\[
K_{2(\ell-1)\ell^{\nu-1}}(\mathbb{R};\,\mathbb{Z}/\ell^{\nu})\cong
K_{2(\ell-1)\ell^{\nu-1}}(\mathbb{C};\,\mathbb{Z}/\ell^{\nu})\cong%
\mathbb{Z}/\ell^{\nu}%
\]
is the image of an integral Bott generator in $K_{2(\ell-1)\ell^{\nu-1}%
}(\mathbb{C})\cong\mathbb{Z}$. According to Convention
\ref{convention: (p,p')}, we shall write $p=2(\ell-1)\ell^{\nu-1}$ (for the
period) and $m=\ell^{\nu}$ (for the order of the coefficient group).

Let us assume now that $A$ is one of the examples of algebras described in the
Introduction. For odd primes, we use the Bloch-Kato conjecture, which is now a
theorem proven by Rost and Voevodsky, \textsl{cf.}~\cite{Rost1}, \cite{Rost2},
\cite{Voevodskymodp}, \cite{HW}, \cite{Suslin1}, \cite{Suslin2}, and
\cite{Weibel}. This implies that the cup-product with the Bott element $b$
induces an isomorphism%
\[
K_{n}(A;\,\mathbb{Z}/m)\overset{\cong}{\longrightarrow}K_{n+p}(A;\,\mathbb{Z}%
/m)
\]
whenever $n\geq d$ for the type of rings $A$ considered in the Introduction.

In order to extend our previous results to hermitian $K$-theory with odd prime
power coefficients, it is convenient to describe more geometrically elements
of the $K$-groups and $KQ$-groups. This description in terms of
\textquotedblleft virtual\textquotedblright\ flat bundles is given in detail
in Appendix 1 of \cite{Karoubiasterisque}. For instance, an element of
${}_{\varepsilon}KQ_{n}(A)$ can be described as a flat $A$-bundle $E$ over an
homology sphere of dimension $n$, provided with a nondegenerate quadratic form
$q$. With this language, we can easily define an involution on the
$KQ$-groups: it is induced by the correspondence%
\[
(E,q)\mapsto(E,-q).
\]
Let us now consider the groups ${}_{\varepsilon}KQ_{n}(A)^{\prime}%
={}_{\varepsilon}KQ_{n}(A)\otimes_{\mathbb{Z}}\mathbb{Z}^{\prime}$. The tensor
product of virtual $A$-bundles (when $A$ is commutative) induces a ring
structure on the direct sum of all these groups (with $\varepsilon=\pm1$). On
the other hand, the previous involution enables us to split each group
${}_{\varepsilon}KQ_{n}(A)^{\prime}$ as a direct sum ${}_{\varepsilon}%
KQ_{n}(A)_{+}^{\prime}\oplus{}_{\varepsilon}KQ_{n}(A)_{-}^{\prime}$.

\begin{lemma}
The sum decomposition of the ${}_{\varepsilon}KQ_{n}(A)^{\prime}$ described
above is a ring product decomposition. In other words, the cup-product map
between elements of $KQ_{+}^{\prime}$ and $KQ_{-}^{\prime}$ is reduced to $0$.
\end{lemma}

\noindent\textbf{Proof.} Since we make $2$ invertible in the $KQ^{\prime}%
$-groups involved, one may think of an element $z$ of $KQ_{+}^{\prime}$ as a
sum $(E,q)+(E,-q)$ and an element $z^{\prime}$ of $KQ_{-}^{\prime}$ as a
difference $(E^{\prime},q^{\prime})-(E^{\prime},-q^{\prime})$, where $q$ and
$q^{\prime}$ are $\varepsilon$- and $\varepsilon^{\prime}$-quadratic forms.
The product $z\cdot z^{\prime}$ is therefore (with $q\otimes q^{\prime}$ an
$\varepsilon\varepsilon^{\prime}$-quadratic form)%
\[
(E\otimes E^{\prime},\,q\otimes q^{\prime})+(E\otimes E^{\prime},\,-q\otimes
q^{\prime})-(E\otimes E^{\prime},\,q\otimes-q^{\prime})-(E\otimes E^{\prime
},\,-q\otimes-q^{\prime}),
\]
which is of course zero. \hfill$\Box$\smallskip

\begin{corollary}
\label{ring product}The ring product decomposition of the direct sum of the
groups ${}_{\varepsilon}KQ_{n}(A)^{\prime}$ induces a ring product
decomposition of the direct sum of the groups $KQ_{n}(A;\,\mathbb{Z}/m)$ when
$m$ is an odd prime power $>3$.
\end{corollary}

\noindent\textbf{Proof.} The corollary follows from general arguments about
cohomology theories $\operatorname{mod}m$ \cite{Araki-Toda}. \hfill$\Box
$\smallskip

\begin{remark}
One also has an involution on the $K$-groups induced by the duality functor,
as was noticed already in Section \ref{Higher KSC-theories}. If we perform the
tensor product by $\mathbb{Z}\left[  1/2\right]  $ or we take coefficients in
the group $\mathbb{Z}/m$ with $m$ odd, the symmetric part is in bijective
correspondence with the symmetric part of the corresponding $KQ$-group. This
correspondence is induced by the forgetful functor or the hyperbolic functor
\cite{KaroubiAnnals2}.
\end{remark}

\begin{remark}
Before introducing the Bott elements in this situation, it is worth mentioning
that the fundamental theorem in hermitian $K$-theory holds for arbitrary rings
(we no longer assume that $1/2\in A$) when we localize away from $2$. The
details may be found in \cite[Lemma 1.1]{Karoubifiltration}. More precisely,
in this case the symmetric part $\mathcal{KQ}(A)_{+}^{\prime}$ of the spectrum
of $\mathcal{KQ}(A)^{\prime}$ is the symmetric part of the spectrum
$\mathcal{K}(A)^{\prime}$, whereas the antisymmetric part $\mathcal{KQ}%
(A)_{-}^{\prime}$ has periodic homotopy groups of period $4$, which are the
higher Witt groups. This remark also applies when we take $\operatorname{mod}%
\ m$ coefficients with $m$ odd.
\end{remark}

Let $b_{+}$ denote the Bott element in $_{1}KQ_{p}(\mathbb{Z};\,\mathbb{Z}%
/m)_{+}=$ $_{1}\overline{KQ}_{p}(\mathbb{Z})_{+}$ corresponding to the usual
Bott element $b_{K}$ in
\[
K_{p}(\mathbb{Z};\,\mathbb{Z}/m)=K_{p}(\mathbb{Z};\,\mathbb{Z}/m)_{+}%
=\overline{K}_{p}(\mathbb{Z})_{+}\text{.}%
\]
The following Theorem is now obvious, since $\overline{K}_{+}\cong%
\overline{KQ}_{+}$.

\begin{theorem}
Let $A$ be any ring such that the cup-product map with the $K$-theory Bott
element $b_{K}$ induces an isomorphism%
\[
\overline{K}_{n}(A)\cong\overline{K}_{n+p}(A)
\]
for $n\geq d$. Then, taking cup-product with $b_{+}$ also induces an
isomorphism%
\[
_{\varepsilon}\overline{KQ}_{n}(A)_{+}\cong{}_{\varepsilon}\overline{KQ}%
_{n+p}(A)_{+}%
\]
for the same values of $n$.$\hfill\Box$
\end{theorem}

On the other hand, there is another\ \textquotedblleft Bott
element\textquotedblright\ $b^{\prime}$ that lies in ${}_{1}\overline{KQ}%
_{p}(\mathbb{Z})_{-}={}_{1}W_{p}(\mathbb{Z};\,\mathbb{Z}/m)$, which is the
higher Witt group $\operatorname{mod}\ m$. It is the image mod\ $m$ of a
suitable power of $u\in{}_{-1}W_{2}(\mathbb{Z})$ constructed in
\cite{KaroubiAnnals2} and \cite[Theorem 1.4]{Karoubifiltration}.

\begin{theorem}
For an odd prime $\ell$, let $p=2(\ell-1)\ell^{\nu-1}$ and $m=$ $\ell^{\nu}$
as in (\ref{convention: (p,p')}). Let $A$ be any ring such that its $K$-theory
$\mathrm{mod}$\ $m$ is periodic for $n\geq d$, the periodicity being given by
the cup-product with the Bott element $b_{K}$. As usual, denote by
$\overline{KQ}$ the $KQ$-groups $\operatorname{mod}\ m$. Then, for $n\geq d$,
there is an isomorphism%
\[
{}_{\varepsilon}\overline{KQ}_{n}(A)\overset{\theta_{+-}}{\longrightarrow}%
{}\underrightarrow{\lim}{}_{\varepsilon}\overline{KQ}_{n+ps}(A)
\]
where $\theta_{+-}$ is induced by the cup-product with the sum $b_{+}%
+b^{\prime}$ of the two previously defined Bott elements in the group
\[
{}_{1}\overline{KQ}_{p}(\mathbb{Z})={}_{1}\overline{KQ}_{p}(\mathbb{Z}%
)_{+}\oplus{}_{1}\overline{KQ}_{p}(\mathbb{Z})_{-}\text{.}%
\]

\end{theorem}

\noindent\textbf{Proof.} If we consider the direct sum%
\[
{}_{\varepsilon}\overline{KQ}_{n}(A)={}_{\varepsilon}\overline{KQ}_{n}%
(A)_{+}\oplus{}_{\varepsilon}\overline{KQ}_{n}(A)_{-}\text{,}%
\]
then the cup-product with this element $b^{\prime}$ is trivial on the summand
${}_{\varepsilon}\overline{KQ}_{n}(A)_{+}$, and induces an isomorphism from
${}_{\varepsilon}\overline{KQ}_{n}(A)_{-}$ to ${}_{\varepsilon}\overline
{KQ}_{n+p}(A)_{-}$ after \cite[Th\'{e}or\`{e}me 3.9]{K:AnnM112hgo}. In a
parallel way, the cup-product with $b_{+}$ is trivial on the term
${}_{\varepsilon}\overline{KQ}_{n}(A)_{-}$ and induces an isomorphism from
${}_{\varepsilon}\overline{KQ}_{n}(A)_{+}$ to ${}_{\varepsilon}\overline
{KQ}_{n+p}(A)_{+}$, as shown in Corollary \ref{ring product}. The theorem
follows. \hfill$\Box$\smallskip

\begin{remark}
By choosing a \textquotedblleft negative\ Bott element\textquotedblright\ in
${}_{1}\overline{KQ}_{-p}(\mathbb{Z})_{-}$, we obtain an equivalent version of
the previous theorem, in analogy with Theorem \ref{Thm: invlimKQdirlim}, as
the following short split exact sequence:%
\[
0\rightarrow\underleftarrow{\lim}{}_{\varepsilon}\overline{KQ}_{n+ps}%
(A)\overset{\theta^{-}}{\longrightarrow}{}_{\varepsilon}\overline{KQ}%
_{n}(A)\overset{\theta^{+}}{\longrightarrow}{}\underrightarrow{\lim}%
{}_{\varepsilon}\overline{KQ}_{n+ps}(A)\rightarrow0\text{.}%
\]
In this exact sequence the inverse system (resp.~direct system) is given by
taking cup-product with the negative (resp.~positive) Bott element in ${}%
_{1}\overline{KQ}_{-p}(\mathbb{Z})_{-}$ (resp. ${}_{1}\overline{KQ}%
_{p}(\mathbb{Z})_{+}$).
\end{remark}

\begin{remark}
For simplicity, we have assumed the ring $A$ commutative in order to define
internal cup-products. However, a closer look at the arguments shows that we
have in fact used \textquotedblleft external\textquotedblright\ cup-products
of the type%
\[
\overline{KQ}(A)\times\overline{KQ}(\mathbb{Z)\longrightarrow}\overline
{KQ}(A)
\]
Therefore, the previous theorem extends easily to noncommutative rings, such
as group rings.
\end{remark}

\section{Generalization to schemes and \'{e}tale theories}

\label{Generalization to schemes}

The proof and the statement of Theorem \ref{Thm: invlimKQdirlim} apply
verbatim to schemes for which $2$ is invertible. This follows since the
fundamental theorem in hermitian $K$-theory has been generalized by
Schlichting -- see his work in progress on exact categories with weak
equivalences and duality \cite{Schlichting}. (Of course, in the context of
CW-spectra, weak equivalences are in fact homotopy equivalences.) In this
section we view this generalization in the context of the \'{e}tale descent
problem for hermitian $K$-theory.

Let $S$ be a regular and separated noetherian scheme of finite Krull dimension
(in the interest of generalizing these assumptions the inclined reader may
compare with \cite{Ostvar-Rosenschon2} and \cite{Schlichting}). 
Throughout, for a fixed prime $\ell$, 
we assume that $\mathcal{O}_{S}$ is a sheaf of $\mathbb{Z}[1/\ell]$-modules.
In particular, 
for the important case $\ell=2$, 
$2$ is invertible in $S$. 
\vspace{0.1in}

For the definition of the hermitian $K$-theory spectrum ${}_{\varepsilon
}\mathcal{KQ}(S)$, and the forgetful and hyperbolic maps between the algebraic
and hermitian $K$-theory of $S$, we refer to Schlichting's work
\cite{Schlichting}. In particular, as for rings, we may form the fibers
${}_{\varepsilon}\mathcal{V}(S)$ and ${}_{\varepsilon}\mathcal{U}(S)$ of the
forgetful and hyperbolic maps, respectively. The generalization of the
fundamental theorem to schemes in \cite{Schlichting} shows that there is a
homotopy equivalence
\begin{equation}
{}_{-\varepsilon}\mathcal{V}(S)\simeq\Omega\,{}_{\varepsilon}\mathcal{U}(S)
\end{equation}
such that the composite map
\[
\Omega^{2}\,{}_{\varepsilon}\mathcal{KQ}(S)\rightarrow\Omega\,{}_{\varepsilon
}\mathcal{U}(S)\rightarrow{}_{-\varepsilon}\mathcal{V}(S)\rightarrow
{}_{-\varepsilon}\mathcal{KQ}(S)
\]
is given by the cup-product with the negative Bott element in ${}_{-1}%
KQ_{-2}(\mathbb{Z}^{\prime})$.

With these results in hand, the proof of Theorem \ref{Thm: invlimKQdirlim}
carries over to the setting of schemes. That is, if cup-product with the Bott
element in $K_{p}(\mathbb{Z};\,\mathbb{Z}/m)$, $m$ and $p$ being $2$-powers
linked by our Convention \ref{convention: (p,p')}, induces an isomorphism
\[
K_{n}(S;\,\mathbb{Z}/m)\overset{\cong}{\longrightarrow}K_{n+p}(S;\,\mathbb{Z}%
/m)
\]
for $n\geq d$, then, for $n\geq d+q-1$ there is a split short exact sequence
\[
0\rightarrow\underleftarrow{\lim}{}_{\varepsilon}\overline{KQ}_{n+ps}%
(S)\overset{\theta^{-}}{\longrightarrow}{}_{\varepsilon}\overline{KQ}%
_{n}(S)\overset{\theta^{+}}{\longrightarrow}{}\underrightarrow{\lim}%
{}_{\varepsilon}\overline{KQ}_{n+ps}(S)\rightarrow0.
\]
We expect that the map $\theta^{+}$ is an isomorphism in many cases of
geometric interest. (For $n\geq d$, and not just $\geq d+q-1$: 
compare with Theorem \ref{Theorem of W-regularity}.)
This question is closely related to the so-called
\'{e}tale descent problem for hermitian $K$-theory and explicit computations,
which are our main concerns in this section. \vspace{0.1in}

Jardine introduced in \cite{Jardinesupercohernence} the Bott periodic
\'{e}tale hermitian $K$-theory spectra of $S$ with mod\ $\ell^{\nu}%
$-coefficients
\[
{}_{\varepsilon}\mathcal{KQ}^{\text{\'{e}t}}/\ell^{\nu}(S)\text{.}%
\]
More precisely, ${}_{\varepsilon}\mathcal{KQ}^{\text{\'{e}t}}/\ell^{\nu}(S)$
is the $\mathcal{KO}$- or equivalently $\mathcal{KU}$-localization of
Jardine's \'{e}tale hermitian $K$-theory. By construction of the Bott periodic
\'{e}tale theory, there exists an induced mod\ $\ell^{\nu}$ comparison map
\[
\Gamma_{S}:{}_{\varepsilon}\mathcal{KQ}/\ell^{\nu}(S)\longrightarrow
{}_{\varepsilon}\mathcal{KQ}^{\text{\'{e}t}}/\ell^{\nu}(S)
\]
obtained by taking global sections of a globally fibrant model for the
presheaf ${}_{\varepsilon}\mathcal{KQ}/\ell^{\nu}(\,\,)$ on some sufficiently
large \'{e}tale site of $S$.

Similarly, the mod\ $\ell^{\nu}$ \'{e}tale self-conjugate $K$-theory
$\mathcal{KSC}^{\text{\'{e}t}}/\ell^{\nu}(S)$ of $S$ is defined by taking a
globally fibrant model of the presheaf $\mathcal{KSC}/\ell^{\nu}(\,\,)$. Later
in this section we shall make use of specific fibrant models. \vspace{0.1in}

Recall that the notation $\mathrm{vcd}_{\ell}$ stands for the mod\ $\ell$
virtual cohomological dimension. When $\ell$ is odd, then $\mathrm{vcd}_{\ell
}$ coincides with the usual mod\ $\ell$ cohomological dimension $\mathrm{cd}%
_{\ell}$. Although $\mathrm{cd}_{2}$ is infinite in the examples of
$\mathbb{Z}[1/2]$ and $\mathbb{R}$, the number $\mathrm{vcd}_{2}$ is finite in
both cases. For a proof of the next result we refer to \cite[Proposition
6.1]{Mitchellprime2}.

\begin{lemma}
\label{hypercohomologymodel} If $\mathrm{vcd}_{\ell}(S)<\infty$ then the
\'{e}tale hypercohomology presheaf
\[
{\mathbb{H}}_{\mathrm{\acute{e}t}}^{\bullet}(-,\,{}_{\varepsilon}%
\mathcal{KQ}/\ell^{\nu}(\,\,))
\]
is a globally fibrant model for ${}_{\varepsilon}\mathcal{KQ}/\ell^{\nu}(-)$.
The same fibrancy result holds for the presheaf $\mathcal{KSC}/\ell^{\nu}%
(-)$.$\hfill\Box$
\end{lemma}

The \'{e}tale descent problem for self-conjugate $K$-theory can be solved
easily using the solution for algebraic $K$-theory \cite{Ostvar-Rosenschon}.
For a point $s\in S$, let $k(s)$ denote the corresponding residue field.

\begin{theorem}
\label{theorem:KSCdescent} The comparison map
\[
\mathcal{KSC}/\ell^{\nu}(S)\longrightarrow\mathcal{KSC}^{\text{$\mathrm{\acute
{e}t}$}}/\ell^{\nu}(S)
\]
is a weak equivalence on $\sup\{\mathrm{vcd}_{\ell}(k(s))-2\}_{s\in S}%
$-connected covers.

Hence, if $\mathrm{vcd}_{\ell}(S)<\infty$ then there is a weak equivalence
\[
L_{\mathcal{KU}}\,\mathcal{KSC}/\ell^{\nu}(S)\longrightarrow\mathcal{KSC}%
^{\text{\'{e}t}}/\ell^{\nu}(S)\simeq{\mathbb{H}}_{\text{\'{e}t}}^{\bullet
}(S,\,L_{\mathcal{KU}}\,\mathcal{KSC}/\ell^{\nu}(\,\,))\text{.}%
\]

\end{theorem}

\noindent\textbf{Proof.} There exists a naturally induced commutative diagram
of fiber sequences of presheaves of spectra%
\[%
\begin{array}
[c]{ccccc}%
\mathcal{KSC}/\ell^{\nu}(-) & \longrightarrow & \mathcal{K}/\ell^{\nu}(-) &
\longrightarrow & \mathcal{K}/\ell^{\nu}(-)\\
\downarrow &  & \downarrow &  & \downarrow\\
\mathcal{KSC}^{\text{\'{e}t}}/\ell^{\nu}(-) & \longrightarrow & \mathcal{K}%
^{\text{\'{e}t}}/\ell^{\nu}(-) & \longrightarrow & \mathcal{K}^{\text{\'{e}t}%
}/\ell^{\nu}(-)
\end{array}
\vspace{-15pt}%
\]
$\hfill\Box\bigskip$

Similarly to Theorem \ref{theorem:KSCdescent}, we do not expect that the
comparison map $\Gamma_{S}$ is a weak equivalence in general, but in many
cases of interest it should be a weak equivalence on some connected cover.
\vspace{0.1in}

Let $\ell$ be an odd prime, and let the subscript $+$ as in $K_{+}$ denote the
symmetric part of $K$-theory. Then the forgetful and hyperbolic functors
induce isomorphisms between the symmetric parts of the algebraic and hermitian
$K$-groups of $S$ at $\ell$. This allows us to infer the following result by
referring to \cite{Ostvar-Rosenschon} and to Corollary
\ref{corollary:KQstalks} for the symmetry of \'{e}tale hermitian $K$-theory at
odd primes.

\begin{theorem}
\label{theorem:oddKQdescent} Let $\ell\neq2$. The comparison map
\[
{}_{\varepsilon}\mathcal{KQ}/\ell^{\nu}(S)_{+}\longrightarrow{}_{\varepsilon
}\mathcal{KQ}^{\text{\'{e}t}}/\ell^{\nu}(S)_{+}\simeq{}_{\varepsilon
}\mathcal{KQ}^{\text{\'{e}t}}/\ell^{\nu}(S)
\]
is a weak equivalence on $\sup\{\mathrm{vcd}_{\ell}(k(s))-2\}_{s\in S}%
$-connected covers.

Hence, if $\mathrm{vcd}_{\ell}(S)<\infty$ then there is a weak equivalence
\[
L_{\mathcal{KU}}\,{}_{\varepsilon}\mathcal{KQ}/\ell^{\nu}(S)_{+}%
\longrightarrow{}_{\varepsilon}\mathcal{KQ}^{\text{\'{e}t}}/\ell^{\nu
}(S)\simeq{\mathbb{H}}_{\text{\'{e}t}}^{\bullet}(S,L_{\mathcal{KU}}%
\,{}_{\varepsilon}\mathcal{KQ}/\ell^{\nu}(\,\,))\text{.}\vspace{-15pt}%
\]
$\hfill\medskip\Box$\vspace{0.1in}
\end{theorem}

At $\ell=2$, we prove Theorem \ref{theorem:splotsurjectivity} stated in the
Introduction. In the following, let $n\geq\sup\{\mathrm{vcd}_{2}%
(k(s))-1\}_{s\in S}+q-1$. Inverting the positive Bott element in the direct
sum decomposition
\[
{}_{\varepsilon}\overline{KQ}_{n,n+1}(S)\cong P{}_{\varepsilon}\overline
{KQ}_{n}(S)\oplus{}_{\varepsilon}\overline{KSC}_{n+p-1}^{(p)}(S)
\]
yields
\[
{}_{\varepsilon}\overline{KQ}_{n,n+1}(S)[\beta^{-1}]\cong{}_{\varepsilon
}\overline{KSC}_{n+p-1}^{(p)}(S)[\beta^{-1}].
\]
In order to simplify the right hand side of this isomorphism, we first use the
fact that if $\mathrm{vcd}_{2}(S)<\infty$ then there is a weak equivalence
\[
\overline{\mathcal{KSC}}^{(r)}(S)[\beta^{-1}]\longrightarrow{\overline
{\mathcal{KSC}}^{(r)}}^{\text{\'{e}t}}(S)\text{.}%
\]
Then, by \'{e}tale descent for self-conjugate $K$-theory shown in Theorem
\ref{theorem:KSCdescent}, the induced comparison map
\[
\overline{\mathcal{KSC}}^{(r)}(S)\longrightarrow{\overline{\mathcal{KSC}%
}^{(r)}}^{\text{\'{e}t}}(S)
\]
is a weak equivalence on $\sup\{\mathrm{vcd}_{2}(k(s))+r-4\}_{s\in S}%
$-connected covers. Hence, there is an isomorphism
\[
{}_{\varepsilon}\overline{KSC}_{n+p-1}^{(p)}(S)[\beta^{-1}]\cong%
{}_{\varepsilon}\overline{KSC}_{n+p-1}^{(p)}(S)\text{.}%
\]
As a result,
\[
{}_{\varepsilon}\overline{KQ}_{n,n+1}(S)={}_{\varepsilon}\overline{KQ}%
_{n}(S)\oplus{}_{\varepsilon}\overline{KQ}_{n+1}(S)
\]
maps by a split surjection onto its Bott localization
\[
{}_{\varepsilon}\overline{KQ}_{n,n+1}(S)[\beta^{-1}]\cong{}_{\varepsilon
}\overline{KQ}_{n}^{\text{\'{e}t}}(S)\oplus{}_{\varepsilon}\overline{KQ}%
_{n+1}^{\text{\'{e}t}}(S)\text{.}%
\]
By looking at one component at a time, we deduce that there is a split
surjection
\[
{}_{\varepsilon}\overline{KQ}_{n}(S)\longrightarrow{}_{\varepsilon}%
\overline{KQ}_{n}(S)[\beta^{-1}]
\]
and an isomorphism
\[
{}_{\varepsilon}\overline{KQ}_{n}(S)[\beta^{-1}]\cong{}_{\varepsilon}%
\overline{KQ}_{n}^{\text{\'{e}t}}(S)
\]
for all $n$. $\hfill\Box$\vspace{0.1in}

Our next result is a local-global
comparison theorem. It will be a consequence of the homotopical setup due to
Jardine, see \textsl{e.g.}~\cite{Jardinebook}, and the rigidity theorem for
hermitian $K$-theory of henselian pairs proven by the second author in
\cite[Theorem 4]{KaroubiJPAA34}, \textsl{cf.}~the unpublished work
\cite{Jardinerigidity} for an approach using homotopy theory of simplicial
presheaves. The specific result we use is as follows.

\begin{theorem}
\label{theorem:rigidity} (\cite{KaroubiJPAA34}) Let $(A,I)$ be a henselian
pair with $q\in A^{\times}\cap\mathbb{Z}$ such that $I$ is invariant by the
involution on $A$, and $\lambda+\overline{\lambda}=1$ for some $\lambda\in A$
if $q$ is even. Then the map of rings with involutions $A\rightarrow A/I$
induces an isomorphism
\[
{}_{\varepsilon}{KQ}_{n}(A;\,\mathbb{Z}/q)\overset{\cong}{\longrightarrow}%
{}_{\varepsilon}{KQ}_{n}(A/I;\,\mathbb{Z}/q)
\]
for all $\varepsilon$ and $n\geq0$.
\end{theorem}

We note that the sharper bound for the connected covers in the theorem below
(relative to that in Theorem \ref{theorem:splotsurjectivity}) equals the one
shown for algebraic $K$-theory in \cite{Ostvar-Rosenschon}.

\begin{theorem}
\label{comparisontheorem} Suppose that $\Gamma_{k(s)}$ is a weak equivalence
on $(\mathrm{vcd}_{2}(k(s))-2)$-connected covers for every residue field
$k(s)$ of $S$. Then the comparison map
\[
\Gamma_{S}:{}_{\varepsilon}\mathcal{KQ}/2^{\nu}(S)\longrightarrow
{}_{\varepsilon}\mathcal{KQ}^{\text{\textrm{$\acute{e}t$}}}/2^{\nu}(S)
\]
is a weak equivalence on $\sup\{\mathrm{vcd}_{2}(k(s))-2\}_{s\in S}$-connected
covers.$\hfill\Box$
\end{theorem}

\noindent
\textbf{Proof.} There is a functorially induced commutative diagram with the
mod\ $2$ comparison map displayed on top:%
\[%
\begin{tabular}
[c]{ccc}%
${}_{\varepsilon}\mathcal{KQ}/2(S)$ & $\longrightarrow$ & ${\mathbb{H}%
}_{\mathrm{\acute{e}t}}^{\bullet}(S,\,L_{\mathcal{KU}}\,{}_{\varepsilon
}\mathcal{KQ}/2(\,\,))$\\
$\downarrow$ &  & $\downarrow$\\
${\mathbb{H}}_{\mathrm{Nis}}^{\bullet}(S,\,{}_{\varepsilon}\mathcal{KQ}%
/2(\,\,))$ & $\longrightarrow$ & ${\mathbb{H}}_{\mathrm{Nis}}^{\bullet
}(S,\,{\mathbb{H}}_{\mathrm{\acute{e}t}}^{\bullet}(S,\,L_{\mathcal{KU}}%
\,{}_{\varepsilon}\mathcal{KQ}/2(\,\,)))$%
\end{tabular}
\ \ \ \ \ \ \ \ \ \label{6.6}%
\]
We claim that the vertical maps are weak equivalences. By the Nisnevich
descent theorem in \cite{Schlichting}, this holds for the left hand side. For
the right hand side, the \'{e}tale topology is finer than the Nisnevich one;
so, the direct image functor maps ${\mathbb{H}}_{\text{\'{e}t}}^{\bullet
}(S,\,L_{\mathcal{KU}}\,{}_{\varepsilon}\mathcal{KQ}/2(\,\,))$ to a globally
fibrant object on the Nisnevich site of $S$. We claim that the mod\ $2$
comparison map is a stalkwise weak equivalence on the given connected cover
for the Nisnevich topology. In fact, let $A$ be a Hensel local ring with
residue field $k$ and consider the functorially induced commutative diagram%
\[%
\begin{tabular}
[c]{ccc}%
${}_{\varepsilon}\mathcal{KQ}/2(A)$ & $\longrightarrow$ & ${}_{\varepsilon
}\mathcal{KQ}/2(k))$\\
$\downarrow$ &  & $\downarrow$\\
${\mathbb{H}}_{\text{\'{e}t}}^{\bullet}(A,\,L_{\mathcal{KU}}\,{}_{\varepsilon
}\mathcal{KQ}/2(\,\,))$ & $\longrightarrow$ & ${\mathbb{H}}_{\text{\'{e}t}%
}^{\bullet}(k,\,L_{\mathcal{KU}}\,{}_{\varepsilon}\mathcal{KQ}/2(\,\,))$%
\end{tabular}
\ \ \ \ \ \ \ \ \
\]
Combining the previous theorem concerning rigidity for hermitian $K$-theory
and the equivalence between the \'{e}tale sites of $A$ and $k$, this reduces
the stalkwise weak equivalence to the assumed case of fields. It follows that
the lower horizontal map in Diagram (\ref{6.6}) is a stalkwise weak
equivalence on the same connected covers between globally fibrant objects, and
hence it is a pointwise weak equivalence on $\sup\left\{  vcd(k(s)-2\right\}
_{s\in S}$-connected covers. $\hfill\Box$\medskip

Motivated by the local hypotheses of Theorem \ref{comparisontheorem}, we make
the following forecast of the outcome of the \'{e}tale descent problem for
hermitian $K$-theory.

\begin{conjecture}
\label{conjecture:comparisonforfields} Suppose that $k$ is a field of
characteristic $\neq2$. Then the comparison map
\[
\Gamma_{k}:{}_{\varepsilon}\mathcal{KQ}/2^{\nu}(k)\longrightarrow
{}_{\varepsilon}\mathcal{KQ}^{\text{\textrm{$\acute{e}t$}}}/2^{\nu}(k)
\]
is a weak equivalence on ($\mathrm{vcd}_{2}(k)-2$)-connected covers.
\end{conjecture}

Conjecture \ref{conjecture:comparisonforfields}, in conjunction with Theorem
\ref{comparisontheorem}, predicts that, in many cases of interest, hermitian
$K$-theory is Bott periodic on some connected cover. Our earlier results on
Bott periodicity can be taken as oblique evidence for this conjecture. For
$n\geq\mathrm{vcd}_{2}(k)+q-1$, recall the exact sequence of $KQ$-groups with
$2$-power coefficients:
\[
0\rightarrow\underleftarrow{\lim}{}_{\varepsilon}\overline{KQ}_{n+ps}%
(k)\overset{\theta^{-}}{\longrightarrow}{}_{\varepsilon}\overline{KQ}%
_{n}(k)\overset{\theta^{+}}{\longrightarrow}{}\underrightarrow{\lim}%
{}_{\varepsilon}\overline{KQ}_{n+ps}(k)\rightarrow0.
\]
By Bott periodicity, Conjecture \ref{conjecture:comparisonforfields} implies
that the inverse limit is trivial, \textsl{i.e.}
\[
\underleftarrow{\lim}{}_{\varepsilon}\overline{KQ}_{n+ps}(k)=0\text{.}%
\]
In other words, the field $k$ should be hermitian regular (Definition
\ref{defn: hermitian regular}). Conversely, if the
inverse limit is trivial, then there is an isomorphism
\[
\theta^{+}\colon{}_{\varepsilon}\overline{KQ}_{n}(k)\overset{\cong%
}{\longrightarrow}{}\underrightarrow{\lim}{}_{\varepsilon}\overline{KQ}%
_{n+ps}(k)
\]
for $n\geq\mathrm{vcd}_{2}(k)-1$,
according to our Theorem \ref{Theorem of W-regularity}.
As noted in the Introduction, a proof of the
above conjecture is to appear in a joint paper with Schlichting
\cite{Schlichting2}.

Lemma \ref{hypercohomologymodel} can be motivated by the conditionally
convergent right half-plane cohomological descent spectral sequence
established by Thomason \cite{Thomason}:
\begin{equation}
{}_{\varepsilon}E_{2}^{p,q}=H_{\text{\textrm{$\acute{e}t$}}}^{p}%
(S,\,\widetilde{\pi}_{q}L_{\mathcal{KU}}\,{}_{\varepsilon}\mathcal{KQ}%
/\ell^{\nu}(\,\,))\Longrightarrow\pi_{q-p}{\mathbb{H}}_{\text{\textrm{$\acute
{e}t$}}}^{\bullet}(S,\,L_{\mathcal{KU}}\,{}_{\varepsilon}\mathcal{KQ}%
/\ell^{\nu}(\,\,)). \label{descent.SS}%
\end{equation}
Here the coefficient sheaf indicated by $\widetilde{\pi}_{\ast}$ is the
\'{e}tale sheafification of the presheaf of stable homotopy groups $\pi_{\ast
}L_{\mathcal{KU}}\,{}_{\varepsilon}\mathcal{KQ}/\ell^{\nu}(\,\,)$. The concept
of \textquotedblleft conditional convergence\textquotedblright\ for spectral
sequences was introduced by Boardman in \cite{Boardman}. A useful consequence
is that the descent spectral sequence (\ref{descent.SS}) converges strongly
provided that there exists only a finite number of nontrivial differentials.
Thus, the spectral sequence (\ref{descent.SS}) is strongly convergent if $S$
has finite mod\ $\ell$ \'{e}tale cohomological dimension. (This will be the
case in all the examples we consider.) The $d_{r}$-differential in
(\ref{descent.SS}) has bidegree $(r,\,1-r)$. \vspace{0.1in}

In order to identify the \'{e}tale stalks of ${}_{\varepsilon}\mathcal{KQ}%
/\ell^{\nu}(\,\,)$, and consequently the $E_{2}$-page of (\ref{descent.SS}),
\textsl{cf.}~\cite{Jardinesupercohernence}, \cite[Theorem 2.6]{Snaith}, we
invoke the Rigidity Theorem \ref{theorem:rigidity} together with the homotopy
equivalences
\begin{equation}
{}_{\varepsilon}\mathcal{KQ}/\ell^{\nu}(A)\simeq{}_{\varepsilon}%
\mathcal{KQ}/\ell^{\nu}(\mathbb{C})\simeq\left\{
\begin{array}
[c]{lll}%
\mathcal{K}/\ell^{\nu}(\mathbb{R}) & \  & \varepsilon=1\\
\Omega^{4}\mathcal{K}/\ell^{\nu}(\mathbb{R}) &  & \varepsilon=-1
\end{array}
\right.
\end{equation}
for a strict Hensel local ring $A$. The above is very similar to the case of
algebraic $K$-theory, where the \'{e}tale sheaf associated to the presheaf
\[
U\mapsto\pi_{n}\mathcal{K}/\ell^{\nu}(U)
\]
is the Tate twisted sheaf of roots of unity $\mu_{\ell^{\nu}}^{\otimes k}$
when $n=2k$ is even, and trivial when $n$ is odd. For $\mathcal{KSC}$ and
${}_{\varepsilon}\mathcal{KQ}$ at $\ell$ we have:

\begin{corollary}
(\cite{Snaith}) The \'{e}tale sheaf associated to the presheaf
\[
U\mapsto\pi_{n}(\mathcal{KSC}/\ell^{\nu}(U))
\]
is given by: \begin{table}[tbh]
\begin{center}%
\begin{tabular}
[c]{p{0.5in}|p{0.5in}|p{0.5in}|}\hline
$n\bmod4$ & $\ell=2$ & $\ell\neq2$\\\hline
$4k$ & $\mu_{2^{\nu}}^{\otimes2k}$ & $\mu_{\ell^{\nu}}^{\otimes2k}$\\
$4k+1$ & $\mu_{2}^{\otimes2k+1}$ & $0$\\
$4k+2$ & $\mu_{2}^{\otimes2k+1}$ & $0$\\
$4k+3$ & $\mu_{2^{\nu}}^{\otimes2k+2}$ & $\mu_{\ell^{\nu}}^{\otimes2k+2}%
$\\\hline
\end{tabular}
\end{center}
\end{table}
\end{corollary}

\begin{corollary}
\label{corollary:KQstalks} (\cite{Jardinesupercohernence}, \cite{Snaith}) The
\'{e}tale sheaf associated to the presheaf
\[
U\mapsto\pi_{n}({}_{\varepsilon}\mathcal{KQ}/\ell^{\nu}(U))
\]
is given as follows.

\begin{enumerate}
\item For $\ell=2$ by: \begin{table}[tbh]
\begin{center}%
\begin{tabular}
[c]{p{0.5in}|p{0.8in}|p{0.9in}|p{0.8in}|p{0.9in}|}\hline
$n\bmod8$ & $\varepsilon=1$, $\nu=1$ & $\varepsilon=-1$, $\nu=1$ &
$\varepsilon=1$, $\nu>1$ & $\varepsilon=-1$, $\nu>1$\\\hline
$8k$ & $\mu_{2}^{\otimes4k}$ & $\mu_{2}^{\otimes4k}$ & $\mu_{2^{\nu}}%
^{\otimes4k}$ & $\mu_{2^{\nu}}^{\otimes4k}$\\
$8k+1$ & $\mu_{2}^{\otimes4k+1}$ & $0$ & $\mu_{2}^{\otimes4k+1}$ & $0$\\
$8k+2$ & $\mu_{4}^{\otimes4k+1}$ & $0$ & $(\mu_{2}^{\otimes4k+1})^{\oplus2}$ &
$0$\\
$8k+3$ & $\mu_{2}^{\otimes4k+1}$ & $0$ & $\mu_{2}^{\otimes4k+1}$ & $0$\\
$8k+4$ & $\mu_{2}^{\otimes4k+2}$ & $\mu_{2}^{\otimes4k+2}$ & $\mu_{2^{\nu}%
}^{\otimes4k+2}$ & $\mu_{2^{\nu}}^{\otimes4k+2}$\\
$8k+5$ & $0$ & $\mu_{2}^{\otimes4k+3}$ & $0$ & $\mu_{2}^{\otimes4k+3}$\\
$8k+6$ & $0$ & $\mu_{4}^{\otimes4k+3}$ & $0$ & $(\mu_{2}^{\otimes
4k+3})^{\oplus2}$\\
$8k+7$ & $0$ & $\mu_{2}^{\otimes4k+3}$ & $0$ & $\mu_{2}^{\otimes4k+3}$\\\hline
\end{tabular}
\end{center}
\end{table}

\item For $\ell\neq2$ and $\varepsilon=\pm1$ by $\mu_{\ell^{\nu}}^{\otimes2k}$
if $n=4k$, and trivial otherwise.
\end{enumerate}
\end{corollary}

\begin{remark}
In Corollary \ref{corollary:KQstalks}, the $(4,2)$-periodicity in the change
of symmetry between the $\varepsilon=1$ and $\varepsilon=-1$ cases in the
table for $\ell=2$ is given by cup-product with a generator of ${}_{-1}%
{KQ}_{4}(\mathbb{C};\,\mathbb{Z}/2^{\nu})$. The case $\ell\neq2$ is similar.
In degrees $8k+2$, recall that $\mathbb{R}P^{2}$ is a $\mathrm{mod}$\ $2$
Moore space and $\widetilde{KO}(\mathbb{R}P^{2})\cong\mathbb{Z}/4$ generated
by the tangent bundle, while the universal coefficient sequence splits for
$\nu>1$.
\end{remark}

As a consequence of Lemma \ref{hypercohomologymodel} and Corollary
\ref{corollary:KQstalks}, we conclude that if $\mathrm{vcd}_{2}(S)<\infty$
then there exist conditionally convergent cohomological spectral sequences
\[
{}_{1}E_{2}^{p,q}=\left\{
\begin{array}
[c]{lll}%
H_{\text{$\mathrm{\acute{e}t}$}}^{p}(S,\mu_{2}^{\otimes4k}) & q=8k & \\
H_{\text{$\mathrm{\acute{e}t}$}}^{p}(S,\mu_{2}^{\otimes4k+1}) & q=8k+1 & \\
H_{\text{$\mathrm{\acute{e}t}$}}^{p}(S,\mu_{4}^{\otimes4k+1}) & q=8k+2 & \\
H_{\text{$\mathrm{\acute{e}t}$}}^{p}(S,\mu_{2}^{\otimes4k+1}) & q=8k+3 & \\
H_{\text{$\mathrm{\acute{e}t}$}}^{p}(S,\mu_{2}^{\otimes4k+2}) & q=8k+4 & \\
0 & q\equiv5,6,7\;(\mathrm{mod}\ 8) &
\end{array}
\right\}  \Longrightarrow{}_{1}{KQ}_{q-p}^{\mathrm{\acute{e}t}}/2(S),
\]
and
\[
{}_{-1}E_{2}^{p,q}=\left\{
\begin{array}
[c]{lll}%
H_{\text{$\mathrm{\acute{e}t}$}}^{p}(S,\mu_{2}^{\otimes4k}) & q=8k & \\
H_{\text{$\mathrm{\acute{e}t}$}}^{p}(S,\mu_{2}^{\otimes4k+2}) & q=8k+4 & \\
H_{\text{$\mathrm{\acute{e}t}$}}^{p}(S,\mu_{2}^{\otimes4k+3}) & q=8k+5 & \\
H_{\text{$\mathrm{\acute{e}t}$}}^{p}(S,\mu_{4}^{\otimes4k+3}) & q=8k+6 & \\
H_{\text{$\mathrm{\acute{e}t}$}}^{p}(S,\mu_{2}^{\otimes4k+3}) & q=8k+7 & \\
0 & q\equiv1,2,3\;(\mathrm{mod}\ 8) &
\end{array}
\right\}  \Longrightarrow{}_{-1}{KQ}_{q-p}^{\mathrm{\acute{e}t}}/2(S).
\]
\vspace{0.1in}

And for $\nu\geq2$, the descent spectral sequences take the forms
\[
{}_{1}E_{2}^{p,q}=\left\{
\begin{array}
[c]{lll}%
H_{\text{$\mathrm{\acute{e}t}$}}^{p}(S,\mu_{2^{\nu}}^{\otimes4k}) & q=8k & \\
H_{\text{$\mathrm{\acute{e}t}$}}^{p}(S,\mu_{2}^{\otimes4k+1}) & q=8k+1 & \\
H_{\text{$\mathrm{\acute{e}t}$}}^{p}(S,\mu_{2}^{\otimes4k+1})^{\oplus2} &
q=8k+2 & \\
H_{\text{$\mathrm{\acute{e}t}$}}^{p}(S,\mu_{2}^{\otimes4k+1}) & q=8k+3 & \\
H_{\text{$\mathrm{\acute{e}t}$}}^{p}(S,\mu_{2^{\nu}}^{\otimes4k+2}) & q=8k+4 &
\\
0 & q\equiv5,6,7\;(\mathrm{mod}\ 8) &
\end{array}
\right\}  \Longrightarrow{}_{1}{KQ}_{q-p}^{\mathrm{\acute{e}t}}/2^{\nu}(S),
\]
and
\[
{}_{-1}E_{2}^{p,q}=\left\{
\begin{array}
[c]{lll}%
H_{\text{$\mathrm{\acute{e}t}$}}^{p}(S,\mu_{2^{\nu}}^{\otimes4k}) & q=8k & \\
H_{\text{$\mathrm{\acute{e}t}$}}^{p}(S,\mu_{2^{\nu}}^{\otimes4k+2}) & q=8k+4 &
\\
H_{\text{$\mathrm{\acute{e}t}$}}^{p}(S,\mu_{2}^{\otimes4k+3}) & q=8k+5 & \\
H_{\text{$\mathrm{\acute{e}t}$}}^{p}(S,\mu_{2}^{\otimes4k+3})^{\oplus2} &
q=8k+6 & \\
H_{\text{$\mathrm{\acute{e}t}$}}^{p}(S,\mu_{2}^{\otimes4k+3}) & q=8k+7 & \\
0 & q\equiv1,2,3\;(\mathrm{mod}\ 8) &
\end{array}
\right\}  \Longrightarrow{}_{-1}{KQ}_{q-p}^{\mathrm{\acute{e}t}}/2^{\nu}(S).
\]
\vspace{0.1in}

For $\ell\neq2$ and $\mathrm{cd}_{\ell}(S)<\infty$, the descent spectral
sequence takes the form
\[
{}_{\varepsilon}E_{2}^{p,q}=\left\{
\begin{array}
[c]{lll}%
H_{\mathrm{\acute{e}t}}^{p}(S,\mu_{\ell^{\nu}}^{\otimes\frac{q}{2}}) &
q\equiv0\;(\mathrm{mod}\ 4) & \\
0 & q\not \equiv 0\;(\mathrm{mod}\ 4) &
\end{array}
\right\}  \Longrightarrow{}_{\varepsilon}{KQ}_{q-p}^{\mathrm{\acute{e}t}}%
/\ell^{\nu}(S).
\]
Note that the $E_{2}$-pages are independent of the symmetry $\varepsilon$.
This is not surprising since on symmetric parts $K$-theory mod\ $\ell^{\nu}$
maps by a weak equivalence to hermitian $K$-theory mod\ $\ell^{\nu}$.
\vspace{0.1in}

The following results are concerned with Bousfield $\ell$-adic completions
(denoted by $\#$) of the self-conjugate and hermitian $K$-theory spectra.
Bousfield introduced this notion in \cite{Bousfield}. First we shall tabulate
the corresponding \'{e}tale sheaves. Let $\mathbb{Z}_{\ell}^{\otimes k}$
denote the $k\,$th Tate twist of the $\ell$-adic integers. \vspace{0.1in}

In the example of self-conjugate $K$-theory the \'etale sheaves are periodic
in the following sense.

\begin{corollary}
\label{corollary:KSCstalk} The \'{e}tale sheaf associated to the presheaf
\[
U\mapsto\pi_{n}(\mathcal{KSC}(U)_{\#})
\]
of $\ell$-adically completed self-conjugate $K$-theory is given by:
\begin{table}[tbh]
\begin{center}%
\begin{tabular}
[c]{p{0.5in}|p{0.5in}|p{0.5in}|}\hline
$n\bmod4$ & $\ell=2$ & $\ell\neq2$\\\hline
$4k$ & $\mathbb{Z}_{2}^{\otimes2k}$ & $\mathbb{Z}_{\ell}^{\otimes2k}$\\
$4k+1$ & $\mu_{2}^{\otimes2k+1}$ & $0$\\
$4k+2$ & $0$ & $0$\\
$4k+3$ & $\mathbb{Z}_{2}^{\otimes2k+2}$ & $\mathbb{Z}_{\ell}^{\otimes2k+2}%
$\\\hline
\end{tabular}
\end{center}
\end{table}
\end{corollary}

For hermitian $K$-theory the \'etale sheaves are periodic in the following sense.

\begin{corollary}
\label{corollary:2KQsheaves} The \'{e}tale sheaf associated to the presheaf
\[
U\mapsto\pi_{n}({}_{\varepsilon}\mathcal{KQ}(U)_{\#})
\]
of $\ell$-adically completed hermitian $K$-theory is given as follows.

\begin{enumerate}
\item For $\ell=2$ by: \vspace{-0.1in} \begin{table}[tbh]
\begin{center}%
\begin{tabular}
[c]{p{0.5in}|p{0.5in}|p{0.5in}|}\hline
$n\bmod8$ & $\varepsilon=1$ & $\varepsilon=-1$\\\hline
$8k$ & $\mathbb{Z}_{2}^{\otimes4k}$ & $\mathbb{Z}_{2}^{\otimes4k}$\\
$8k+1$ & $\mu_{2}^{\otimes4k+1}$ & $0$\\
$8k+2$ & $\mu_{2}^{\otimes4k+1}$ & $0$\\
$8k+3$ & $0$ & $0$\\
$8k+4$ & $\mathbb{Z}_{2}^{\otimes4k+2}$ & $\mathbb{Z}_{2}^{\otimes4k+2}$\\
$8k+5$ & $0$ & $\mu_{2}^{\otimes4k+3}$\\
$8k+6$ & $0$ & $\mu_{2}^{\otimes4k+3}$\\
$8k+7$ & $0$ & $0$\\\hline
\end{tabular}
\end{center}
\end{table}

\item For $\ell\neq2$ and $\varepsilon=\pm1$ by $\mathbb{Z}_{\ell}^{\otimes
2k}$ if $n=4k$, and trivial otherwise.
\end{enumerate}
\end{corollary}

\vspace{0.1in}

In the following, \'{e}tale cohomology is continuous \'{e}tale cohomology
\cite{DF}, \cite{UJ}. \vspace{0.1in}

As a result of the previous corollaries, the descent spectral sequences for
the $2$-completed \'{e}tale self-conjugate \'{e}tale $K$-theory and hermitian
$K$-theory of $S$ take the forms
\[
{}_{1}E_{2}^{p,q}=\left\{
\begin{array}
[c]{lll}%
H_{\text{$\mathrm{\acute{e}t}$}}^{p}(S,\,\mathbb{Z}_{2}^{\otimes2k}) & q=4k &
\\
H_{\text{$\mathrm{\acute{e}t}$}}^{p}(S,\,\mu_{2}^{\otimes2k+1}) & q=4k+1 & \\
H_{\text{$\mathrm{\acute{e}t}$}}^{p}(S,\,\mathbb{Z}_{2}^{\otimes2k+2}) &
q=4k+3 & \\
0 & q=4k+2 &
\end{array}
\right\}  \Longrightarrow{KSC}_{q-p}^{\text{$\mathrm{\acute{e}t}$}}(S)_{\#},
\]%
\[
{}_{1}E_{2}^{p,q}=\left\{
\begin{array}
[c]{lll}%
H_{\text{$\mathrm{\acute{e}t}$}}^{p}(S,\,\mathbb{Z}_{2}^{\otimes4k}) & q=8k &
\\
H_{\text{$\mathrm{\acute{e}t}$}}^{p}(S,\,\mu_{2}^{\otimes4k+1}) & q=8k+1 & \\
H_{\text{$\mathrm{\acute{e}t}$}}^{p}(S,\,\mu_{2}^{\otimes4k+1}) & q=8k+2 & \\
H_{\text{$\mathrm{\acute{e}t}$}}^{p}(S,\,\mathbb{Z}_{2}^{\otimes4k+2}) &
q=8k+4 & \\
0 & q\equiv3,5,6,7\;(\mathrm{mod}\ 8) &
\end{array}
\right\}  \Longrightarrow{}_{1}{KQ}_{q-p}^{\text{$\mathrm{\acute{e}t}$}%
}(S)_{\#},
\]
and
\[
{}_{-1}E_{2}^{p,q}=\left\{
\begin{array}
[c]{lll}%
H_{\text{$\mathrm{\acute{e}t}$}}^{p}(S,\,\mathbb{Z}_{2}^{\otimes4k}) & q=8k &
\\
H_{\text{$\mathrm{\acute{e}t}$}}^{p}(S,\,\mathbb{Z}_{2}^{\otimes4k+2}) &
q=8k+4 & \\
H_{\text{$\mathrm{\acute{e}t}$}}^{p}(S,\,\mu_{2}^{\otimes4k+3}) & q=8k+5 & \\
H_{\text{$\mathrm{\acute{e}t}$}}^{p}(S,\,\mu_{2}^{\otimes4k+3}) & q=8k+6 & \\
0 & q\equiv1,2,3,7\;(\mathrm{mod}\ 8) &
\end{array}
\right\}  \Longrightarrow{}_{-1}{KQ}_{q-p}^{\text{$\mathrm{\acute{e}t}$}%
}(S)_{\#}.
\]
\vspace{0.1in}

For $\ell\neq2$ the descent spectral sequences take the forms
\[
{}_{1}E_{2}^{p,q}=\left\{
\begin{array}
[c]{lll}%
H_{\text{$\mathrm{\acute{e}t}$}}^{p}(S,\,\mathbb{Z}_{\ell}^{\otimes2k}) &
q=4k & \\
H_{\text{$\mathrm{\acute{e}t}$}}^{p}(S,\,\mathbb{Z}_{\ell}^{\otimes2k+2}) &
q=4k+3 & \\
0 & q\equiv1,2\;(\mathrm{mod}\ 4) &
\end{array}
\right\}  \Longrightarrow{KSC}_{q-p}^{\text{$\mathrm{\acute{e}t}$}}(S)_{\#},
\]
and
\[
{}_{\varepsilon}E_{2}^{p,q}=\left\{
\begin{array}
[c]{lll}%
H_{\text{$\mathrm{\acute{e}t}$}}^{p}(S,\,\mathbb{Z}_{\ell}^{\otimes\frac{q}%
{2}}) & q\equiv0\;(\mathrm{mod}\ 4) & \\
0 & q\not \equiv 0\;(\mathrm{mod}\ 4) &
\end{array}
\right\}  \Longrightarrow{}_{\varepsilon}{KQ}_{q-p}^{\text{$\mathrm{\acute
{e}t}$}}(S)_{\#}.
\]
\vspace{0.1in}

Our next objective is to compute $\ell$-adically completed \'{e}tale
self-conjugate and hermitian $K$-groups in terms of \'{e}tale cohomology
groups. To this end we need some more notation. \vspace{0.1in}

Let $A\bullet B$ denote an abelian group extension of $B$ by $A$, so that
there exists a short exact sequence
\[
0\rightarrow A\rightarrow A\bullet B\rightarrow B\rightarrow0.
\]

\begin{lemma}
\label{lemma:2cd2computation} If $\mathrm{cd}_{2}(S)=2$ and $H_{\text{\'{e}t}%
}^{0}(S,\mathbb{Z}_{2}^{\otimes i})=0$ for $i>0$ then the $2$-completed
\'{e}tale hermitian $K$-groups of $S$ are computed up to extensions in the
following table. \begin{table}[tbh]
\begin{center}%
\begin{tabular}
[c]{p{0.5in}|p{2.0in}|p{2.0in}|}\hline
$n\bmod8$ & ${}_{1}{KQ}_{n}^{\mathrm{\acute{e}t}}(S)_{\#}$ & ${}_{-1}{KQ}%
_{n}^{\mathrm{\acute{e}t}}(S)_{\#}$\\\hline
$8k>0$ & $H_{\mathrm{\acute{e}t}}^{2}(S,\mu_{2}^{\otimes4k+1})\bullet
H_{\mathrm{\acute{e}t}}^{1}(S,\mu_{2}^{\otimes4k+1})$ & $0$\\
$8k+1$ & $H_{\mathrm{\acute{e}t}}^{1}(S,\mu_{2}^{\otimes4k+1})\bullet
H_{\mathrm{\acute{e}t}}^{0}(S,\mu_{2}^{\otimes4k+1})$ & $0$\\
$8k+2$ & $H_{\mathrm{\acute{e}t}}^{2}(S,\mathbb{Z}_{2}^{\otimes4k+2})\bullet
H_{\mathrm{\acute{e}t}}^{0}(S,\mu_{2}^{\otimes4k+1})$ & $H_{\mathrm{\acute
{e}t}}^{2}(S,\mathbb{Z}_{2}^{\otimes4k+2})$\\
$8k+3$ & $H_{\mathrm{\acute{e}t}}^{1}(S,\mathbb{Z}_{2}^{\otimes4k+2})$ &
$H_{\mathrm{\acute{e}t}}^{2}(S,\mu_{2}^{\otimes4k+3})\bullet H_{\mathrm{\acute
{e}t}}^{1}(S,\mathbb{Z}_{2}^{\otimes4k+2})$\\
$8k+4$ & $0$ & $H_{\mathrm{\acute{e}t}}^{2}(S,\mu_{2}^{\otimes4k+3})\bullet
H_{\mathrm{\acute{e}t}}^{1}(S,\mu_{2}^{\otimes4k+3})$\\
$8k+5$ & $0$ & $H_{\mathrm{\acute{e}t}}^{1}(S,\mu_{2}^{\otimes4k+3})\bullet
H_{\mathrm{\acute{e}t}}^{0}(S,\mu_{2}^{\otimes4k+3})$\\
$8k+6$ & $H_{\mathrm{\acute{e}t}}^{2}(S,\mathbb{Z}_{2}^{\otimes4k+4})$ &
$H_{\mathrm{\acute{e}t}}^{2}(S,\mathbb{Z}_{2}^{\otimes4k+4})\bullet
H^{0}_{\mathrm{\acute{e}t}}(S,\mu_{2}^{\otimes4k+3})$\\
$8k+7$ & $H_{\mathrm{\acute{e}t}}^{2}(S,\mu_{2}^{\otimes4k+5})\bullet
H^{1}_{\mathrm{\acute{e}t}}(S,\mathbb{Z}_{2}^{\otimes4k+4})$ & $H^{1}%
_{\mathrm{\acute{e}t}}(S,\mathbb{Z}_{2}^{\otimes4k+4})$\\\hline
\end{tabular}
\end{center}
\end{table}
\end{lemma}

\begin{remark}
The assumption on the vanishing of $H_{\text{\'{e}t}}^{0}(S,\mathbb{Z}%
_{2}^{\otimes i})$ for $i>0$ in Lemma \ref{lemma:2cd2computation} is a
commonplace and holds for the examples considered in Section
\ref{section:applications}. We note, however, that the assumptions in Lemma
\ref{lemma:2cd2computation} are not satisfied for the field of real numbers,
and for number fields with at least one real embedding.
\end{remark}

\begin{lemma}
\label{secondoddprimarycomputation} If $\ell$ is an odd prime and
$\mathrm{cd}_{\ell}(S)\leq7$ the $\ell$-completed \'{e}tale hermitian
$K$-groups of $S$ are computed up to extensions in the following table.
\begin{table}[tbh]
\begin{center}%
\begin{tabular}
[c]{p{0.5in}|p{2.0in}|}\hline
$n\bmod4$ & ${}_{\varepsilon}{KQ}_{n}^{\mathrm{\acute{e}t}}(S)_{\#}$\\\hline
$4k>0$ & $H^{4}_{\mathrm{\acute{e}t}}(S,\mathbb{Z}_{\ell}^{\otimes
2k+2})\bullet H^{0}_{\mathrm{\acute{e}t}}(S,\mathbb{Z}_{\ell}^{\otimes2k})$\\
$4k+1$ & $H^{7}_{\mathrm{\acute{e}t}}(S,\mathbb{Z}_{\ell}^{\otimes
2k+4})\bullet H^{3}_{\mathrm{\acute{e}t}}(S,\mathbb{Z}_{\ell}^{\otimes2k+2}%
)$\\
$4k+2$ & $H^{6}_{\mathrm{\acute{e}t}}(S,\mathbb{Z}_{\ell}^{\otimes
2k+4})\bullet H^{2}_{\mathrm{\acute{e}t}}(S,\mathbb{Z}_{\ell}^{\otimes2k+2}%
)$\\
$4k+3$ & $H^{5}_{\mathrm{\acute{e}t}}(S,\mathbb{Z}_{\ell}^{\otimes
2k+4})\bullet H^{1}_{\mathrm{\acute{e}t}}(S,\mathbb{Z}_{\ell}^{\otimes2k+2}%
)$\\\hline
\end{tabular}
\end{center}
\end{table}
\end{lemma}

Corollary \ref{corollary:2KQsheaves} and the corresponding result for
algebraic $K$-theory imply the next result by inspection.

\begin{corollary}
The \'{e}tale sheaf associated to the presheaf
\[
U\mapsto\pi_{n}({}_{\varepsilon}\mathcal{V}(U)_{\#})
\]
of $\ell$-adically completed hermitian $V$-theory is given as follows.

\begin{enumerate}
\item For $\ell=2$ by: \begin{table}[tbh]
\begin{center}%
\begin{tabular}
[c]{p{0.5in}|p{0.5in}|p{0.5in}|}\hline
$n\bmod8$ & $\varepsilon=1$ & $\varepsilon=-1$\\\hline
$8k$ & $0$ & $0$\\
$8k+1$ & $\mathbb{Z}_{2}^{\otimes4k+1}$ & $\mathbb{Z}_{2}^{\otimes4k+1}$\\
$8k+2$ & $\mu_{2}^{\otimes4k+1}$ & $0$\\
$8k+3$ & $\mu_{2}^{\otimes4k+2}$ & $0$\\
$8k+4$ & $0$ & $0$\\
$8k+5$ & $\mathbb{Z}_{2}^{\otimes4k+3}$ & $\mathbb{Z}_{2}^{\otimes4k+3}$\\
$8k+6$ & $0$ & $\mu_{2}^{\otimes4k+3}$\\
$8k+7$ & $0$ & $\mu_{2}^{\otimes4k+4}$\\\hline
\end{tabular}
\end{center}
\end{table}

\item For $\ell\neq2$ by $\mathbb{Z}_{\ell}^{\otimes2k+1}$ if $n=4k+1$, and
trivial otherwise.
\end{enumerate}
\end{corollary}

\vspace{0.1in}

The previous corollary allows us to immediately identify the $E_{2}$-page of
the descent spectral sequence for $\ell$-adically completed \'{e}tale
hermitian $V$-groups in terms of \'{e}tale cohomology. As a consequence,
imposing a commonplace restriction on the \'{e}tale cohomological dimension
yields the following computation.

\begin{lemma}
\label{lemma:completedcd2lemma} If $\mathrm{cd}_{2}(S)=2$ and
$H_{\text{\'{e}t}}^{0}(S,\,\mathbb{Z}_{2}^{\otimes i})=0$ for $i>0$, then the
$2$-completed \'{e}tale $V$-groups of $S$ are computed up to extensions in the
following table. \begin{table}[tbh]
\begin{center}%
\begin{tabular}
[c]{p{0.5in}|p{2.0in}|p{2.0in}|}\hline
$n\bmod8$ & ${}_{1}{V}_{n}^{\mathrm{\acute{e}t}}(S)_{\#}$ & ${}_{-1}{V}%
_{n}^{\mathrm{\acute{e}t}}(S)_{\#}$\\\hline
$8k\geq0$ & $H_{\mathrm{\acute{e}t}}^{2}(S,\mu_{2}^{\otimes4k+1})\bullet
H_{\mathrm{\acute{e}t}}^{1}(S,\mathbb{Z}_{2}^{\otimes4k+1})$ &
$H_{\mathrm{\acute{e}t}}^{1}(S,\mathbb{Z}_{2}^{\otimes4k+1})$\\
$8k+1$ & $H_{\mathrm{\acute{e}t}}^{2}(S,\mu_{2}^{\otimes4k+2})\bullet
H_{\mathrm{\acute{e}t}}^{1}(S,\mu_{2}^{\otimes4k+1})$ & $0$\\
$8k+2$ & $H_{\mathrm{\acute{e}t}}^{1}(S,\mu_{2}^{\otimes4k+2})\bullet
H_{\mathrm{\acute{e}t}}^{0}(S,\mu_{2}^{\otimes4k+1})$ & $0$\\
$8k+3$ & $H_{\mathrm{\acute{e}t}}^{2}(S,\mathbb{Z}_{2}^{\otimes4k+3})\bullet
H_{\mathrm{\acute{e}t}}^{0}(S,\mu_{2}^{\otimes4k+2})$ & $H_{\mathrm{\acute
{e}t}}^{2}(S,\mathbb{Z}_{2}^{\otimes4k+3})$\\
$8k+4$ & $H_{\mathrm{\acute{e}t}}^{1}(S,\mathbb{Z}_{2}^{\otimes4k+3})$ &
$H_{\mathrm{\acute{e}t}}^{2}(S,\mu_{2}^{\otimes4k+3})\bullet H_{\mathrm{\acute
{e}t}}^{1}(S,\mathbb{Z}_{2}^{\otimes4k+3})$\\
$8k+5$ & $0$ & $H_{\mathrm{\acute{e}t}}^{2}(S,\mu_{2}^{\otimes4k+4})\bullet
H_{\mathrm{\acute{e}t}}^{1}(S,\mu_{2}^{\otimes4k+3})$\\
$8k+6$ & $0$ & $H_{\mathrm{\acute{e}t}}^{1}(S,\mu_{2}^{\otimes4k+4})\bullet
H^{0}_{\mathrm{\acute{e}t}}(S,\mu_{2}^{\otimes4k+3})$\\
$8k+7$ & $H_{\mathrm{\acute{e}t}}^{2}(S,\mathbb{Z}_{2}^{\otimes4k+5})$ &
$H_{\mathrm{\acute{e}t}}^{2}(S,\mathbb{Z}_{2}^{\otimes4k+5})\bullet
H^{0}_{\mathrm{\acute{e}t}}(S,\mu_{2}^{\otimes4k+4})$\\\hline
\end{tabular}
\end{center}
\end{table}
\end{lemma}

The previous computations are supplemented by more specialized examples in the
next section. For the earliest \'{e}tale $K$-theory computations we refer the
reader to \cite{DF} and \cite{Thomason}. It is worthwhile to point out that
the difference between the \'{e}tale $K$-theory, in \textsl{loc.~cit.}, and
the \'{e}tale hermitian $K$-theory computations in this paper is reminiscent
of the situation for the classical Atiyah-Hirzebruch spectral sequences based
on complex and real topological $K$-theory. This analogy is evident on the
level of \'{e}tale stalks by comparison with the complex and real $K$-theories
of a point.

\section{Applications to finite fields, local and global fields}

\label{section:applications} In this section we point out some computational
consequences of the above results. The examples are geometric in nature and
relate to finite fields, and to local and global number fields. Our main
interest and focus are on the $2$-primary computations.

\begin{ex}
In what follows, we apply Lemma \ref{secondoddprimarycomputation} to some
classes of examples. Throughout, $\ell$ is an odd prime number.

\begin{enumerate}
\item If $S$ is a $d$-dimensional smooth complex variety then $\mathrm{cd}%
_{\ell}(S)\leq2d$. Lemma \ref{secondoddprimarycomputation} computes the group
${}_{\varepsilon}KQ_{n}^{\text{\'{e}t}}(S)_{\#}$ up to extensions if $S$ is of
dimension at most $3$. For curves and surfaces there are no undetermined
extensions. Detailed computations of the algebraic $K$-theory of $S$ were
worked out in \cite{PW1} and \cite{PW2}.

\item If $S$ is a smooth curve over a number field then $\mathrm{cd}_{\ell
}(S)=4$. In this case, for $n>0$, there are no undetermined extensions in the
computation of ${}_{\varepsilon}KQ_{n}^{\text{\'{e}t}}(S)_{\#}$ given in Lemma
\ref{secondoddprimarycomputation}. For detailed computations of the algebraic
$K$-theory of $S$ we refer to \cite{ROcurves}.

\item The ring $\mathcal{O}_{F}\left[  1/\ell\right]  $ of $\ell$-integers in
any number field $F$ has cohomological dimension $\mathrm{cd}_{\ell
}(\mathcal{O}_{F}\left[  1/\ell\right]  )=\mathrm{cd}_{\ell}(F)=2$. In
particular, ${}_{\varepsilon}KQ_{n}^{\text{\'{e}t}}(\mathcal{O}_{F}\left[
1/\ell\right]  )_{\#}$ is trivial when $0<n\equiv0,1\;(\mathrm{mod}$\ $4)$ and
finite when $n\equiv2\;(\mathrm{mod}$\ $4)$. The same cohomological dimension
bound holds for local number fields and their valuation rings, \textsl{e.g.}%
~the field of $\ell$-adic numbers.
\end{enumerate}
\end{ex}

Let $F$ be a field of characteristic $\neq2$ and $\zeta_{r}$ be a primitive
$r$\thinspace th root of unity. For $i\in\mathbb{Z}$, let $w_{i}(F)$ be the
maximal $2$-power $2^{n}$ such that the exponent of the Galois group of
$F(\zeta_{2^{n}})/F$ divides $i$. If $F$ contains $\zeta_{4}$ and
$i=2^{\lambda}k$ with $k$ odd, then $w_{i}(F)=2^{r+\lambda}$ where $r$ is
maximal such that $F$ contains a primitive $2^{r}$-root of unity. If $i$ is
odd, $w_{i}(\mathbb{Q}(\sqrt{-1}))=4$, while if $\sqrt{-1}\not \in F$ then
$w_{i}(F)=2$. In all our examples the number $w_{i}(F)$ is finite.
\vspace{0.1in}

Using Lemma \ref{lemma:2cd2computation} and the \'{e}tale cohomology groups of
finite fields, we tabulate the $2$-completed \'{e}tale hermitian $K$-groups of
$\mathbb{F}_{t}$ for $t$ odd. Our findings are in agreement with Friedlander's
computation of the hermitian $K$-groups of $\mathbb{F}_{t}$ in
\cite{Friedlander}.

\begin{ex}
Let $\mathbb{F}_{t}$ be a finite field with an odd number of elements $t$. The
$2$-completed \'{e}tale hermitian $K$-groups of $\mathbb{F}_{t}$ are computed
in the following table. \begin{table}[tbh]
\begin{center}%
\begin{tabular}
[c]{p{0.5in}|p{0.9in}|p{0.9in}|}\hline
$n\bmod8$ & ${}_{1}{KQ}_{n}^{\mathrm{\acute{e}t}}(\mathbb{F}_{t})_{\#}$ &
${}_{-1}{KQ}_{n}^{\mathrm{\acute{e}t}}(\mathbb{F}_{t})_{\#}$\\\hline
$8k>0$ & $\mathbb{Z}/2$ & $0$\\
$8k+1$ & $(\mathbb{Z}/2)^{2}$ & $0$\\
$8k+2$ & $\mathbb{Z}/2$ & $0$\\
$8k+3$ & $\mathbb{Z}/w_{4k+2}(\mathbb{F}_{t})$ & $\mathbb{Z}/w_{4k+2}%
(\mathbb{F}_{t})$\\
$8k+4$ & $0$ & $\mathbb{Z}/2$\\
$8k+5$ & $0$ & $(\mathbb{Z}/2)^{2}$\\
$8k+6$ & $0$ & $\mathbb{Z}/2$\\
$8k+7$ & $\mathbb{Z}/w_{4k+4}(\mathbb{F}_{t})$ & $\mathbb{Z}/w_{4k+4}%
(\mathbb{F}_{t})$\\\hline
\end{tabular}
\end{center}
\end{table}

The extension problem for ${}_{1}KQ_{8k+1}^{\text{\'{e}t}}(\mathbb{F}%
_{t})_{\#}$ can be resolved using the homotopy fibration
\cite{KaroubiAnnals2}
\[
\mathcal{KSC}(\mathbb{F}_{t})\longrightarrow\Omega{}_{\varepsilon}%
\mathcal{KQ}(\mathbb{F}_{t})\overset{\sigma^{(2)}}{\longrightarrow}\Omega
^{-1}{}_{-\varepsilon}\mathcal{KQ}(\mathbb{F}_{t}).
\]
The group $KSC_{8k}(\mathbb{F}_{t})_{\#}$ has order $2$, \textsl{cf.}~Example
\ref{ex:KSCfinitefields}, and is a direct summand of ${}_{1}KQ_{8k+1}%
(\mathbb{F}_{t})_{\#}$. This also resolves the extension problem in degree
$8k+5$ for $\varepsilon=-1$.
\end{ex}

\begin{ex}
\label{ex:dyadicKQ} Lemma \ref{lemma:2cd2computation} applies to every dyadic
local field $F$, \textsl{i.e.}~every finite extension of the $2$-adic numbers
$\mathbb{Q}_{2}$. If the field extension degree $[F:\mathbb{Q}_{2}]=d$, then
$H_{\text{\'{e}t}}^{0}(F,\,\mathbb{Z}_{2}^{\otimes i})=\delta_{i0}%
\mathbb{Z}_{2}$,\ $H_{\text{\'{e}t}}^{1}(F,\,\mathbb{Z}_{2}^{\otimes
1})=\mathbb{Z}_{2}^{d+1}\oplus\mathbb{Z}/2$,\ $H_{\text{\'{e}t}}%
^{2}(F,\,\mathbb{Z}_{2}^{\otimes1})=\mathbb{Z}_{2}$,
\[
H_{\text{\'{e}t}}^{1}(F,\,\mathbb{Z}_{2}^{\otimes i})=\left\{
\begin{array}
[c]{ll}%
\mathbb{Z}_{2}^{d}\oplus\mathbb{Z}/2 & {i>1\text{ odd,}}\\
\mathbb{Z}_{2}^{d}\oplus\mathbb{Z}/w_{i}(F) & {i\text{ even,}}%
\end{array}
\right.
\]
and
\[
H_{\text{\'{e}t}}^{2}(F,\,\mathbb{Z}_{2}^{\otimes i})=\left\{
\begin{array}
[c]{ll}%
\mathbb{Z}/w_{i-1}(F) & {i>1\text{ odd,}}\\
\mathbb{Z}/2 & {i\text{ even.}}%
\end{array}
\right.
\]

The $2$-completed \'{e}tale hermitian $K$-groups of $F$ are computed up to
extensions in the following table. \vspace{-0.04in} \begin{table}[tbh]
\begin{center}%
\begin{tabular}
[c]{p{0.5in}|p{1.5in}|p{1.5in}|}\hline
$n\bmod8$ & ${}_{1}{KQ}_{n}^{\mathrm{\acute{e}t}}(F)_{\#}$ & ${}_{-1}{KQ}%
_{n}^{\mathrm{\acute{e}t}}(F)_{\#}$\\\hline
$8k>0$ & $\mathbb{Z}/2\bullet(\mathbb{Z}/2)^{d+2}$ & $0$\\
$8k+1$ & $(\mathbb{Z}/2)^{d+2}\bullet\mathbb{Z}/2$ & $0$\\
$8k+2$ & $\mathbb{Z}/2\bullet\mathbb{Z}/2$ & $\mathbb{Z}/2$\\
$8k+3$ & $\mathbb{Z}_{2}^{d}\oplus\mathbb{Z}/w_{4k+2}(F)$ & $\mathbb{Z}%
/2\bullet(\mathbb{Z}_{2}^{d}\oplus\mathbb{Z}/w_{4k+2}(F))$\\
$8k+4$ & $0$ & $\mathbb{Z}/2\bullet(\mathbb{Z}/2)^{d+2}$\\
$8k+5$ & $0$ & $(\mathbb{Z}/2)^{d+2}\bullet\mathbb{Z}/2$\\
$8k+6$ & $\mathbb{Z}/2$ & $\mathbb{Z}/2\bullet\mathbb{Z}/2$\\
$8k+7$ & $\mathbb{Z}/2\bullet(\mathbb{Z}_{2}^{d}\oplus\mathbb{Z}/w_{4k+4}(F))$
& $\mathbb{Z}_{2}^{d}\oplus\mathbb{Z}/w_{4k+4}(F)$\\\hline
\end{tabular}
\end{center}
\end{table}
\end{ex}

For a non-dyadic local field, \textsl{i.e.}~a finite extension of the $p$-adic
numbers $\mathbb{Q}_{p}$ for some odd prime $p$, the $2$-completed \'{e}tale
hermitian $K$-groups are comprised of finite groups in positive degrees. The
\'{e}tale cohomology computation leading to this conclusion is given in
\cite[Proposition 7.3.10]{NSW}.

\begin{ex}
\label{ex:nondyadicKQ} The $2$-completed \'{e}tale hermitian $K$-groups of a
finite extension $F$ of $\mathbb{Q}_{p}$ for $p$ odd are computed up to
extensions in the following table. \begin{table}[tbh]
\begin{center}%
\begin{tabular}
[c]{p{0.5in}|p{1.2in}|p{1.2in}|}\hline
$n\bmod8$ & ${}_{1}{KQ}_{n}^{\mathrm{\acute{e}t}}(F)_{\#}$ & ${}_{-1}{KQ}%
_{n}^{\mathrm{\acute{e}t}}(F)_{\#}$\\\hline
$8k>0$ & $\mathbb{Z}/2\bullet(\mathbb{Z}/2)^{2}$ & $0$\\
$8k+1$ & $(\mathbb{Z}/2)^{2}\bullet\mathbb{Z}/2$ & $0$\\
$8k+2$ & $\mathbb{Z}/2\bullet\mathbb{Z}/2$ & $\mathbb{Z}/2$\\
$8k+3$ & $\mathbb{Z}/w_{4k+2}(F)$ & $\mathbb{Z}/2\bullet\mathbb{Z}%
/w_{4k+2}(F)$\\
$8k+4$ & $0$ & $\mathbb{Z}/2\bullet(\mathbb{Z}/2)^{2}$\\
$8k+5$ & $0$ & $(\mathbb{Z}/2)^{2}\bullet\mathbb{Z}/2$\\
$8k+6$ & $\mathbb{Z}/2$ & $\mathbb{Z}/2\bullet\mathbb{Z}/2$\\
$8k+7$ & $\mathbb{Z}/2\bullet\mathbb{Z}/w_{4k+4}(F)$ & $\mathbb{Z}%
/w_{4k+4}(F)$\\\hline
\end{tabular}
\end{center}
\end{table}
\end{ex}

A totally imaginary number field $F$ is called $2$\emph{-regular} if the
$2$-Sylow subgroup of $K_{2}(\mathcal{O}_{F})$ is trivial. The Gaussian
numbers $\mathbb{Q}(\sqrt{-1})$ is an example of such a number field. For the
\'{e}tale cohomology of $\mathcal{O}_{F}\left[  1/2\right]  $ the $2$-regular
assumption implies that $H_{\text{\'{e}t}}^{2}(\mathcal{O}_{F}\left[
1/2\right]  ,\,\mathbb{Z}_{2}^{\otimes i})$ is trivial for $i\neq0,1$
\cite[Proposition 2.2]{RO2regular}. Moreover, $H_{\text{\'{e}t}}%
^{1}(\mathcal{O}_{F}\left[  1/2\right]  ,\,\mathbb{Z}_{2}^{\otimes i})$
identifies with $\mathbb{Z}_{2}^{c}\oplus\mathbb{Z}/w_{i}(F)$ for $i\neq0$,
and $H_{\text{\'{e}t}}^{1}(\mathcal{O}_{F}\left[  1/2\right]  ,\,\mu
_{2}^{\otimes i})\cong(\mathbb{Z}/2)^{c+1}$ where $c$ denotes the number of
pairs of complex embeddings of the number field $F$. By way of example, the
number $w_{i}(\mathbb{Q}(\sqrt{-1}))=2^{2+(i)_{2}}$ for all $i$ , where
$(i)_{2}$ is the $2$-adic valuation of $i$. \vspace{0.1in}

With these preliminaries in hand, we are ready to state the following computation.

\begin{ex}
\label{ex:totallyimaginary2regular} Let $F$ be a totally imaginary $2$-regular
number field with $c$ pairs of complex embeddings. The $2$-completed \'{e}tale
hermitian $K$-groups of its ring of $2$-integers $\mathcal{O}_{F}\left[
1/2\right]  $ are computed up to extensions in the following table.
\begin{table}[tbh]
\begin{center}%
\begin{tabular}
[c]{p{0.5in}|p{1.2in}|p{1.2in}|}\hline
$n\bmod8$ & ${}_{1}{KQ}_{n}^{\mathrm{\acute{e}t}}(\mathcal{O}_{F}\left[
1/2\right]  )_{\#}$ & ${}_{-1}{KQ}_{n}^{\mathrm{\acute{e}t}}(\mathcal{O}%
_{F}\left[  1/2\right]  )_{\#}$\\\hline
$8k>0$ & $(\mathbb{Z}/2)^{c+1}$ & $0$\\
$8k+1$ & $(\mathbb{Z}/2)^{c+1}\bullet\mathbb{Z}/2$ & $0$\\
$8k+2$ & $\mathbb{Z}/2$ & $0$\\
$8k+3$ & $\mathbb{Z}_{2}^{c}\oplus\mathbb{Z}/w_{4k+2}(F)$ & $\mathbb{Z}%
_{2}^{c}\oplus\mathbb{Z}/w_{4k+2}(F)$\\
$8k+4$ & $0$ & $(\mathbb{Z}/2)^{c+1}$\\
$8k+5$ & $0$ & $(\mathbb{Z}/2)^{c+1}\bullet\mathbb{Z}/2$\\
$8k+6$ & $0$ & $\mathbb{Z}/2$\\
$8k+7$ & $\mathbb{Z}_{2}^{c}\oplus\mathbb{Z}/w_{4k+4}(F)$ & $\mathbb{Z}%
_{2}^{c}\oplus\mathbb{Z}/w_{4k+4}(F)$\\\hline
\end{tabular}
\end{center}
\end{table}
\end{ex}

\begin{remark}
We expect that ${}_{1}KQ_{8k+1}^{\text{\'{e}t}}(\mathcal{O}_{F}\left[
1/2\right]  )_{\#}$ and ${}_{-1}KQ_{8k+5}^{\text{\'{e}t}}(\mathcal{O}%
_{F}\left[  1/2\right]  )_{\#}$ are elementary abelian $2$-groups of rank
equal to $c+2$.
\end{remark}

In the following discussion of \'{e}tale $V$-theory we shall specialize Lemma
\ref{lemma:completedcd2lemma} to the previous examples of finite fields, local
fields and totally imaginary $2$-regular number fields. As for hermitian
\'etale $K$-theory, it turns out that the \'{e}tale $V$-groups of dyadic and
non-dyadic local number fields are completely different; although in some
degrees we are only able to compute these groups up to extensions, we can
conclude that the former allow free summands in some degrees while the latter
are always finite abelian groups.

\vspace{0.1in}

Our computation of the $2$-completed \'{e}tale $V$-groups of finite fields is
in agreement with Hiller's results for $V$-groups in \cite{Hiller}.

\begin{ex}
Let $\mathbb{F}_{t}$ be a finite field with an odd number of elements $t$. The
$2$-completed \'{e}tale $V$-groups of $\mathbb{F}_{t}$ are computed in the
following table. \begin{table}[tbh]
\begin{center}%
\begin{tabular}
[c]{p{0.5in}|p{0.8in}|p{0.8in}|}\hline
$n\bmod8$ & ${}_{1}{V}_{n}^{\mathrm{\acute{e}t}}(\mathbb{F}_{t})_{\#}$ &
${}_{-1}{V}_{n}^{\mathrm{\acute{e}t}}(\mathbb{F}_{t})_{\#}$\\\hline
$8k\geq0$ & $\mathbb{Z}/w_{4k+1}(\mathbb{F}_{t})$ & $\mathbb{Z}/w_{4k+1}%
(\mathbb{F}_{t})$\\
$8k+1$ & $\mathbb{Z}/2$ & $0$\\
$8k+2$ & $(\mathbb{Z}/2)^{2}$ & $0$\\
$8k+3$ & $\mathbb{Z}/2$ & $0$\\
$8k+4$ & $\mathbb{Z}/w_{4k+3}(\mathbb{F}_{t})$ & $\mathbb{Z}/w_{4k+3}%
(\mathbb{F}_{t})$\\
$8k+5$ & $0$ & $\mathbb{Z}/2$\\
$8k+6$ & $0$ & $(\mathbb{Z}/2)^{2}$\\
$8k+7$ & $0$ & $\mathbb{Z}/2$\\\hline
\end{tabular}
\end{center}
\end{table}

The extension problem in degree $8k+2$ can be resolved using that ${}%
_{1}KQ_{8k+2}(\mathbb{F}_{t})_{\#}$ has order $2$ and is a direct summand of
${}_{1}V_{8k+2}(\mathbb{F}_{t})_{\#}$. Likewise, this also resolves the
extension problem in degree $8k+6$ for $\varepsilon=-1$.
\end{ex}

Next we turn to local number fields. We find it convenient to distinguish
between dyadic and non-dyadic local fields.

\begin{ex}
The $2$-completed \'{e}tale $V$-groups of a dyadic local number field $F$ of
degree $d$ are computed up to extensions in the following table.
\begin{table}[tbh]
\begin{center}%
\begin{tabular}
[c]{p{0.5in}|p{1.2in}|p{1.2in}|}\hline
$n\bmod8$ & ${}_{1}{V}_{n}^{\mathrm{\acute{e}t}}(F)_{\#}$ & ${}_{-1}{V}%
_{n}^{\mathrm{\acute{e}t}}(F)_{\#}$\\\hline
$8k\geq0$ & $\mathbb{Z}/2\bullet(\mathbb{Z}_{2}^{d}\oplus\mathbb{Z}/2)$ &
$\mathbb{Z}_{2}^{d}\oplus\mathbb{Z}/2$\\
$8k+1$ & $\mathbb{Z}/2\bullet(\mathbb{Z}/2)^{d+2}$ & $0$\\
$8k+2$ & $(\mathbb{Z}/2)^{d+2}\bullet\mathbb{Z}/2$ & $0$\\
$8k+3$ & $\mathbb{Z}/w_{4k+2}(F)\bullet\mathbb{Z}/2$ & $\mathbb{Z}%
/w_{4k+2}(F)$\\
$8k+4$ & $\mathbb{Z}_{2}^{d}\oplus\mathbb{Z}/2$ & $\mathbb{Z}/2\bullet
(\mathbb{Z}_{2}^{d}\oplus\mathbb{Z}/2)$\\
$8k+5$ & $0$ & $\mathbb{Z}/2\bullet(\mathbb{Z}/2)^{d+2}$\\
$8k+6$ & $0$ & $(\mathbb{Z}/2)^{d+2}\bullet\mathbb{Z}/2$\\
$8k+7$ & $\mathbb{Z}/w_{4k+4}(F)$ & $\mathbb{Z}/w_{4k+4}(F)\bullet
\mathbb{Z}/2$\\\hline
\end{tabular}
\end{center}
\end{table}

If $i$ is even, the number $w_{i}(\mathbb{Q}_{2})=2^{2+(i)_{2}}$.
\end{ex}

For non-dyadic local number fields the \'{e}tale $V$-groups turn out to be
torsion abelian groups.

\begin{ex}
The $2$-completed \'{e}tale $V$-groups of a non-dyadic local number field $F$
are computed up to extensions in the following table. \begin{table}[tbh]
\begin{center}%
\begin{tabular}
[c]{p{0.5in}|p{1.2in}|p{1.2in}|}\hline
$n\bmod8$ & ${}_{1}{V}_{n}^{\mathrm{\acute{e}t}}(F)_{\#}$ & ${}_{-1}{V}%
_{n}^{\mathrm{\acute{e}t}}(F)_{\#}$\\\hline
$8k\geq0$ & $\mathbb{Z}/2\bullet\mathbb{Z}/2$ & $\mathbb{Z}/2$\\
$8k+1$ & $\mathbb{Z}/2\bullet(\mathbb{Z}/2)^{2}$ & $0$\\
$8k+2$ & $(\mathbb{Z}/2)^{2}\bullet\mathbb{Z}/2$ & $0$\\
$8k+3$ & $\mathbb{Z}/w_{4k+2}(F)\bullet\mathbb{Z}/2$ & $\mathbb{Z}%
/w_{4k+2}(F)$\\
$8k+4$ & $\mathbb{Z}/2$ & $\mathbb{Z}/2\bullet\mathbb{Z}/2$\\
$8k+5$ & $0$ & $\mathbb{Z}/2\bullet(\mathbb{Z}/2)^{2}$\\
$8k+6$ & $0$ & $(\mathbb{Z}/2)^{2}\bullet\mathbb{Z}/2$\\
$8k+7$ & $\mathbb{Z}/w_{4k+4}(F)$ & $\mathbb{Z}/w_{4k+4}(F)\bullet
\mathbb{Z}/2$\\\hline
\end{tabular}
\end{center}
\end{table}
\end{ex}

Our last example concerning \'{e}tale $V$-theory deals with totally imaginary
$2$-regular number fields. We refer to the discussion prior to Example
\ref{ex:totallyimaginary2regular} for some of the salient features of these
number fields.

\begin{ex}
Let $F$ be a totally imaginary $2$-regular number field with $c$ pairs of
complex embeddings. The $2$-completed \'{e}tale $V$-groups of its ring of
$2$-integers $\mathcal{O}_{F}\left[  1/2\right]  $ are computed up to
extensions in the following table. \begin{table}[tbh]
\begin{center}%
\begin{tabular}
[c]{p{0.5in}|p{1.2in}|p{1.2in}|}\hline
$n\bmod8$ & ${}_{1}{V}_{n}^{\mathrm{\acute{e}t}}(\mathcal{O}_{F}\left[
1/2\right]  )_{\#}$ & ${}_{-1}{V}_{n}^{\mathrm{\acute{e}t}}(\mathcal{O}%
_{F}\left[  1/2\right]  )_{\#}$\\\hline
$8k\geq0$ & $\mathbb{Z}_{2}^{c}\oplus\mathbb{Z}/w_{4k+1}(F)$ & $\mathbb{Z}%
_{2}^{c}\oplus\mathbb{Z}/w_{4k+1}(F)$\\
$8k+1$ & $(\mathbb{Z}/2)^{c+1}$ & $0$\\
$8k+2$ & $(\mathbb{Z}/2)^{c+1}\bullet\mathbb{Z}/2$ & $0$\\
$8k+3$ & $\mathbb{Z}/2$ & $0$\\
$8k+4$ & $\mathbb{Z}_{2}^{c}\oplus\mathbb{Z}/w_{4k+3}(F)$ & $\mathbb{Z}%
_{2}^{c}\oplus\mathbb{Z}/w_{4k+3}(F)$\\
$8k+5$ & $0$ & $(\mathbb{Z}/2)^{c+1}$\\
$8k+6$ & $0$ & $(\mathbb{Z}/2)^{c+1}\bullet\mathbb{Z}/2$\\
$8k+7$ & $0$ & $\mathbb{Z}/2$\\\hline
\end{tabular}
\end{center}
\end{table}
\end{ex}

\begin{remark}
We expect that ${}_{1}V_{8k+2}^{\text{\'{e}t}}(\mathcal{O}_{F}\left[
1/2\right]  )_{\#}$ and ${}_{-1}V_{8k+6}^{\text{\'{e}t}}(\mathcal{O}%
_{F}\left[  1/2\right]  )_{\#}$ are elementary abelian $2$-groups of rank
equal to $c+2$.
\end{remark}

The next examples concern self-conjugate algebraic $K$-theory.

\begin{ex}
\label{ex:KSCfinitefields} The $2$-completed $KSC$-groups of a finite field
$\mathbb{F}_{t}$ of odd characteristic are given in the following table.
\begin{table}[tbh]
\begin{center}%
\begin{tabular}
[c]{p{0.7in}|p{1.0in}|}\hline
$n\bmod4$ & ${KSC}_{n}(\mathbb{F}_{t})_{\#}$\\\hline
$4k>0$ & $\mathbb{Z}/2$\\
$4k+1$ & $\mathbb{Z}/2$\\
$4k+2$ & $\mathbb{Z}/w_{2k+2}(\mathbb{F}_{t})$\\
$4k+3$ & $\mathbb{Z}/w_{2k+2}(\mathbb{F}_{t})$\\\hline
\end{tabular}
\end{center}
\end{table}

Recall that $\tau$ denotes the duality functor in algebraic $K$-theory. The
map%
\[
\pi_{n}(1-\tau):K_{n}(\mathbb{F}_{t})_{\#}\longrightarrow{K}_{n}%
(\mathbb{F}_{t})_{\#}%
\]
is multiplication by $2$ if $n\equiv1\;($\textrm{mod}$\ 4)$ and trivial otherwise.
\end{ex}

\begin{ex}
Let $F$ be a totally imaginary $2$-regular number field with $c$ pairs of
complex embeddings. The $2$-completed $KSC$-groups of its ring of $2$-integers
$\mathcal{O}_{F}\left[  1/2\right]  $ are given in the following table.
\begin{table}[tbh]
\begin{center}%
\begin{tabular}
[c]{p{0.7in}|p{1.5in}|}\hline
$n\bmod4$ & ${KSC}_{n}(\mathcal{O}_{F}\left[  1/2\right]  )_{\#}$\\\hline
$4k>0$ & $(\mathbb{Z}/2)^{c+1}$\\
$4k+1$ & $\mathbb{Z}/2$\\
$4k+2$ & $\mathbb{Z}_{2}^{c}\oplus\mathbb{Z}/w_{2k+2}(F)$\\
$4k+3$ & $\mathbb{Z}_{2}^{c}\oplus\mathbb{Z}/w_{2k+2}(F)$\\\hline
\end{tabular}
\end{center}
\end{table}

The map
\[
\pi_{n}(1-\tau):K_n(\mathcal{O}_F\left[  1/2\right]
)_\#\longrightarrow K_n(\mathcal{O}_F\left[  1/2\right]  )_\#
\]
is multiplication by $2$ if $n\equiv1\;(\mathrm{mod}$\ $4)$ and trivial
otherwise. There is an exact sequence (with $A=\mathcal{O}_{F}\left[
1/2\right]  $)
\[
0\rightarrow KSC_{2n+1}(A)_{\#}\rightarrow K_{2n+1}(A)_{\#}\rightarrow
K_{2n+1}(A)_{\#}\rightarrow KSC_{2n}(A)_{\#}\rightarrow0.
\]

\end{ex}

In the examples above, the assertions concerning $\pi_{n}(1-\tau)$ follow by
inspection, using the computations of ${K}_{n}(\mathbb{F}_{t})_{\#}$
\cite{Quillenfinitefields} and ${K}_{n}(\mathcal{O}_{F}\left[  1/2\right]
)_{\#}$ \cite[Theorem 3.1]{RO2regular}. \vspace{0.1in}

A systematic approach to the $KSC$-computations is to first compute the
descent spectral sequence for \'{e}tale $KSC$-theory obtained from Corollary
\ref{corollary:KSCstalk}, and then invoke Theorem \ref{theorem:KSCdescent}. In
general, if $\mathrm{cd}_{2}(S)<\infty$, this approach gives \textquotedblleft
in sufficiently high degrees\textquotedblright\ a strongly convergent
cohomological spectral sequence
\[
{}_{1}E_{2}^{p,q}=\left\{
\begin{array}
[c]{ll}%
H_{\text{$\mathrm{\acute{e}t}$}}^{p}(S,\,\mathbb{Z}_{2}^{\otimes2k}) & q=4k\\
H_{\text{$\mathrm{\acute{e}t}$}}^{p}(S,\,\mu_{2}^{\otimes2k+1}) & q=4k+1\\
H_{\text{$\mathrm{\acute{e}t}$}}^{p}(S,\,\mathbb{Z}_{2}^{\otimes2k+2}) &
q=4k+3\\
0 & q=4k+2
\end{array}
\right\}  \;\;``\Longrightarrow\text{\textquotedblright}\;\;{KSC}%
_{q-p}(S)_{\#}.
\]

\section{Applications to group rings, complex varieties and commutative Banach algebras}

\label{section:applications:2}Let $G$ be a finite group and $R$ be the ring of
$S$-integers in a number field. Let $m$ be a prime power which is prime to the
order of $G$. In \cite{WeibelGrouprings}, Weibel proves a periodicity theorem
for the higher algebraic $K$-theory of the group ring $A=R\left[  G\right]  $,
with coefficients in $\mathbb{Z}/m$. More precisely, the cup-product with the
Bott element in $K$-theory induces an isomorphism%
\[
K_{i}(A;\,\mathbb{Z}/m)\cong K_{i+p}(A;\,\mathbb{Z}/m)
\]
for $i>0$. Here, the integers $m$ and $p$ are linked according to our
Convention \ref{convention: (p,p')}. We can apply our periodicity theorems of
Section 4 in order to show that, for $1/2\in A$ and any involution on $A$, for
instance the one induced by $g\mapsto g^{-1}$, we have a split short exact
sequence for $m$ a power of $2$ and $i>q-1.$
\[
0\longrightarrow\underleftarrow{\lim}{}_{\varepsilon}KQ_{i+ps}(A;\,\mathbb{Z}%
/m)\overset{\theta_{n}^{-}}{\longrightarrow}{}_{\varepsilon}KQ_{i}%
(A;\,\mathbb{Z}/m)\overset{\theta_{n}^{+}}{\longrightarrow}{}\underrightarrow
{\lim}{}_{\varepsilon}KQ_{i+ps}(A;\,\mathbb{Z}/m)\longrightarrow0\text{.}%
\]
Here the number $q$ is given by our convention \ref{stronger convention}. We
note that Weibel's theorem is also true for $i\geq0$ if we replace the number
ring $R$ by a local field.

In particular, the $KQ$-groups are also periodic, \textsl{i.e.}%
\[
{}_{\varepsilon}KQ_{i}(A;\,\mathbb{Z}/m)\cong{}_{\varepsilon}KQ_{i+p}%
(A;\,\mathbb{Z}/m)
\]
at least for $i>q-1.$ In fact, a more careful analysis forces us to
distinguish two cases according to the parity of $m$.

If $m$ is even, the condition that the order to $G$ is prime to $m$ implies
that $G$ is of odd order. According to the famous theorem of Feit and Thompson
\cite{FeitThompson}, this implies that $G$ is solvable. Therefore, for
$2$-primary coefficients, we have a periodicity statement only for a special
class of solvable groups.

If $m$ is odd, we already know, without the hypothesis $1/2\in A$, that the
group ${}_{\varepsilon}KQ_{i}(A;\,\mathbb{Z}/m)$ splits into the direct sum%
\[
{}_{\varepsilon}KQ_{i}(A;\,\mathbb{Z}/m)\cong{}_{\varepsilon}KQ_{i}%
(A;\,\mathbb{Z}/m)_{+}\oplus{}_{\varepsilon}KQ_{i}(A;\,\mathbb{Z}/m)_{-}.
\]
In this direct sum decomposition, the group ${}_{\varepsilon}KQ_{i}%
(A;\,\mathbb{Z}/m)_{-}$ is the higher Witt group with $\mathbb{Z}/m$
coefficients and we have the periodicity isomorphism%
\[
{}_{\varepsilon}KQ_{i}(A;\,\mathbb{Z}/m)_{-}\cong{}_{-e}KQ_{i+2}%
(A;\,\mathbb{Z}/m)_{-}%
\]
for all values $i\in\mathbb{Z}$. On the other hand, the group ${}%
_{\varepsilon}KQ_{i}(A;\,\mathbb{Z}/m)_{+}$ may be identified with
$K_{i}(A;\,\mathbb{Z}/m)_{+}$, the symmetric part of $K_{i}(A;\,\mathbb{Z}/m)$
with respect to the involution given by the duality. Since this involution is
compatible with the Bott map, Weibel's theorem \cite{WeibelGrouprings} implies
another periodicity isomorphism%
\[
{}_{\varepsilon}KQ_{i}(A;\,\mathbb{Z}/m)_{+}\cong{}_{\varepsilon}%
KQ_{i+p}(A;\,\mathbb{Z}/m)_{+}%
\]
but only if $i>0$. Summarizing, we have proved the
following theorem.

\begin{theorem}
Let $G$ be a finite group and let $R$ be the ring of $S$-integers in a number
field. If $m$ is a $2$-power and if $G$ is of odd order, then we have a
periodicity isomorphism%
\[
{}_{\varepsilon}KQ_{i}(R\left[  G\right]  ;\,\mathbb{Z}/m)\cong{}%
_{\varepsilon}KQ_{i+p}(R\left[  G\right]  ;\,\mathbb{Z}/m)
\]
if $1/2\in R$ and if $i>q-1$. On the other hand, if $m$ is an odd prime power
and if $G$ is an arbitrary finite group whose order is prime to $m$, we have
the same periodicity isomorphism, with only the restriction that $i>0$.
\end{theorem}

\begin{remark}
As in Section \ref{Generalization to schemes}, we may conjecture that, in the
case where $m$ is a $2$-power, the inverse limit%
\[
\underleftarrow{\lim}{}_{\varepsilon}KQ_{i+ps}(A;\,\mathbb{Z}/m)
\]
is reduced to $0$. In other words, we conjecture that the ring $A$ is
hermitian regular according to Definition \ref{defn: hermitian regular}%
. This will imply that the positive Bott map%
\[
{}_{\varepsilon}KQ_{i}(R\left[  G\right]  ;\,\mathbb{Z}/m)\longrightarrow
{}_{\varepsilon}KQ_{i+p}(R\left[  G\right]  ;\,\mathbb{Z}/m)
\]
is an isomorphism for $i>0$ according to Theorem \ref{Theorem of W-regularity}.
\end{remark}

\medskip

Let us now turn our attention to a smooth complex variety $S$ of dimension
$n$. As we briefly mentioned in Section \ref{Generalization to schemes}, the
\'{e}tale dimension of $S$ is $2n$. As a consequence of Artin-Grothendieck
theory, it is well-known that the Betti cohomology of $S$ with coefficients
$\mathbb{Z}/m$ is isomorphic to the mod\ $m$ \'{e}tale cohomology. The same
result is valid for any cohomology theory by the method initiated by Dwyer and
Friedlander \cite{DF}. For instance, the mod\ $m$ \'{e}tale $K$-theory of $S$
coincides with the mod\ $m$ complex topological $K$-theory of Atiyah and
Hirzebruch. By the same argument, the mod\ $m$ \'{e}tale {}$_{1}KQ$-theory
coincides with the mod\ $m$ $K$-theory of complex vector bundles provided with
a nondegenerate symmetric bilinear form. This theory is well understood and is
detailed for instance in Appendix B to \cite{BKO}: it is the usual mod\ $m$
topological real $K$-theory. In the same way, the mod\ $m$ \'{e}tale {}%
$_{-1}KQ$-theory coincides with the mod\ $m$ $K$-theory of complex vector
bundles provided with a nondegenerate antisymmetric bilinear form$.$ This
theory is also well understood: it is the usual mod\ $m$ topological
symplectic $K$-theory. In both cases, we shall write ${}_{\varepsilon
}\overline{KQ}_{n}^{\mathrm{top}}(S)$, with $\varepsilon=1,-1$ if we consider
symmetric or antisymmetric bilinear forms respectively.

\begin{theorem}
Let $S$ be a smooth complex variety of dimension $n$. Then the mod\ $m$ 
\'{e}tale $_{\varepsilon}KQ$-theory of $S$ coincides with the mod\ $m$
topological $K$-theory of its complex points, real or symplectic according to
$\varepsilon$. Moreover, the canonical map%
\[
{}_{\varepsilon}\overline{KQ}_{i}(S)\longrightarrow{}_{\varepsilon}%
\overline{KQ}_{i}^{\mathrm{\acute{e}t}}(S)\cong{}_{\varepsilon}\overline
{KQ}_{i}^{\mathrm{top}}(S(\mathbb{C}))
\]
is split surjective\footnote{In \cite{Schlichting2}, we show that in fact this
canonical map is an isomorphism when $i\geq 2n-1$.} when $i\geq2n+q-1$.
Moreover, for odd prime power coefficients, it is an isomorphism for
$i\geq2n-1$, with an identification of ${}_{\varepsilon}\overline{KQ}%
_{i}^{\mathrm{top}}(S)$ with $\overline{K}_{i}^{\mathrm{top}}(S)_{+}$, the
symmetric part of $\overline{K}_{i}^{\mathrm{top}}(S)$ with respect to the
involution induced by the duality functor.
\end{theorem}

\noindent\textbf{Proof.} This theorem is mostly a consequence of the general
results in Section \ref{Generalization to schemes}. What remains to be shown
is that $_{\varepsilon}\overline{KQ}_{i}(S)_{-}$ is zero for odd prime power
coefficients: this is a consequence of the fact that $-1$ has a square root in
$\mathbb{C}$. Therefore, the classical Witt group and also the higher Witt
groups have only $2$-torsion.$\hfill\Box\smallskip$

Let us now consider a real or complex commutative Banach algebra $A$. It is a
theorem of Fisher \cite{Fisher} and Prasolov \cite{Prasolov} that the natural
map%
\[
K_{i}^{\mathrm{alg}}(A;{}\mathbb{Z}/m)\longrightarrow K_{i}^{\mathrm{top}%
}(A;{}\mathbb{Z}/m)
\]
is an isomorphism for $i\geq1$. In particular, the groups $K_{i}%
^{\mathrm{alg}}(A;\mathbb{Z}/m)$ are periodic of period $2$ if $A$ is complex
and of period $8$ if $A$ is real. In this context, it is natural to state the
following conjecture.

\begin{conjecture}
Let $A$ be a real or complex commutative Banach algebra with involution. Then
the natural map
\[
_{\varepsilon}KQ_{i}^{\mathrm{alg}}(A;{}\mathbb{Z}/m)\longrightarrow
{}_{\varepsilon}KQ_{i}^{\mathrm{top}}(A;{}\mathbb{Z}/m)
\]
is an isomorphism for $i\geq1$.
\end{conjecture}

Applying the general arguments in this paper, we can prove a theorem that
would also be a consequence of this conjecture, namely the periodicity of the
groups \thinspace$_{\varepsilon}KQ_{i}^{\mathrm{alg}}(A;{}\mathbb{Z}/m)$,
which we simply write $_{\varepsilon}KQ_{i}(A;{}\mathbb{Z}/m)$. More
precisely, the theorem of Fisher and Prasolov implies that the Bott map
$K_{i}(A;{}\mathbb{Z}/m)\rightarrow K_{i+p}(A;{}\mathbb{Z}/m)$ is an
isomorphism for $i\geq1$, with $m$ and $p$ being related by our Convention
\ref{convention: (p,p')}. From Theorem \ref{theorem:splotsurjectivity}, we
therefore deduce the following periodicity pattern for the groups
{}$_{\varepsilon}KQ_{i}(A;{}\mathbb{Z}/m)$.

\begin{theorem}
Let $A$ be a real or complex commutative Banach algebra with involution. Then
we have an isomorphism
\[
{}_{\varepsilon}KQ_{i}(A;{}\mathbb{Z}/m)\cong{}_{\varepsilon}KQ_{i+p}%
(A;{}\mathbb{Z}/m)
\]
for $i\geq q$, where $m$, $p$ and $q$ are $2$-powers related by our
Conventions \ref{convention: (p,p')} and \ref{stronger convention}.$\hfill
\Box$
\end{theorem}

As a matter of fact, if $m$ is an odd prime power, we can prove a much better
result. For, we know already by our general theory that the subgroup
${}_{\varepsilon}KQ_{i}(A;{}\mathbb{Z}/m)_{-}$ is periodic of period $4$ for
all values of $i$. Moreover, ${}_{\varepsilon}KQ_{i}(A;{}\mathbb{Z}/m)_{+}$ is
isomorphic to $K_{i}(A;{}\mathbb{Z}/m)_{+}$, the symmetric part of $K$-theory
which (as a direct consequence of the theorem of Fisher and Prasolov) is
periodic of period $4$ if $A$ is complex or real. Summarizing, we get the
following more precise theorem for $m$ an odd prime power.

\begin{theorem}
Let $A$ be a real or complex commutative Banach algebra with involution and
let $m$ be an odd prime power. Then, for $i\geq1$, we have an isomorphism
given by the cup-product with a Bott element%
\[
{}_{\varepsilon}KQ_{i}(A;{}\mathbb{Z}/m)\overset{\cong}{\longrightarrow}%
{}_{\varepsilon}KQ_{i+4}(A;{}\mathbb{Z}/m)\text{.}%
\]
$\hfill\bigskip$
\end{theorem}

\vspace{0.2in}

\begin{center}
A.\thinspace Jon Berrick \\[0pt]Department of Mathematics, National University
of Singapore, Singapore. \\[0pt]e-mail: berrick@math.nus.edu.sg \vspace{0.2in}

Max Karoubi \\[0pt]UFR de Math{\'{e}}matiques, Universit{\'{e}} Paris 7,
France. \\[0pt]e-mail: max.karoubi@gmail.com \vspace{0.2in}

Paul Arne {\O }stv{\ae }r \\[0pt]Department of Mathematics, University of
Oslo, Norway. \\[0pt]e-mail: paularne@math.uio.no
\end{center}

\end{document}